\documentclass[a4paper, reqno]{amsart}      
\usepackage{amsaddr}
\usepackage{amssymb}
\usepackage[hidelinks]{hyperref}
\usepackage[maxbibnames=10, 
  url=false, 
  doi=false, 
  date=year,
  isbn=false]{biblatex}
\usepackage{enumitem}                      
\usepackage{tikz, tikz-cd}
\newcommand\addvmargin[1]{
  \node[fit=(current bounding box),inner ysep=#1,inner xsep=0]{};
}
\usepackage{algorithm, algpseudocode}

\usetikzlibrary{babel, calc, topaths, fit}
\usepackage{multicol, multirow, float, scrhack, caption, longtable}
\makeatletter
\def\blfootnote{\gdef\@thefnmark{}\@footnotetext}
\makeatother

\title{Nuclearity of Hypergraph $C^\ast$-Algebras}
\author{Björn Schäfer and Moritz Weber}
\email{bschaefer@math.uni-sb.de}
\address{Department of Mathematics, Saarland University, D-66123 Saarbr\"ucken, Germany}
\email{weber@math.uni-sb.de}
\date{\today}

\newtheorem{thm}{Theorem}[section]
\newtheorem{definition}[thm]{Definition}
\newtheorem{lemma}[thm]{Lemma}
\newtheorem{corollary}[thm]{Corollary}
\newtheorem{proposition}[thm]{Proposition}
\newtheorem*{proposition*}{Proposition}
\newtheorem*{lemma*}{Lemma}
\newtheorem*{definition*}{Definition}
\newtheorem{example}[thm]{Example}
\newtheorem{remark}[thm]{Remark}
\newtheorem{problem}[thm]{Problem}

\newcommand{\Ha}{\mathrm{H}}

\makeatletter
\newcommand{\leqnomode}{\tagsleft@true}
\newcommand{\reqnomode}{\tagsleft@false}
\makeatother

\newcommand{\C}{\mathbb{C}}

\newcommand{\N}{\mathbb{N}}
\newcommand{\Z}{\mathbb{Z}}

\tikzset{
    partial ellipse/.style args={#1:#2:#3}{
        insert path={+ (#1:#3) arc (#1:#2:#3)}
    }
}

\begin{document}

\begin{abstract}
  We partially characterize nuclearity for the recently introduced class of hypergraph $C^\ast$-algebras using a tailor-made hypergraph minor relation. The latter is generated by certain operations on hypergraphs which resemble the moves on directed graphs used by Eilers, Restorff, Ruiz and S{\o}rensen to classify unital graph $C^\ast$-algebras. In particular, we obtain a new proof of the fact that every graph $C^\ast$-algebra associated to a finite graph is nuclear.
\end{abstract}


\maketitle

\setcounter{tocdepth}{1}
\tableofcontents

\section{Introduction}

Graph $C^\ast$-algebras are an important class of $C^\ast$-algebras comprising all finite-dimensional $C^\ast$-algebras, the Toeplitz algebra as well as the Cuntz algebras and the Cuntz-Krieger algebras. For graph $C^\ast$-algebras there is a rich theory with many connections between properties of the $C^\ast$-algebra  and easily accessible properties of the graph. A good general reference for this topic is \cite{raeburnGraphAlgebras2005}.
    
Hypergraph $C^\ast$-algebras naturally generalize the concept of graph $C^\ast$-algebras by passing from directed graphs to directed hypergraphs where an edge can have multiple vertices in its range or source. These algebras have recently been introduced by Mirjam Trieb, Dean Zenner and the second author in \cite{triebHypergraphAlgebras2024}. Further, \cite{farossQuantumAutomorphismGroups2024} introduces quantum automorphism groups of hypergraphs.
Special cases of hypergraph $C^\ast$-algebras are algebras associated to ultragraphs that were introduced earlier in \cite{tomfordeUnifiedApproachExelLaca2003}. It emerged that hypergraph $C^\ast$-algebras truly extend the class of graph $C^\ast$-algebras. In particular, unlike the latter, hypergraph $C^\ast$-algebras can be non-nuclear with an example given by the unital free product $C(S^1) *_\C \C^2$ \cite[Proposition 4.2]{triebHypergraphAlgebras2024}.
Let us mention that slightly similar free product constructions are achieved by separated graph $C^\ast$-algebras as constructed in \cite{araAlgebrasSeparatedGraphs2011} (or in \cite{duncanCertainFreeProducts2010} as $C^\ast$-algebras of edge-labelled graphs). 

The literature on hypergraphs distinguishes between directed and undirected hypergraphs. For a general reference on hypergraphs we refer to \cite{bergeGraphsHypergraphs1976}, and for a survey on directed hypergraphs we refer to \cite{ausielloDirectedHypergraphsIntroduction2017}. In this work, hypergraphs are generally assumed to be directed. However, our definition of a hypergraph encompasses both directed and undirected hypergraphs (see Definition \ref{pre_hyp_def}).
    
\begin{samepage}
The present paper continues the study of hypergraph $C^\ast$-algebras and 
aims at answering the following question:
\begin{center}
    \textit{For which hypergraphs $\Ha\Gamma$ is the $C^\ast$-algebra $C^\ast(\Ha\Gamma)$ nuclear?}
\end{center}
To this end, we introduce a tailor-made hypergraph minor relation $\leq$ and obtain
%
a result of the following form (see the next Section \ref{mainres_sec}): For any hypergraph $\Ha\Gamma$ one can construct a hypergraph minor $\Ha\Delta$ of $\Ha\Gamma$ such that $C^\ast(\Ha\Gamma)$ is nuclear if, and only if, the same holds for $C^\ast(\Ha\Delta)$. Further, if the minors of $\Ha\Delta$ include one of four forbidden minors, then $C^\ast(\Ha\Delta)$ is not nuclear. 
Thus, in the latter case we obtain a negative answer to the question of nuclearity of $C^\ast(\Ha\Gamma)$. On the other hand, if $\Ha\Delta$ has none of the forbidden minors then 
it is often found that $C^\ast(\Ha\Delta)$ (and hence $C^\ast(\Ha\Gamma)$) is obviously nuclear. 
\end{samepage}

The relation $\leq$ is generated by seven (hypergraph) \emph{minor operations} which transform a given hypergraph $\Ha\Gamma$ into a ``simpler'' hypergraph $\Ha\Delta$. These operations resemble the moves for directed graphs used in \cite{eilersCompleteClassificationUnital2021} to classify unital graph $C^\ast$-algebras. However, (graph) moves aim at not changing the graph $C^\ast$-algebra in a suitable sense whereas the (hypergraph) minor operations aim at breaking down the hypergraph $C^\ast$-algebra in a controlled way. We remark that our definition of a (hypergraph) minor operation is mostly motivated by its behavior on the $C^\ast$-algebra. In particular, untypically some minor operations do increase the number of vertices or edges. For undirected hypergraphs, minors has been studied in \cite{adlerHypertreedepthMinorsHypergraphs2012}. 



In \textbf{Section \ref{mainres_sec}} we present the main results of this paper. Next, \textbf{Section \ref{min_sec}} introduces the notion of a hypergraph minor and investigates the minor operations on the $C^\ast$-algebra side. Then, in \textbf{Section \ref{norm_sec}} these results are used to obtain a normalization procedure which transforms a hypergraph without changing its $C^\ast$-algebra up to Morita equivalence. The normalization procedure is complemented with a reduction procedure in \textbf{Section \ref{red_sec}}. The latter is a procedure which simplifies a hypergraph without changing nuclearity of its associated $C^\ast$-algebra. Reduced hypergraphs will turn out to be more amenable to combinatorial reasoning. \textbf{Section \ref{refor_sec}} combines the previous results to obtain a proof of the main theorem. Finally, in \textbf{Section \ref{exa_sec}} we demonstrate in three cases how to determine whether a hypergraph $C^\ast$-algebra is nuclear. In particular, we retain the fact that every graph $C^\ast$-algebra associated to a finite graph is nuclear. We end the paper with a discussion of open problems in \textbf{Section \ref{pro_sec}}.

\subsection*{Acknowledgements}

This work has been part of the first author's master's thesis under the supervision of the second author. Both authors were supported by the SFB-TRR 195 Symbolic Tools in Mathematics and their Applications of the German Research Foundation (DFG). MW was further supported by the Heisenberg program of the DFG and a joint OPUS-LAP grant with Adam Skalski.

\section{Main Definitions and Main Results}
\label{mainres_sec}

\subsection{Notation}

Throughout this paper, the variables $A, B, C, D$ denote $C^\ast$-algebras. 
Every map between $C^\ast$-algebras is a $\ast$-homomorphism.
A map $\varphi: A \to B$ is called an embedding if it is injective.
Whenever there is an embedding of $B$ into $A$ we write $B \subset A$, and if $A$ and $B$ are $\ast$-isomorphic (Morita equivalent) we write $A = B$ ($A =_M B$).
An ideal $I \subset A$ is generally two-sided and closed. 
If $S \subset A$ is a subset of a $C^\ast$-algebra, then $(S)$ denotes the ideal generated by $S$. 
$M_k$ is the matrix algebra of dimension $k$, and we denote its standard matrix units consistently $E_{ij}$ for $i,j \leq k$. 
We write $S^1$ for the unit circle in the complex plane, $\Z_n$ for the group $\Z / n \Z$ and $\mathbb{F}_n$ for the free group with $n$ generators.
For unital $C^\ast$-algebras $A$ and $B$, $A *_\C B$ denotes their full unital free product.
If $p$ ($q$) is a projection in $A$ ($B$), then we write $A *_{p=q} B$ for the full amalgamated free product of $A$ and $B$ where the amalgamation is over the $C^\ast$-algebra $\C$ with the embeddings of $\C$ into $A$ and $B$ given by $z \mapsto zp$ and $z \mapsto zq$, respectively.
We refer to \cite{blackadarOperatorAlgebrasTheory2006} for the definition of full (amalgamated) free products of $C^\ast$-algebras.

\subsection{Main Definitions}

In this section, we present the definition of a hypergraph and its associated $C^\ast$-algebra.

\begin{definition}[hypergraph] \label{pre_hyp_def}
    A hypergraph $\Ha\Gamma$ is a tuple $(E^0, E^1, r, s)$, where
    \begin{itemize}
        \item $E^0 = E^0(\Ha\Gamma)$ is the (countable) set of vertices of $\Ha\Gamma$,
        \item $E^1 = E^1(\Ha\Gamma)$ is the (countable) set of edges of $\Ha\Gamma$,
        \item $r = r_{\Ha\Gamma}: E^1 \to \mathcal{P}(E^0)$ maps every edge to its \emph{range} (set),
        \item $s = s_{\Ha\Gamma}: E^1 \to \mathcal{P}(E^1) \setminus \{\emptyset\}$ maps every edge to its \emph{source} (set).
    \end{itemize}
    We call $\Ha\Gamma$ \emph{finite} if both $E^0$ and $E^1$ are finite sets, and we call $\Ha\Gamma$ \emph{undirected} if every edge $e \in E^1(\Ha\Gamma)$ has empty range.
\end{definition}

Throughout this paper we will only be interested in finite hypergraphs. Our notion of a hypergraph differs from definitions found in the literature. For instance, the directed hypergraphs from \cite{ausielloDirectedHypergraphsIntroduction2017} have (in our language) the additional requirement that $|r(e)| = 1$ holds for all edges $e \in E^1$, while the directed hypergraphs from \cite{galloDirectedHypergraphsApplications1993} are required to satisfy $s(e) \cap r(e) = \emptyset$ for all $e \in E^1$.

We say that an edge $e \in E^1$ starts from a set $V \subset E^0$ if $s(e) \cap V \neq \emptyset$. A vertex $v \in E^0$ is called a sink if there is no edge $e \in E^1$ with $v \in s(e)$. A path of length $n$ is a sequence $\mu = e_1 \dots e_n$ of edges with $r(e_i) \cap s(e_{i+1}) \neq \emptyset$ for all $i < n$ and the vertices are considered as paths of length zero. One sets $s_\mu = s_{e_1} \cdots s_{e_n} \in C^\ast(\Ha\Gamma)$ and $s(\mu) = s(e_1), r(\mu) = r(e_n)$. Further, $\mu$ is called \emph{closed} if $\mu$ has non-zero length and $r(\mu) \cap s(\mu) \neq \emptyset$. 
The path $\mu$ is a cycle if $\mu$ is closed and if $r(e_i) \cap s(e_j) = \emptyset$ unless $j = i+1 \leq n$ or $1 = j \leq i = n$.

The following definition of a hypergraph $C^\ast$-algebra is a slight generalization of the definition introduced in \cite{triebHypergraphAlgebras2024}.

\begin{definition}[hypergraph $C^\ast$-algebra]
    Let $\Ha\Gamma$ be a finite hypergraph. The hypergraph $C^\ast$-algebra $C^\ast(\Ha\Gamma)$ is the universal $C^\ast$-algebra generated by pairwise orthogonal projections $p_v$ and partial isometries $s_e$ for $v \in E^0$, $e \in E^1$, respectively, such that the following holds:
    \begingroup
    \leqnomode
    \begin{align*}
        \tag{HR1} s_e^\ast s_f &= \begin{cases}
            \delta_{ef} \sum_{v \in r(e)} p_v,  &r(e) \neq \emptyset, \\
            \delta_{ef} s_e,                    &\text{otherwise},
        \end{cases} \quad \text{ for all } e, f \in E^1, \\
        \tag{HR2} s_e s_e^\ast &\leq \sum_{v \in s(e)} p_v \text{ for all } e \in E^1, \\
        \tag{HR3} p_v &\leq \sum_{e: v \in s(e)} s_e s_e^\ast \text{ for all } v \in E^0 \text{ such that } v \text{ is not a sink}. 
    \end{align*}
    \endgroup
\end{definition}

Note that the hypergraph $C^\ast$-algebra $C^\ast(\Ha\Gamma)$ exists as a universal $C^\ast$-algebra since it is generated by projections and partial isometries.

The following proposition from \cite{triebHypergraphAlgebras2024} shows that hypergraph $C^\ast$-algebras truly generalize graph $C^\ast$-algebras. Here, we use the definition of a graph $C^\ast$-algebra from \cite{raeburnGraphAlgebras2005} with the role of ranges and sources swapped.

\begin{proposition}[{\cite[Proposition 3.4]{triebHypergraphAlgebras2024}}]
    Let $\Gamma = (E^0, E^1, r, s)$ be a finite, directed graph. 
    Further, obtain the hypergraph $\Ha\Gamma = (E^0, E^1, \hat{r}, \hat{s})$ by setting $\hat{s}(e) = \{s(e)\}$ and $\hat{r}(e) = \{r(e)\}$ for all edges $e \in E^1$. 
    Then the graph $C^\ast$-algebra $C^\ast(\Gamma)$ is $\ast$-isomorphic to the hypergraph $C^\ast$-algebra $C^\ast(\Ha\Gamma)$.
\end{proposition}

\subsection{Main Results}


Leaving the details to Section \ref{min_sec}, we call $\Ha\Delta$ a (hypergraph) minor of $\Ha\Gamma$, written $\Ha\Delta \leq \Ha\Gamma$, if the former is obtained from the latter by a combination of the following minor operations:
\begin{itemize}
    \item edge/vertex deletion
    \item forward/backward edge contraction
    \item edge cutting
    \item source separation
    \item range decomposition
\end{itemize}

In Table \ref{main_res_for_min_tbl} we list and sketch four concrete hypergraphs $\Ha\Gamma_1, \Ha\Gamma_2, \Ha\Gamma_3, \Ha\Gamma_4$ which we will call the \emph{forbidden minors}. These hypergraphs turn out to account for non-nuclearity of a large portion of hypergraph $C^\ast$-algebras. Their $C^\ast$-algebras are easily determined.

\begingroup
\renewcommand{\arraystretch}{2.5}
\begin{table}
    \centering
    \begin{tabular}{l|lll|c}
    \multirow{2}{*}{$\Ha\Gamma_1$} & $E^0 = \{v_1, v_2, v_3\}$, & $s(e) = s(f) = E^0$,          & \multirow{2}{*}{}        & \multirow{2}{*}{\begin{tikzpicture}[baseline=0]
        \node [fill=black, circle, inner sep=1pt] at (-.5,.3) {};
        \node [fill=black, circle, inner sep=1pt] at (.5, .3) {};
        \node [fill=black, circle, inner sep=1pt] at (1.5, .3) {};
        \draw[-] (.5, .3) ellipse [x radius = 45pt, y radius=10pt];
        \draw[-] (-.2, .6) -- (-.6, 1.5);
        \draw[-] (1.2, .6) -- (1.6, 1.5);
    \end{tikzpicture}} \\
                                   & $E^1 = \{e, f\}$,          & $r(e) = r(f) = \emptyset$                                               &                                                    &                                                                                                    \\ 
                                   \hline
    \multirow{2}{*}{$\Ha\Gamma_2$} & $E^0 = \{v_1, v_2\}$,      & $s(e) = s(f) = s(g) = E^0$,         & \multirow{2}{*}{} &   \multirow{2}{*}{\begin{tikzpicture}
        \node [fill=black, circle, inner sep=1pt] at (-.5,.3) {};
        \node [fill=black, circle, inner sep=1pt] at (.5, .3) {};
        \draw[-] (0, .3) ellipse [x radius = 30pt, y radius=10pt];
        \draw[-] (-.5, .6) -- (-.8, 1.5);
        \draw[-] (0, .65) -- (0, 1.6);        
        \draw[-] (.5, .6) -- (.8, 1.5);
    \end{tikzpicture} }                                                                                       \\
                                   & $E^1 = \{e, f, g\}$,       & $r(e) = r(f) = r(g) = \emptyset$                                               &                                                    &                                                                                                    \\ \hline
    \multirow{2}{*}{$\Ha\Gamma_3$} & $E^0 = \{v, w\}$,          & $s(e) = s(f) = E^0$,         &                                &         \multirow{2}{*}{\begin{tikzpicture}
        \node [fill=black, circle, inner sep=1pt] at (-.5,.3) {};
        \node [fill=black, circle, inner sep=1pt] at (.5, .3) {};
        \draw[-] (0, .3) ellipse [x radius = 30pt, y radius=15pt];
        \draw[-] (-.6, .73) -- (-.8, 1.5);
        \draw[<-] (0,.6) [partial ellipse=-10:170:.5cm and 1cm];
    \end{tikzpicture} }                                                                                           \\
                                   & $E^1 = \{e, f\}$,          & $r(e) = \emptyset, \; r(f) = \{w\}$                                              &                                    &                                                                                                    \\ \hline
    \multirow{2}{*}{$\Ha\Gamma_4$} & $E^0 = \{v_1, v_2, w\}$,   & $s(e) = s(f) = \{v_1, v_2\}$, & \multirow{2}{*}{}            &                                  \multirow{2}{*}{\begin{tikzpicture}
        \node [fill=black, circle, inner sep=1pt] at (-.5,.3) {};
        \node [fill=black, circle, inner sep=1pt] at (.5, .3) {};
        \draw[-] (0, .3) ellipse [x radius = 30pt, y radius=10pt];
        \node [fill=black, circle, inner sep=1pt] at (0,1.6) {};
        \draw[->] (-.5, .6) -- (-.1, 1.5);       
        \draw[->] (.5, .6) -- (.1, 1.5);
    \end{tikzpicture} }                                                                  \\
                                   & $E^1 = \{e, f\}$,          & $r(e) = r(f) = \{w\}$                                               &                                                    &                                                                                                   
    \end{tabular}
    \caption{The Forbidden Minors $\Ha\Gamma_1, \Ha\Gamma_2, \Ha\Gamma_3, \Ha\Gamma_4$} \label{main_res_for_min_tbl}
\end{table}

\endgroup

\begin{samepage}
\begin{proposition} \label{main_res_forbidden_min_algs}
    We have
    \begin{enumerate}
        \item $C^\ast(\Ha\Gamma_1) = C^\ast(\Ha\Gamma_2) = \C^2 *_\C \C^3$,
        \item $C^\ast(\Ha\Gamma_3)$ is the universal unital $C^\ast$-algebra generated by one partial isometry,
        \item $C^\ast(\Ha\Gamma_4)$ is Morita-equivalent to $M_2 *_\C \C^2 $.
    \end{enumerate}
    The $C^\ast$-algebras $C^\ast(\Ha\Gamma_1), C^\ast(\Ha\Gamma_2), C^\ast(\Ha\Gamma_3)$ are not exact while $C^\ast(\Ha\Gamma_4)$ is not nuclear.
\end{proposition}
\end{samepage}

\begin{proof}
    Ad (1): One readily checks that $C^\ast(\Ha\Gamma_1)$ is
    \begin{align*}
        &C^\ast\left(p_{v_1}, p_{v_2}, p_{v_3}, s_{e}, s_{f} \; \left| \;
        \begin{array}{ll}
            p_{v_1}, p_{v_2}, p_{v_3} \text{ pairwise orthogonal projections}, \\
            s_{e}, s_{f} \text{ orthogonal projections}, \\
            p_{v_1} + p_{v_2} + p_{v_3} = s_{e} + s_{f}
        \end{array} 
        \right. \right) \\
    &\qquad= \C^2 *_\C \C^3.  
    \end{align*}
    Analogously, one sees $C^\ast(\Ha\Gamma_2) = \C^2 *_\C \C^3$.
    It is a well-known fact from group theory that the free group with two generators $\mathbb{F}_2$ is a subgroup of the free product $\Z_2 * \Z_3$. Therefore,
    \begin{align}
        C^\ast(\mathbb{F}_2) \subset C^\ast(\Z_2 * \Z_3) = \C^2 *_\C \C^3.
    \end{align}
    Since $C^\ast(\mathbb{F}_2)$ is not exact (see \cite{wassermannTensorProductsCertain1976}) and exactness transfers to subalgebras, it follows that $\C^2 *_\C \C^3$ is not exact.

    Ad (2): An inspection of $C^\ast(\Ha\Gamma_3)$ reveals that all its generators can be expressed in terms of the unit and the partial isometry $s_f$. Therefore, $C^\ast(\Ha\Gamma_3)$ is the universal unital $C^\ast$-algebra generated by one partial isometry. In \cite{brenkenAlgebraPartialIsometry2011} it is proved that the universal unital $C^\ast$-algebra $\mathcal{P}$ generated by one partial isometry is not exact. Clearly, $\mathcal{P} \subset C^\ast(\Ha\Gamma_3)$. Hence, $C^\ast(\Ha\Gamma_3)$ is not exact since exactness transfers to subalgebras.

    Ad (3): Let $\Ha\Delta$ be the graph given by
    \begin{multicols}{2}

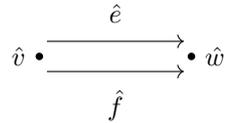
    \begin{itemize}
        \item $E^0(\Delta) = \{\hat{v}, \hat{w}\}$,
        \item $E^1(\Delta) = \{\hat{e}, \hat{f}\}$,
        \item $s_\Delta(\hat{e}) = s_\Delta(\hat{f}) = \{\hat{v}\}$,
        \item $r_\Delta(\hat{e}) = r_\Delta(\hat{f}) = \{\hat{w}\}$.
    \end{itemize}
    \columnbreak

    \begin{tikzpicture}[baseline=0]
        \node [fill=black, circle, inner sep=1pt, label=180:$\hat{v}$] at (0,0) {};
        \node [fill=black, circle, inner sep=1pt, label=0:$\hat{w}$] at (2,0) {};
        \draw[->] (.1, .2) -- (1.9, .2);
        \draw[->] (.1, -.2) -- (1.9, -.2);
        \node [label=90:$\hat{e}$] at (1,.2) {};
        \node [label=-90:$\hat{f}$] at (1,-.2) {};
    \end{tikzpicture}

\end{multicols}
    On the right-hand side above, we sketch the hypergraph $\Delta$.
    From \cite[Proposition 1.18]{raeburnGraphAlgebras2005} we obtain $C^\ast(\Delta) = M_3$ where one makes the identifications $p_{\hat{w}} = E_{11}, p_{\hat{v}} = E_{22} + E_{33}, s_{\hat{e}} = E_{21}, \hat{f} = E_{31}$. Now, one verifies
    \begin{align*}
        &C^\ast(\Ha\Gamma_4) \\
            &= C^\ast\left(p_{{v_1}}, p_{{v_2}}, p_w, s_e, s_f \; \left| \;
            \begin{array}{ll}
                p_{{v_1}}, p_{{v_2}}, p_w \text{ are pairwise orthogonal projections}, \\
                s_e, s_f \text{ are partial isometries}, \\
                s_e^\ast s_e = s_f^\ast s_f = p_w, \\
                s_e^\ast s_f = 0, \\
                p_{{v_1}} + p_{{v_2}} = s_e s_e^\ast + s_f s_f^\ast
            \end{array} 
            \right. \right)
            \\
            &= C^\ast\left(p_v, p_{{v_1}}, p_{{v_2}}, p_w, s_e, s_f \; \left| \;
            \begin{array}{ll}
                p_v, p_w \text{ are orthogonal projections}, \\
                p_{{v_1}}, p_{{v_2}} \text{ are orthogonal projections}, \\
                s_e, s_f \text{ are partial isometries}, \\
                s_e^\ast s_e = s_f^\ast s_f = p_w, \\
                s_e^\ast s_f = 0, \\
                p_v = s_e s_e^\ast + s_f s_f^\ast, \\
                p_v = p_{{v_1}} + p_{{v_2}}
            \end{array} 
            \right. \right)
            \\
            &= C^\ast\left(p_{\hat{v}}, p_{\hat{w}}, s_{\hat{e}}, s_{\hat{f}} \; \left| \;
            \begin{array}{ll}
                p_{\hat{v}}, p_{\hat{w}} \text{ are orthogonal projections}, \\
                s_{\hat{e}}, s_{\hat{f}} \text{ are partial isometries}, \\
                s_{\hat{e}}^\ast s_{\hat{e}} = s_{\hat{f}}^\ast s_{\hat{f}} = p_w, \\
                s_{\hat{e}}^\ast s_{\hat{f}} = 0, \\
                p_{\hat{v}} = s_{\hat{e}} s_{\hat{e}}^\ast + s_{\hat{f}} s_{\hat{f}}^\ast
            \end{array} \right. \right)
            \\ & \qquad
            *_{p_{\hat{v}} = p_{{v_1}} + p_{{v_2}}}
            C^\ast\left(p_{{v_1}}, p_{{v_2}} \; \left| \;
            \begin{array}{ll}
                p_{{v_1}}, p_{{v_2}} \text{ are orthogonal projections}
            \end{array} 
            \right. \right)
            \\
            &= C^\ast(\Delta) *_{p_{\hat{v}} = 1_{\C^2}} \C^2
            \\
            &= M_3 *_{E_{22} + E_{33} = 1_{\C^2}} \C^2.
        \end{align*}
    The latter $C^\ast$-algebra contains 
    $$ M_2 *_\C \C^2 = (E_{22} + E_{33}) (M_3 *_{E_{22} + E_{33} = 1_{\C^2}} \C^2) (E_{22} + E_{33}) $$ 
    as a full corner. Therefore, $C^\ast(\Ha\Gamma_4)$ is Morita-equivalent to $M_2 *_\C \C^2$. However, this $C^\ast$-algebra is non-nuclear by \cite[Proposition 3+6]{albeverioWildnessClassesAlgebras2006}. The claim follows as nuclearity is preserved under Morita equivalence.
\end{proof}

Finally, let us state the main theorem.

\begin{samepage}
\begin{thm} \label{main_res}
    Let $\Ha\Gamma$ be a hypergraph. One can construct a 
    hypergraph $\Ha\Delta \leq \Ha\Gamma$ such that $C^\ast(\Ha\Gamma)$ is nuclear if, and only if, the same holds for $C^\ast(\Ha\Delta)$. Further, the following is true:
    \begin{enumerate}
        \item If $\Ha\Gamma_i \leq \Ha\Delta$ holds for some $i \leq 3$, then $C^\ast(\Ha\Gamma)$ is not exact.
        \item If $\Ha\Gamma_4 \leq \Ha\Delta$ holds, then $C^\ast(\Ha\Gamma)$ is not nuclear.
        \item If none of the above holds, then $\Ha\Delta$ is an undirected hypergraph, i.e. all edges of $\Ha\Delta$ have empty range.
    \end{enumerate}
    Here $\Ha\Gamma_1, \Ha\Gamma_2, \Ha\Gamma_3$ and $\Ha\Gamma_4$ are the forbidden minors from Table \ref{main_res_for_min_tbl}.
\end{thm}
\end{samepage}

Crucially, there is an explicit procedure for constructing $\Ha\Delta$ starting from $\Ha\Gamma$ which we describe in Section \ref{red_sec} as Algorithm \ref{pro_nuc_reduction_algo}.
The proof of the previous theorem is done in Section \ref{refor_mainres_proof_subsec}. 

\begin{remark} \label{mainres_remark_undirected_hyp_algs}
    If $\Ha\Gamma$ has none of the forbidden minors, then one checks that $\Ha\Delta$ satisfies the following:
\begin{itemize}
    \item Every edge has empty range.
    \item For any distinct vertices $v,w \in E^0(\Ha\Delta)$ there are at most two distinct edges $e, f$ with $\{v,w\} \subset s(e) \cap s(f)$.
    \item For any distinct edges $e, f \in E^1(\Ha\Delta)$ there are at most two distinct vertices $v, w$ with $\{v,w\} \subset s(e) \cap s(f)$.
\end{itemize}

We believe that for any hypergraph $\Ha\Delta$ with these properties the associated $C^\ast$-algebra $C^\ast(\Ha\Delta)$ is nuclear. Then it would follow that any hypergraph $C^\ast$-algebra $C^\ast(\Ha\Gamma)$ is nuclear if, and only if, $\Ha\Delta$ has none of the forbidden minors.
\end{remark}

\section{Hypergraph Minors}
\label{min_sec}

\subsection{Definition and Main Result of this Section}

This section introduces the notion of a hypergraph minor and discusses the effect of the minor operations on the $C^\ast$-algebra side.

\begin{definition}[hypergraph minor] \label{min_minor_def}
    We say that $\Ha \Delta$ is obtained from $\Ha \Gamma$ by 
    \begin{enumerate}
        \item \emph{vertex deletion}\index{vertex deletion} if there is a vertex $v$ in $\Ha \Gamma$ such that 
        \begin{itemize}
            \item $E^0(\Ha\Delta) = E^0(\Ha\Gamma) \setminus \{v\}$,
            \item $E^1(\Ha\Delta) = E^1(\Ha\Gamma) \setminus \{e \in E^1(\Ha\Gamma): \, s(e) = \{v\} \}$,
            \item $r_{\Ha\Delta}(e) = r_{\Ha\Gamma}(e) \setminus \{v\}$ for all $e \in E^1(\Ha\Delta)$,
            \item $s_{\Ha\Delta}(e) = s_{\Ha\Gamma}(e) \setminus \{v\}$ for all $e \in E^1(\Ha\Delta)$,
        \end{itemize}
        \item \emph{edge deletion}\index{edge!deletion} if there is an edge $f$ in $\Ha\Gamma$ such that
        \begin{itemize}
            \item $E^0(\Ha\Delta) = E^0(\Ha\Gamma)$,
            \item $E^1(\Ha\Delta) = E^1(\Ha\Gamma) \setminus \{f\}$,
            \item $r_{\Ha\Delta}(e) = r_{\Ha\Gamma}(e)$ for all $e \in E^1(\Ha\Delta)$,
            \item $s_{\Ha\Delta}(e) = s_{\Ha\Gamma}(e)$ for all $e \in E^1(\Ha\Delta)$,
        \end{itemize}
        \item \emph{forward edge contraction}\index{edge!contraction} if there is an edge $f$ and a vertex $w$ in $\Ha\Gamma$ with $s_{\Ha\Gamma}(f) = \{w\}$ and 
        \begin{itemize}
            \item $f$ is the only edge starting from $w$ in $\Ha\Gamma$, i.e. 
            $$ s_{\Ha\Gamma}(e) \cap s_{\Ha\Gamma}(f) \neq \emptyset \implies e = f \quad \text{for all } e \in E^1(\Ha\Gamma), $$
            \item there is no edge $e \in E^1(\Ha\Gamma)$ with $w \in r_{\Ha\Gamma}(e)$ and $r_{\Ha\Gamma}(e) \cap r_{\Ha\Gamma}(f) \neq \emptyset$,
            \item $E^0(\Ha\Delta) = E^0(\Ha\Gamma) \setminus \{w\}$,
            \item $E^1(\Ha\Delta) = E^1(\Ha\Gamma) \setminus \{f\}$,
            \item $r_{\Ha\Delta}(e) = \begin{cases}
                r_{\Ha\Gamma}(e), &w \not \in r_{\Ha\Gamma}(e), \\
                (r_{\Ha\Gamma}(e) \setminus \{w\}) \cup r_{\Ha\Gamma}(f), &\text{otherwise},
            \end{cases} \quad$ for all $e \in E^1(\Ha\Delta)$, 
            \item $s_{\Ha\Delta}(e) = s_{\Ha\Gamma}(e)$ for all $e \in E^1(\Ha\Delta)$,
        \end{itemize}
    \begin{samepage}
        \item \emph{backward edge contraction} if there is an edge $f$ and a vertex $w$ in $\Ha\Gamma$ with $r_{\Ha\Gamma}(f) = \{w\}$ and 
        \begin{itemize}
            \item $f$ is the only edge starting from $s_{\Ha\Gamma}(f)$ in $\Ha\Gamma$, i.e. 
            $$ s_{\Ha\Gamma}(e) \cap s_{\Ha\Gamma}(f) \neq \emptyset \implies e = f \quad \text{for all } e \in E^1(\Ha\Gamma), $$
            \item there is no edge $e \in E^1(\Ha\Gamma)$ with $r_{\Ha\Gamma}(e) \cap s_{\Ha\Gamma}(f) \neq \emptyset$ and $w \in r_{\Ha\Gamma}(e)$,
            \item $E^0(\Ha\Delta) = E^0(\Ha\Gamma) \setminus \{w\}$,
            \item $E^1(\Ha\Delta) = E^1(\Ha\Gamma) \setminus \{f\}$,
            \item $r_{\Ha\Delta}(e) = \begin{cases}
                r_{\Ha\Gamma}(e), &w \not \in r_{\Ha\Gamma}(e), \\
                (r_{\Ha\Gamma}(e) \setminus \{w\}) \cup s_{\Ha\Gamma}(f), &\text{otherwise},
            \end{cases} \quad$ for all $e \in E^1(\Ha\Delta)$,
            \item $s_{\Ha\Delta}(e) = \begin{cases}
                s_{\Ha\Gamma}(e), &w \not \in s_{\Ha\Gamma}(e), \\
                (s_{\Ha\Gamma}(e) \setminus \{w\}) \cup s_{\Ha\Gamma}(f), &\text{otherwise},
            \end{cases} \quad$ for all $e \in E^1(\Ha\Delta)$,
        \end{itemize}
    \end{samepage}
        \item \emph{edge cutting}\index{edge!cutting} if there is an edge $f$ in $\Ha\Gamma$ such that
        \begin{itemize}
            \item $E^0(\Ha\Delta) = E^0(\Ha\Gamma)$,
            \item $E^1(\Ha\Delta) = E^1(\Ha\Gamma)$,
            \item $r_{\Ha\Delta}(e) = \begin{cases}
                r_{\Ha\Gamma}(e), &e \neq f, \\
                \emptyset, &e = f,
            \end{cases} \quad$ for all $e \in E^1(\Ha\Delta)$,
            \item $s_{\Ha\Delta}(e) = s_{\Ha\Gamma}(e)$ for all $e \in E^1(\Ha\Delta)$,
        \end{itemize}
        \item \emph{source separation}\index{source separation} if there is a set $F \subset E^1(\Ha\Gamma)$, a vertex $w \in E^0(\Ha\Gamma)$ and some vertex $w^\prime \in E^0(\Ha\Delta) \setminus E^0(\Ha\Gamma)$ such that \label{min_source_sep_def}
        \begin{itemize}
            \item $\emptyset \neq F \subsetneq \{e \in E^1(\Ha\Gamma): w \in s_{\Ha\Gamma}(e)\}$
            \item $E^0(\Ha\Delta) = E^0(\Ha\Gamma) \cup \{w^\prime\}$
            \item $E^1(\Ha\Delta) = E^1(\Ha\Gamma)$,
            \item $r_{\Ha\Delta}(e) = \begin{cases}
                r_{\Ha\Gamma}(e),                   &w \not \in r_{\Ha\Gamma}(e), \\
                r_{\Ha\Gamma}(e) \cup \{w^\prime\}, &\text{otherwise},
            \end{cases} \quad$ for all $e \in E^1(\Ha\Delta)$,
            \item $s_{\Ha\Delta}(e) = \begin{cases}
                s_{\Ha\Gamma}(e),                                         &e \not \in F, \\
                (s_{\Ha\Gamma}(e) \setminus \{w\}) \cup \{w^\prime\},     &e \in F,
            \end{cases} \quad$ for all $e \in E^1(\Ha\Delta)$,
        \end{itemize}
        \item \emph{range decomposition}\index{range decomposition} if there is an edge $f$ in $\Ha\Gamma$ with nonempty range such that
        \begin{itemize}
            \item $E^0(\Ha\Delta) = E^0(\Ha\Gamma)$,
            \item $E^1(\Ha\Delta) = (E^1(\Ha\Gamma) \setminus \{f\}) \cup \{(f,v): v \in r_{\Ha\Gamma}(f) \}$,
            \item $r_{\Ha\Delta}(e) = \begin{cases}
                r_{\Ha\Gamma}(e),           &e \not \in \{(f,v): v \in r_{\Ha\Gamma}(f)\}, \\
                \{v\},                      &e = (f,v),
            \end{cases} \quad$ for all $e \in E^1(\Ha\Delta)$,
            \item $s_{\Ha\Delta}(e) = s_{\Ha\Gamma}(e)$ for all $e \in E^1(\Ha\Delta)$.
        \end{itemize}
    \end{enumerate}
    The hypergraph $\Ha\Delta$ is a \emph{minor}\index{hypergraph!minor} of $\Ha\Gamma$, written $\Ha\Delta \leq \Ha\Gamma$, if it is obtained from $\Ha\Gamma$ by any finite combination of these operations.
\end{definition}

\begin{example}
    In Table \ref{min_example_minor_operations} we give an example for every minor operation except for edge and vertex deletion. We trust that the sketches are self-explanatory and do not give explicit definitions of the involved hypergraphs. Vertices or edges relevant to the respective minor operation are highlighted in red.
    
    \begingroup
    \renewcommand{\arraystretch}{6}
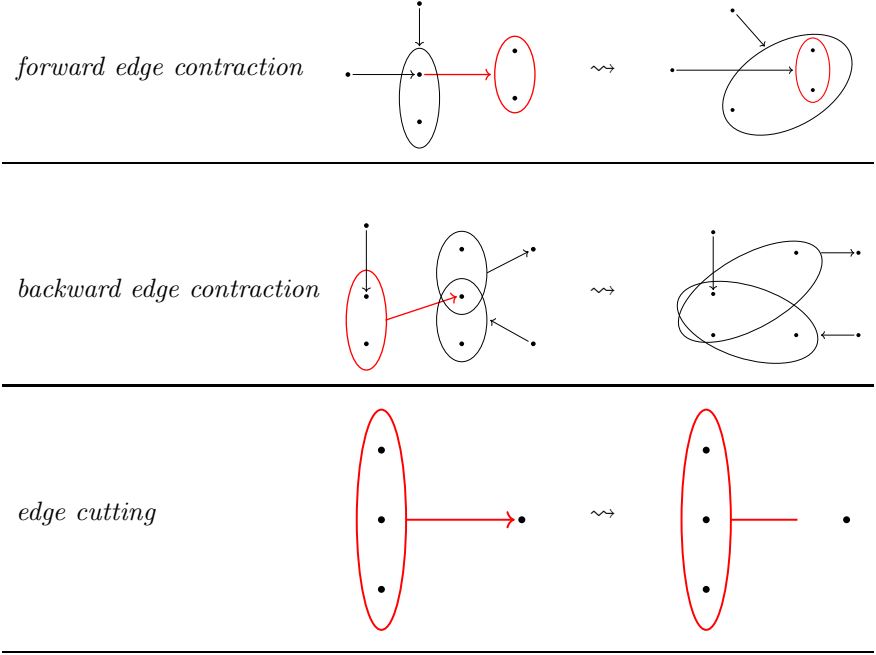
    \begin{center}
    \begin{longtable}{l  c c c}
        forward edge contraction
            &\resizebox{2.5cm}{!}{\begin{tikzpicture}[baseline = (current bounding box.center)]
                \node [fill=black, circle, inner sep=1pt] at (-1.5, 0) {};
                \node [fill=black, circle, inner sep=1pt] at (0, 1.5) {};
                \node [fill=black, circle, inner sep=1pt] at (0, 0) {};
                \node [fill=black, circle, inner sep=1pt] at (0, -1) {};
                \node [fill=black, circle, inner sep=1pt] at (2, .5) {};
                \node [fill=black, circle, inner sep=1pt] at (2, -.5) {};
                \draw[->, red, thick] (.1, 0) -- (1.5, 0);
                \draw[-, red, thick] (2, 0) ellipse [x radius = 12pt, y radius=23pt];
                \draw[->] (-1.4, 0) -- (-.1, 0);
                \draw[-] (0, -.5) ellipse [x radius = 12pt, y radius=30pt];
                \draw[->] (0, 1.4) -- (0, .6);
                \addvmargin{3mm}
            \end{tikzpicture}}
            &$\quad \rightsquigarrow \quad$
            &\resizebox{2.5cm}{!}{\begin{tikzpicture}[baseline = (current bounding box.center)]
                \node [fill=black, circle, inner sep=1pt] at (-1.5, 0) {};
                \node [fill=black, circle, inner sep=1pt] at (0, 1.5) {};
                \node [fill=black, circle, inner sep=1pt] at (0, -1) {};
                \node [fill=black, circle, inner sep=1pt] at (2, .5) {};
                \node [fill=black, circle, inner sep=1pt] at (2, -.5) {};
                \draw[-, red, thick] (2, 0) ellipse [x radius = 12pt, y radius=23pt];
                \draw[->] (-1.4, 0) -- (1.5, 0);
                \begin{scope}[rotate=30]
                     \draw[-] (1, -1) ellipse [x radius = 50pt, y radius=30pt];
                \end{scope}
                \draw[->] (.1, 1.4) -- (.8, .6);
                \addvmargin{3mm}
            \end{tikzpicture}} \\
        \hline 
        backward edge contraction
            &\resizebox{2.5cm}{!}{\begin{tikzpicture}[baseline = (current bounding box.center)]
                \node [fill=black, circle, inner sep=1pt] at (0, 1.5) {};
                \node [fill=black, circle, inner sep=1pt] at (0, 0) {};
                \node [fill=black, circle, inner sep=1pt] at (0, -1) {};
                \node [fill=black, circle, inner sep=1pt] at (2, 1) {};
                \node [fill=black, circle, inner sep=1pt] at (2, 0) {};
                \node [fill=black, circle, inner sep=1pt] at (2, -1) {};
                \node [fill=black, circle, inner sep=1pt] at (3.5, -1) {};
                \draw[-, red, thick] (0, -.5) ellipse [x radius = 12pt, y radius=30pt];
                \draw[->, red, thick] (.4, -.5) -- (1.9, 0);
                \draw[-] (2, .5) ellipse [x radius = 15pt, y radius=25pt];
                \draw[->] (2.53, .5) -- (3.4, .95);
                \node [fill=black, circle, inner sep=1pt] at (3.5, 1) {};
                \draw[-] (2, -.5) ellipse [x radius = 15pt, y radius=25pt];
                \draw[<-] (2.6, -.5) -- (3.4, -.95);
                \draw[->] (0, 1.4) -- (0, .1);
                \addvmargin{3mm}
            \end{tikzpicture}}
            &$\quad \rightsquigarrow \quad$
            &\resizebox{2.5cm}{!}{\begin{tikzpicture}[baseline = (current bounding box.center)]
                \node [fill=black, circle, inner sep=1pt] at (0, 1.5) {};
                \node [fill=black, circle, inner sep=1pt] at (0, 0) {};
                \node [fill=black, circle, inner sep=1pt] at (0, -1) {};
                \node [fill=black, circle, inner sep=1pt] at (2, 1) {};
                \node [fill=black, circle, inner sep=1pt] at (2, -1) {};
                \node [fill=black, circle, inner sep=1pt] at (3.5, -1) {};
                \begin{scope}[rotate=30]
                    \draw[-] (.8, -.4) ellipse [x radius = 55pt, y radius=25pt];
                \end{scope}
                \draw[->] (2.6, 1) -- (3.4, 1);
                \node [fill=black, circle, inner sep=1pt] at (3.5, 1) {};
                \begin{scope}[rotate=-18]
                    \draw[-] (1, -.4) ellipse [x radius = 50pt, y radius=25pt];
                \end{scope}
                \draw[<-] (2.6, -1) -- (3.4, -1);
                \draw[->] (0, 1.4) -- (0, .1);
                \addvmargin{3mm}
            \end{tikzpicture}} \\
        \hline 
        edge cutting
            & \resizebox{!}{1.75cm}{\begin{tikzpicture}[baseline = (current bounding box.center)]
                \node [fill=black, circle, inner sep=1pt] at (0, 1) {};
                \node [fill=black, circle, inner sep=1pt] at (0, 0) {};
                \node [fill=black, circle, inner sep=1pt] at (0, -1) {};
                \node [fill=black, circle, inner sep=1pt] at (2, 0) {};
                \draw[-, red, thick] (0, 0) ellipse [x radius = 10pt, y radius=45pt];
                \draw[->, red, thick] (.35, 0) -- (1.9, 0);
                \addvmargin{3mm}
            \end{tikzpicture}}
            &$\quad \rightsquigarrow \quad$
            &\resizebox{!}{1.75cm}{\begin{tikzpicture}[baseline = (current bounding box.center)]
                \node [fill=black, circle, inner sep=1pt] at (0, 1) {};
                \node [fill=black, circle, inner sep=1pt] at (0, 0) {};
                \node [fill=black, circle, inner sep=1pt] at (0, -1) {};
                \node [fill=black, circle, inner sep=1pt] at (2, 0) {};
                \draw[-, red, thick] (0, 0) ellipse [x radius = 10pt, y radius=45pt];
                \draw[-, red, thick] (.35, 0) -- (1.3, 0);
                \addvmargin{3mm}
            \end{tikzpicture}}\\
        \hline 
        source separation
            & \resizebox{2.5cm}{!}{\begin{tikzpicture}[baseline = (current bounding box.center)]
                \node [fill=black, circle, inner sep=1pt] at (0, 0) {};
                \node [fill=black, circle, inner sep=1pt] at (2, 1) {};
                \node [fill=red, circle, inner sep=1.5pt] at (2, 0) {};
                \node [fill=black, circle, inner sep=1pt] at (2, -1.5) {};
                \node [fill=black, circle, inner sep=1pt] at (3.5, -1) {};
                \draw[->] (.1, 0) -- (1.9, 0);
                \draw[-] (2, .5) ellipse [x radius = 15pt, y radius=25pt];
                \draw[->] (2.53, .5) -- (3.4, .95);
                \node [fill=black, circle, inner sep=1pt] at (3.5, 1) {};
                \draw[-, red, thick] (2, -.75) ellipse [x radius = 15pt, y radius=30pt];
                \draw[->, red, thick] (2.53, -.75) -- (3.4, -.95);
                \addvmargin{3mm}
            \end{tikzpicture}}
            &$\quad \rightsquigarrow \quad$
            &\resizebox{2.5cm}{!}{\begin{tikzpicture}[baseline = (current bounding box.center)]
                \node [fill=black, circle, inner sep=1pt] at (0, 0) {};
                \node [fill=black, circle, inner sep=1pt] at (2, 1) {};
                \node [fill=red, circle, inner sep=1.5pt] at (2, 0) {};
                \node [fill=red, circle, inner sep=1.5pt] at (2, -.75) {};
                \node [fill=black, circle, inner sep=1pt] at (2, -1.5) {};
                \node [fill=black, circle, inner sep=1pt] at (3.5, -1) {};
                \draw[->] (.1, -.01) -- (1.6, -.375);
                \draw[-] (2, -.375) ellipse [x radius = 10pt, y radius=20pt];
                \draw[-] (2, .5) ellipse [x radius = 15pt, y radius=25pt];
                \draw[->] (2.53, .5) -- (3.4, .95);
                \node [fill=black, circle, inner sep=1pt] at (3.5, 1) {};
                \draw[-, red, thick] (2, -1.15) ellipse [x radius = 15pt, y radius=20pt];
                \draw[->, red, thick] (2.53, -1) -- (3.4, -1);
                \addvmargin{3mm}
            \end{tikzpicture}}\\
        \hline 
        range decomposition
            & \resizebox{!}{1.25cm}{\begin{tikzpicture}[baseline = (current bounding box.center)]
                \node [fill=black, circle, inner sep=1pt] at (0, 0) {};
                \node [fill=black, circle, inner sep=1pt] at (2, 1) {};
                \node [fill=black, circle, inner sep=1pt] at (2, 0) {};
                \node [fill=black, circle, inner sep=1pt] at (2, -1) {};
                \draw[->, red, thick] (.1, 0) -- (1.6, 0);
                \draw[-, red, thick] (2, 0) ellipse [x radius = 10pt, y radius=45pt];
            \end{tikzpicture}}
            &$\quad \rightsquigarrow \quad$
            &\resizebox{!}{1.25cm}{\begin{tikzpicture}[baseline = (current bounding box.center)]
                \node [fill=black, circle, inner sep=1pt] at (0, 0) {};
                \node [fill=black, circle, inner sep=1pt] at (2, 1) {};
                \node [fill=black, circle, inner sep=1pt] at (2, 0) {};
                \node [fill=black, circle, inner sep=1pt] at (2, -1) {};
                \draw[->, red, thick] (.1, 0) -- (1.9, 1);
                \draw[->, red, thick] (.1, 0) -- (1.9, 0);
                \draw[->, red, thick] (.1, 0) -- (1.9, -1);
            \end{tikzpicture}} \\
            \caption{Examples for the Minor Operations}
            \label{min_example_minor_operations}
    \end{longtable}
    \end{center}
    \endgroup
        
\end{example}

The next theorem describes the effect of the minor operations on the associated $C^\ast$-algebras.

\begin{thm} \label{min_big_thm}
    Let $\Ha\Gamma$ and $\Ha\Delta$ be two hypergraphs. The algebra $C^\ast(\Ha\Delta)$ is
    \begin{itemize}
        \item isomorphic to $C^\ast(\Ha\Gamma)$ if $\Ha\Delta$ is obtained from $\Ha\Gamma$ by range decomposition,
        \item a quotient of $C^\ast(\Ha\Gamma)$ if $\Ha\Delta$ is obtained from $\Ha\Gamma$ by source separation,
        \item a subalgebra of $C^\ast(\Ha\Gamma)$ if $\Ha\Delta$ is obtained from $\Ha\Gamma$ by edge cutting,
        \item a quotient of a subalgebra of $C^\ast(\Ha\Gamma)$ if $\Ha\Delta$ is obtained from $\Ha\Gamma$ by edge or vertex deletion,
        \item a full corner of $C^\ast(\Ha\Gamma)$ if $\Ha\Delta$ is obtained from $\Ha\Gamma$ by forward or backward edge contraction.
    \end{itemize}
    In particular, $\Ha\Delta \leq \Ha\Gamma$ implies that $C^\ast(\Ha\Delta)$ is -- up to Morita equivalence -- obtained from $C^\ast(\Ha\Gamma)$ by alternatingly taking subalgebras and quotients. If $C^\ast(\Ha\Gamma)$ is exact then the same holds for $C^\ast(\Ha\Delta)$.
\end{thm}

In the next subsections we will investigate the different minor operations separately. We will not deal with range decomposition and instead refer to \cite[Proposition 3.18]{triebHypergraphAlgebras2024}, where this operation has been studied. This will combine to the proof of the previous theorem at the end of Section \ref{min_edge_contraction_subsec}.

\subsection{Source Separation}

On the $C^\ast$-algebra side, source separation corresponds to taking a quotient. Under special conditions this improves to an isomorphism of the associated hypergraph $C^\ast$-algebras.

\begin{proposition}[source separation] \label{min_edge_sep}
    Let $\Ha\Delta$ be obtained from $\Ha\Gamma$ by source separation of a nonempty set $F \subset E^1(\Ha\Gamma)$ 
    at $w \in \bigcap_{f \in F} s_{\Ha\Gamma}(f)$. Then $C^\ast(\Ha\Delta)$ is a quotient of $C^\ast(\Ha\Gamma)$. Moreover, if  
    \begin{align*}
        w \in s_{\Ha\Gamma}(g) \implies \{w\} = s_{\Ha\Gamma}(f) \cap s_{\Ha\Gamma}(g) \qquad \text{for all } g \not \in F, \tag{$\ast$}
    \end{align*}
    then $C^\ast(\Ha\Delta) = C^\ast(\Ha\Gamma)$.
\end{proposition}

\begin{proof}
    Let $w^\prime \in E^0(\Ha\Delta)$ be as in Definition \ref{min_minor_def}(\ref{min_source_sep_def}).
    Use the universal property of $C^\ast(\Ha\Gamma)$ to obtain $\varphi: C^\ast(\Ha\Gamma) \to C^\ast(\Ha\Delta)$ with
    \begin{align*}
        \varphi: 
        \left\{
        \begin{aligned}
            p_v &\mapsto \hat{p}_v,                     && v \in E^0(\Ha\Gamma) \setminus \{w\}, \\
            p_v &\mapsto \hat{p}_w + \hat{p}_{w^\prime},  && v = w,  \\
            s_e &\mapsto \hat{s}_e,                       && e \in E^1(\Ha\Gamma),
        \end{aligned}
        \right.
    \end{align*}
    
    

    Clearly the $\varphi(p_v)$ are pairwise orthogonal projections and the $\varphi(s_e)$ are partial isometries. One checks the following relations:
    (HR1): For $e, e^\prime \in E^1(\Ha\Gamma)$ it is
    \begin{align*}
        \varphi(s_e)^\ast \varphi(s_{e^\prime})
            = \hat{s}_e^\ast \hat{s}_{e^\prime}
            = \begin{cases}
                \delta_{e e^\prime} \sum_{v \in r_{\Ha\Delta}(e)} \hat{p}_v = \delta_{e e^\prime} \varphi(\sum_{v \in r_{\Ha\Gamma}(e)} p_v),     & r_{\Ha\Gamma}(e) \neq \emptyset, \\
                \delta_{e e^\prime} \hat{s}_e = \delta_{e e^\prime} \varphi(s_e),                                   &\text{otherwise}.
            \end{cases}
    \end{align*}
    (HR2): For $e \in E^1(\Ha\Gamma)$ we have
    \begin{align*}
        \varphi(s_e) \varphi(s_e)^\ast
            = \hat{s}_e \hat{s}_e^\ast
            \leq \sum_{v \in s_{\Ha\Delta}(e)} \hat{p}_v
            \leq \varphi\left(\sum_{v \in s_{\Ha\Gamma}(e)} p_v\right).
    \end{align*}
    (HR3): Let $v \in E^0(\Ha\Gamma) \setminus \{w\}$ not be a sink in $\Ha\Gamma$. Then $v$ is not a sink in $\Ha\Delta$ as well and therefore
    \begin{align*}
        \varphi(p_v)
            = \hat{p}_v
            \leq \sum_{e \in E^1(\Ha\Delta): v \in s_{\Ha\Delta}(e)} \hat{s}_e \hat{s}_e^\ast 
            = \sum_{e \in E^1(\Ha\Gamma): v \in s_{\Ha\Gamma}(e)} \varphi(s_e) \varphi(s_e)^\ast.
    \end{align*}
    For the vertex $w$ observe that neither $w$ nor $w^\prime$ is a sink in $\Ha\Delta$ and thus
    \begin{align*}
        \varphi(p_w)
            &= \hat{p}_w + \hat{p}_{w^\prime} \\
            &\leq \sum_{e \in E^1(\Ha\Delta): w \in s_{\Ha\Delta}(e) \vee w^\prime \in s_{\Ha\Delta}(e)} \hat{s}_e \hat{s}_e^\ast \\
            &= \sum_{e \in E^1(\Ha\Gamma): w \in s_{\Ha\Gamma}(e)} \varphi(s_e) \varphi(s_e)^\ast.
    \end{align*}

    Finally, observe
    \begin{align*}
        \hat{p}_{w^\prime} = (\hat{p}_w + \hat{p}_{w^\prime}) \sum_{f \in F} \hat{s}_f \hat{s}_f^\ast = \varphi(p_w) \sum_{f \in F} \varphi(s_f) \varphi(s_f)^\ast \in \varphi(C^\ast(\Ha\Gamma)).
    \end{align*}
    Therefore, the image of $\varphi$ contains all generators of $C^\ast(\Ha\Delta)$. Hence, $\varphi$ is surjective and $C^\ast(\Ha\Delta) = C^\ast(\Ha\Gamma) / \mathrm{ker}(\varphi)$.

    Under the additional assumption $(\ast)$ the universal property of $C^\ast(\Ha\Delta)$ yields a map $\psi: C^\ast(\Ha\Delta) \to C^\ast(\Ha\Gamma)$ with
    \begin{align*}
        \psi: 
        \left\{
        \begin{aligned}
            \hat{p}_v &\mapsto p_v,                                                        && v \in E^0(\Ha\Delta) \setminus \{w, w^\prime\}, \\
            \hat{p}_v &\mapsto \left(\sum_{g \in G} s_g s_g^\ast\right) p_w \left(\sum_{g \in G} s_g s_g^\ast\right),    && v = w,  \\
            \hat{p}_v &\mapsto \left(\sum_{f \in F} s_f s_f^\ast\right) p_w \left(\sum_{f \in F} s_f s_f^\ast\right),    && v = w^\prime,  \\
            \hat{s}_e &\mapsto s_e,                     && e \in E^1(\Ha\Delta).
        \end{aligned}
        \right.
    \end{align*}
    Indeed, observe for any $f \in F$ and $g \in G$
    \begin{align*}
        0 = s_f s_f^\ast s_g s_g^\ast = s_f s_f^\ast \left(\sum_{v \in s_{\Ha\Gamma}(f)} p_v \right) \left(\sum_{v \in s_{\Ha\Gamma}(g)} p_v \right) s_g s_g^\ast
            = s_f s_f^\ast p_w s_g s_g^\ast.
    \end{align*}
    Therefore, we have
    \begin{align*}
        p_w &= \left( \sum_{e: w \in s_{\Ha\Gamma}(w)} s_e s_e^\ast \right) p_w \left( \sum_{e: w \in s_{\Ha\Gamma}(w)} s_e s_e^\ast \right) \\
            &= \left(\sum_{f \in F} s_f s_f^\ast\right) p_w \left(\sum_{f \in F} s_f s_f^\ast\right)
                + \left(\sum_{g \in G} s_g s_g^\ast\right) p_w \left(\sum_{g \in G} s_g s_g^\ast\right)
    \end{align*}
    and 
    \begin{align*}
        \psi(\hat{p}_{w^\prime})^2 
            &= \left( \left(\sum_{f \in F} s_f s_f^\ast\right) p_w \left(\sum_{f \in F} s_f s_f^\ast\right) \right)^2 \\
            &= \left(\sum_{f \in F} s_f s_f^\ast\right) p_w \left(\sum_{f \in F} s_f s_f^\ast\right) p_w \left(\sum_{f \in F} s_f s_f^\ast\right) \\
            &= \left(\sum_{f \in F} s_f s_f^\ast\right) p_w \left(\sum_{e \in F \cup G} s_e s_e^\ast\right) p_w \left(\sum_{f \in F} s_f s_f^\ast\right) \\
            &= \left(\sum_{f \in F} s_f s_f^\ast\right) p_w \left(\sum_{f \in F} s_f s_f^\ast\right) \\
            &= \psi(\hat{p}_{w^\prime})
    \end{align*}
    Similarly, one sees that $\psi(\hat{p}_w)$ is a projection, and it is easily checked that $\psi(\hat{p}_w)$ and $\psi(\hat{p}_{w^\prime})$ are orthogonal. Thus, the $\psi(\hat{p}_v)$ are pairwise orthogonal projections and evidently the $\psi(\hat{s}_e)$ are partial isometries. We check the relations (HR1), (HR2) and (HR3).
        
    (HR1): Since in $\Ha\Delta$ every edge $e$ with $w \in r_{\Ha\Delta}(e)$ has also $w^\prime$ in its range we have for any edges $e, f \in E^1(\Ha\Delta)$
    \begin{align*}
        \psi(\hat{s}_e)^\ast \psi(\hat{s}_f)
            &= s_e^\ast s_f \\
            &= \begin{cases}
                \delta_{ef} \sum_{v \in r_{\Ha\Gamma}(e)} p_v,   &r_{\Ha\Gamma}(e) \neq \emptyset, \\
                \delta_{ef} s_e,                  &\text{otherwise},
            \end{cases} \\
            &= \begin{cases}
                \delta_{ef} \psi\left(\sum_{v \in r_{\Ha\Delta}(e)} \hat{p}_v\right),   &r_{\Ha\Delta}(e) \neq \emptyset, \\
                \delta_{ef} \psi(\hat{s}_e),                  &\text{otherwise},
            \end{cases} 
    \end{align*}
    using $\psi(\hat{p}_w) + \psi(\hat{p}_{w^\prime}) = p_w$.
        
    (HR2): For $e \in E^1(\Ha\Delta) \setminus (F \cup G)$ we have
    \begin{align*}
        \psi(\hat{s}_e) \psi(\hat{s}_e)^\ast 
            = s_e s_e^\ast 
            \leq \sum_{v \in s_{\Ha\Gamma}(e)} p_v = \psi\left(\sum_{v \in s_{\Ha\Delta}(e)} \hat{p}_v\right).
    \end{align*}
    For $f \in F$ it is 
    \begin{align*}
        \psi(\hat{s}_f) \psi(\hat{s}_f)^\ast 
            &= s_f s_f^\ast  \\
            &= \left(1 - \sum_{g \in G} s_g s_g^\ast\right) (s_f s_f^\ast) \left(1 - \sum_{g \in G} s_g s_g^\ast\right) \\
            &\leq \left(1 - \sum_{g \in G} s_g s_g^\ast\right) \left( \sum_{v \in s_{\Ha\Gamma}(f)} p_v \right) \left(1 - \sum_{g \in G} s_g s_g^\ast\right) \\
            &= \sum_{v \in s_{\Ha\Gamma}(f) \setminus \{w\}} p_v + \left(1 - \sum_{g \in G} s_g s_g^\ast\right) p_w \left(1 - \sum_{g \in G} s_g s_g^\ast\right) \\
            &= \sum_{v \in s_{\Ha\Gamma}(f) \setminus \{w\}} p_v + \left(\sum_{f \in F} s_f s_f^\ast\right) p_w \left(\sum_{f \in F} s_f s_f^\ast\right) \\
            &= \psi\left(\sum_{v \in s_{\Ha\Delta}(e)} \hat{p}_v\right).
    \end{align*}
    Similarly, one sees $\psi(\hat{s}_g) \psi(\hat{s}_g)^\ast \leq \psi\left(\sum_{v \in s_{\Ha\Delta}(g)} \hat{p}_v\right)$ for all $g \in G$.
    
    (HR3): Let $v \in E^0(\Ha\Delta) \setminus \{w, w^\prime\}$ not be a sink in $\Ha\Delta$. Then $v$ is not a sink in $\Ha\Gamma$ neither and therefore
    \begin{align*}
        \psi(\hat{p}_v) 
            = p_v 
            \leq \sum_{e \in E^1(\Ha\Gamma): v \in s_{\Ha\Gamma}(e)} s_e s_e^\ast 
            = \sum_{e \in E^1(\Ha\Delta): v \in s_{\Ha\Delta}(e)} \psi(\hat{s}_e) \psi(\hat{s}_e)^\ast.
    \end{align*}
    For the vertex $w$ we have 
    \begin{align*}
        \psi(\hat{p}_w)
            &= \left(\sum_{g \in G} s_g s_g^\ast\right) p_w \left(\sum_{g \in G} s_g s_g^\ast\right) \\
            &\leq \left(\sum_{g \in G} s_g s_g^\ast\right) \left( \sum_{e \in E^1(\Ha\Gamma): w \in s_{\Ha\Gamma}(e)} s_e s_e^\ast \right) \left(\sum_{g \in G} s_g s_g^\ast\right) \\
            &= \left(\sum_{g \in G} s_g s_g^\ast\right) \left( \sum_{e \in F \cup G} s_e s_e^\ast \right) \left(\sum_{g \in G} s_g s_g^\ast\right) \\
            &= \sum_{g \in G} s_g s_g^\ast \\
            &= \sum_{e \in E^1(\Ha\Delta): w \in s_{\Ha\Delta}(e)} \psi(\hat{s}_e) \psi(\hat{s}_e)^\ast.
    \end{align*}
    Similarly, one obtains $\psi(\hat{p}_{w^\prime}) \leq \sum_{e \in E^1(\Ha\Delta): w^\prime \in s_{\Ha\Delta}(e)} \psi(\hat{s}_e) \psi(\hat{s}_e)^\ast$.

    It is not hard to check that $\varphi$ and $\psi$ are inverse to each other. Therefore, we obtain $C^\ast(\Ha\Gamma) = C^\ast(\Ha\Delta)$.
\end{proof}

\begin{remark} \label{min_rem_separate_source_of_e}
    We will say that $\Ha\Delta$ is obtained from $\Ha\Gamma$ by separating the source of an edge $f \in E^1(\Ha\Gamma)$ if $\Ha\Delta$ is obtained from $\Ha\Gamma$ by applying successively source separation on the non-empty set $\{f\} \subset E^1(\Ha\Gamma)$ at all vertices $w \in s_{\Ha\Gamma}(f)$ for which there is another edge $e \in E^1(\Ha\Gamma) \setminus \{f\}$ with $w \in s_{\Ha\Gamma}(e)$.
\end{remark}

\subsection{Edge Cutting}

On the $C^\ast$-algebra side, edge cutting corresponds to taking a subalgebra.

\begin{proposition}[edge cutting] \label{min_edge_cut}
    Assume that $\Ha\Delta$ is obtained from $\Ha\Gamma$ by cutting an edge $f$. Then $C^\ast(\Ha\Delta)$ is a subalgebra of $C^\ast(\Ha\Gamma)$.
\end{proposition}

\begin{proof}
    Use the universal property of $C^\ast(\Ha\Delta)$ to obtain $\varphi: C^\ast(\Ha\Delta) \to C^\ast(\Ha\Gamma)$ with
    \begin{align*}
        \varphi: 
        \left\{
        \begin{aligned}
            p_v &\mapsto \hat{p}_v,            && v \in E^0(\Ha\Delta), \\
            s_e &\mapsto \hat{s}_e,            && e \in E^1(\Ha\Delta) \setminus \{f\}, \\
            s_e &\mapsto \hat{s}_f \hat{s}_f^\ast,&& e = f,
        \end{aligned}
        \right.
    \end{align*}
    where we write again $\hat{p}_v$ and $\hat{s}_e$ for the generators of $C^\ast(\Ha\Gamma)$ to avoid confusion. The hypergraph relations are checked by routine calculations.

    We need to show that $\varphi$ is injective. Let $\rho$ be the universal representation of $C^\ast(\Ha\Delta)$ on a Hilbert space $\mathcal{H}$ given by the GNS-construction. Further, let $\kappa$ be a cardinal larger than the dimension of $\mathcal{H}$ and let $\sigma: C^\ast(\Ha\Delta) \to B(\mathcal{H}^\kappa)$ be $\kappa$ times the representation $\rho$. Then $\sigma(s_f) \mathcal{H}^\kappa$ and $\sigma(\sum_{v \in r_{\Ha\Delta}(f)} p_v) \mathcal{H}^\kappa$ have the same dimension and therefore $B(\mathcal{H}^\kappa)$ contains a partial isometry $V$ with $V V^\ast = \sigma(s_f)$ and $V^\ast V = \sigma(\sum_{v \in r_{\Ha\Delta}(f)} p_v)$.
    Now, one readily checks that the universal property of $C^\ast(\Ha\Gamma)$ yields a representation $\tau$ of $C^\ast(\Ha\Gamma)$ on $\mathcal{H}^\kappa$ with 
    \begin{align*}
        \tau: 
        \left\{
        \begin{aligned}
            \hat{p}_v &\mapsto \sigma(p_v),            && v \in E^0(\Ha\Gamma), \\
            \hat{s}_e &\mapsto \sigma(s_e),            && e \in E^1(\Ha\Gamma) \setminus \{f\}, \\
            \hat{s}_e &\mapsto V,                      && e = f,
        \end{aligned}
        \right.
    \end{align*}
    Evidently, $\sigma = \tau \circ \varphi$, so $x \in \mathrm{ker}(\varphi) \implies x \in \mathrm{ker}(\sigma) = \mathrm{ker}(\rho)$. Since $\rho$ is faithful, the latter entails $x = 0$. Thus, $\varphi$ is injective.
\end{proof}

\subsection{Edge and Vertex Deletion}

On the $C^\ast$-algebra side, edge and vertex deletion correspond to taking the quotient of a subalgebra. Under special conditions this improves to taking the quotient, or to isomorphy of the associated $C^\ast$-algebras.

To prove this, we need some more general considerations. 
Let $\Ha\Gamma$ be hypergraph. Given a set $S \subset E^0(\Ha\Gamma) \cup E^1(\Ha\Gamma)$, we say that $\Ha\Delta$ is obtained from $\Ha\Gamma$ by deleting the set $S$, if one gets $\Ha\Delta$ from $\Ha\Gamma$ be successively deleting the edges and vertices in $\Ha\Gamma$. Clearly, the order in which these operations are performed is irrelevant.

Sometimes deleting a set $S$ of vertices and edges successively, turns out to behave better than an arbitrary single deletion. This is the case if the set $S$ has the property given by the next definition.

\begin{definition}[ideally closed set]  \label{min_ideally_closed_def}
    Let $S \subset E^0(\Ha\Gamma) \cup E^1(\Ha\Gamma)$ be a subset of the edges and vertices of a hypergraph $\Ha\Gamma$. The set $S$ is called 
    \emph{ideally closed}\index{ideally closed} if
    \begin{itemize}
        \item whenever an edge $e$ is in $S$, then $r(e) \subset S$,
        \item whenever an edge $e \in E^1(\Ha\Gamma)$ satisfies $s(e) \subset S$ or $\emptyset \neq r(e) \subset S$, then $e \in S$,
        \item whenever a vertex $v \in E^0(\Ha\Gamma)$ is not a sink and satisfies 
        $$ v \in s(e) \implies e \in S \quad \text{ for all edges } e \in E^1(\Ha\Gamma),$$ 
        then $v \in S$.
    \end{itemize}
\end{definition}


\begin{lemma} \label{min_lemma_delete_ideally_closed_set}
    Assume that $\Ha\Delta$ is obtained from $\Ha\Gamma$ by deleting an ideally closed set $S \subset E^0(\Ha\Gamma) \cup E^1(\Ha\Gamma)$. Then $C^\ast(\Ha\Delta)$ is isomorphic to the quotient $C^\ast(\Ha\Gamma) / (S)$, where $(S) \subset C^\ast(\Ha\Gamma)$ is the ideal generated by 
    $$ \{p_v, s_e: v \in S \cap E^0(\Ha\Gamma), e \in S \cap E^1(\Ha\Gamma)\}. $$
\end{lemma}

\begin{proof}
    Use the universal property of $C^\ast(\Ha\Gamma)$ to obtain $\varphi: C^\ast(\Ha\Gamma) \to C^\ast(\Ha\Delta)$ with
    \begin{align*}
        \varphi: 
        \left\{
        \begin{aligned}
            p_v &\mapsto \hat{p}_v,            && v \in E^0(\Ha\Gamma) \setminus S, \\
            p_v &\mapsto 0,                  && v \in E^0(\Ha\Gamma) \cap S,   \\
            s_e &\mapsto \hat{s}_e,            && e \in E^1(\Ha\Gamma) \setminus S, \\
            s_e &\mapsto 0,                  && e \in E^1(\Ha\Gamma) \cap S,
        \end{aligned}
        \right.
    \end{align*}
    where we write $\hat{p}_v$ and $\hat{s}_e$ for the generators of $C^\ast(\Ha\Delta)$ to avoid confusion. Using that the set $S$ is ideally closed one checks the hypergraph relations. 

    
    Clearly the $\varphi(p_v)$ are pairwise orthogonal projections and the $\varphi(s_e)$ are partial isometries.

    (HR1): Let $e, f \in E^1(\Ha\Gamma) \setminus S$. If $e$ has nonempty range in $\Ha\Delta$, then the same holds in $\Ha\Gamma$ since otherwise the edge $e$ would have been deleted, too. Thus, we have 
    \begin{align*}
        \varphi(s_e)^\ast \varphi(s_f) 
            &= \hat{s}_e^\ast \hat{s}_f  \\
            &= \begin{cases}
                \delta_{ef} \sum_{v \in r_{\Ha\Delta}(e)} \hat{p}_v,             &r_{\Ha\Delta}(e) \neq \emptyset, \\
                \delta_{ef} \hat{s}_e,                            &\text{otherwise},
            \end{cases} \\
            &= \begin{cases}
                \delta_{ef} \varphi\left(\sum_{v \in r_{\Ha\Gamma}(e)} p_v\right),             &r_{\Ha\Gamma}(e) \neq \emptyset, \\
                \delta_{ef} \varphi(s_e),                            &\text{otherwise}.
            \end{cases}
    \end{align*}
    For $e \in E^1(\Ha\Gamma) \cap S$ and $f \in E^1(\Ha\Gamma)$ it is
    \begin{align*}
        \varphi(s_e)^\ast \varphi(s_f) = 0 = \begin{cases}
            \delta_{ef} \varphi\left(\sum_{v \in r_{\Ha\Gamma}(e)} p_v\right),  &r_{\Ha\Gamma}(e) \neq \emptyset, \\
            \delta_{ef} \varphi(s_e),                 &\text{otherwise},
        \end{cases}
    \end{align*}
    since together with $e$ all vertices in $r_{\Ha\Gamma}(e)$ are deleted. By taking adjoints one obtains the same equality for $e \in E^1(\Ha\Gamma)$ and $f \in E^1(\Ha\Gamma) \cap S$.
    
    (HR2): For $e \in E^1(\Ha\Gamma) \setminus S$ observe
        $$ \varphi(s_e) \varphi(s_e)^\ast = \hat{s}_e \hat{s}_e^\ast \leq \sum_{v \in s_{\Ha\Delta}(e)} \hat{p}_v = \varphi\left(\sum_{v \in s_{\Ha\Gamma}(e)} p_v\right). $$
    If $e \in E^1(\Ha\Gamma) \cap S$, then $\varphi(s_e) \varphi(s_e)^\ast = 0 \leq \varphi\left(\sum_{v \in s_{\Ha\Gamma}(e)} p_v \right)$ is trivial.
    
    (HR3): If $v \in E^0(\Ha\Gamma) \setminus S$ were a sink in $\Ha\Delta$ but not in $\Ha\Gamma$ then it would have been deleted since $S$ is ideally closed. Thus, we have for every $v \in E^0(\Ha\Gamma) \setminus S$ that is not a sink in $\Ha\Gamma$
    \begin{align*}
        \varphi(p_v) = \hat{p}_v 
            &\leq \sum_{e \in E^1(\Ha\Delta): v \in s_{\Ha\Delta}(e)} \hat{s}_e \hat{s}_e^\ast \\
            &= \sum_{e \in E^1(\Ha\Delta): v \in s_{\Ha\Delta}(e)} \varphi(s_e) \varphi(s_e)^\ast \\
            &= \sum_{e \in E^1(\Ha\Gamma): v \in s_{\Ha\Gamma}(e)} \varphi(s_e) \varphi(s_e)^\ast.
    \end{align*}
    If $v \in E^0(\Ha\Gamma) \cap S$, then $ \varphi(p_v) = 0 \leq \sum_{e \in E^1(\Ha\Gamma): v \in s_{\Ha\Gamma}(e)} \varphi(s_e) \varphi(s_e)^\ast$ is trivial.

    Evidently, the image of $\varphi$ contains all generators of $C^\ast(\Ha\Delta)$. Hence, it is $C^\ast(\Ha\Delta) = C^\ast(\Ha\Gamma) / \mathrm{ker}(\varphi)$. 
    

    Moreover, it turns out that $C^\ast(\Ha\Delta)$ satisfies the universal property of the quotient $C^\ast(\Ha\Gamma) / (S)$ with the quotient map given by $\varphi$.
    Indeed, let $A$ be any $C^\ast$-algebra and $\chi: C^\ast(\Ha\Gamma) \to A$ a $\ast$-homomorphism with $S \subset \mathrm{ker}(\chi)$. The universal property of $C^\ast(\Ha\Delta)$ yields a map $\psi: C^\ast(\Ha\Delta) \to A$ with
    \begin{align*}
        \psi: 
        \left\{
        \begin{aligned}
            \hat{p}_v &\mapsto \chi(p_v),            && v \in E^0(\Ha\Delta), \\
            \hat{s}_e &\mapsto \chi(s_e),            && e \in E^1(\Ha\Delta).
        \end{aligned}
        \right.
    \end{align*}
    Evidently, the $\psi(\hat{p}_v)$ are pairwise orthogonal projections and the $\psi(\hat{s}_e)$ are partial isometries. We check the hypergraph relations.

    (HR1): For $e, f \in E^1(\Ha\Delta)$ we have 
    \begin{align*}
        \psi&(\hat{s}_e)^\ast \psi(\hat{s}_f)
            = \chi(s_e)^\ast \chi(s_f) \\
            &= \chi(s_e^\ast s_f) \\
            &= \begin{cases}
                \chi(\delta_{ef} \sum_{v \in r_{\Ha\Gamma}(e)} p_v) 
                    = \delta_{ef} \sum_{v \in E^0(\Ha\Gamma) \cap S \cap r_{\Ha\Gamma}(e)} \chi(p_v), &r_{\Ha\Gamma}(e) \neq \emptyset, \\
                \chi( \delta_{ef} s_e), &\text{otherwise},
            \end{cases} \\
            &= \begin{cases}
                \delta_{ef} \psi(\sum_{v \in r_{\Ha\Delta}(e)} \hat{p}_v), &r_{\Ha\Delta}(e) \neq \emptyset, \\
                 \delta_{ef} \psi(\hat{s}_e),                &\text{otherwise},
            \end{cases}
    \end{align*}
    using $\chi(p_v) = 0$ for $v \in E^0(\Ha\Gamma) \cap S$.

    (HR2): For $e \in E^1(\Ha\Delta)$ observe 
    \begin{align*}
        \psi(\hat{s}_e) \psi(\hat{s}_e)^\ast 
            &= \chi(s_e) \chi(s_e)^\ast  \\
            &= \chi(s_e s_e^\ast) \\
            &\leq \chi(\sum_{v \in s_{\Ha\Gamma}(e)} p_v) \\
            &= \sum_{v \in E^0(\Ha\Gamma) \cap S \cap s_{\Ha\Gamma}(e)} \chi( p_v) \\
            &= \psi(\sum_{v \in s_{\Ha\Delta}(e)} \hat{p}_v).
    \end{align*}
    (HR3): If $v \in E^0(\Ha\Delta)$ is not a sink, then $v$ is not a sink in $\Ha\Gamma$ neither, and we have
    \begin{align*}
        \psi(\hat{p}_v)
            &= \chi(p_v) \\
            &\leq \chi\left( \sum_{e \in E^1(\Ha\Gamma): v \in s_{\Ha\Gamma}(e)} s_e s_e^\ast \right) \\
            &= \sum_{e \in E^1(\Ha\Gamma) \cap S: v \in s_{\Ha\Gamma}(e)} \chi\left( s_e s_e^\ast \right) \\
            &= \sum_{e \in E^1(\Ha\Delta): v \in s_{\Ha\Delta}(e)} \psi(\hat{s}_e) \psi(\hat{s}_e)^\ast.
    \end{align*}

    One readily checks $\chi = \psi \circ \varphi$. Clearly, $\psi$ is the unique map from $C^\ast(\Ha\Delta)$ into $A$ with this property. It follows $C^\ast(\Ha\Delta) = C^\ast(\Ha\Gamma) / (S)$ as desired.
\end{proof}

The next preparatory lemma allows removing a vertex from the source of an edge without changing the associated hypergraph $C^\ast$-algebra.

\begin{lemma} \label{min_edge_deletion_trivial}
    Let $\Ha\Gamma$ be a hypergraph and let $w \in E^0(\Ha\Gamma)$, $f \in E^1(\Ha\Gamma)$ with $w \in s_{\Ha\Gamma}(f)$ and $r_{\Ha\Gamma}(f) = \emptyset$. Assume that $f$ is the only edge starting from $w$, i.e. for all $e \in E^1(\Ha\Gamma) \setminus \{f\}$ it is
    $
        w \not \in s_{\Ha\Gamma}(e)
    $,
    and obtain $\Ha\Delta$ from $\Ha\Gamma$ by removing $w$ from the source of $f$, possibly deleting the edge $f$ if $\{w\} = s_{\Ha\Gamma}(f)$, i.e.
    \begin{itemize}
        \item $E^0(\Ha\Delta) = E^0(\Ha\Gamma)$,
        \item $E^1(\Ha\Delta) = \begin{cases}
            E^1(\Ha\Gamma) \setminus \{f\},     &\text{if } \{w\} = s_{\Ha\Gamma}(f), \\
            E^1(\Ha\Gamma),                     &\text{otherwise},
        \end{cases}$
        \item $s_{\Ha\Delta}(e) = \begin{cases}
            s_{\Ha\Gamma}(e),                   &e \neq f, \\
            s_{\Ha\Gamma}(e) \setminus \{w\},   &e = f,
        \end{cases}
        \quad$ for all $e \in E^1(\Ha\Delta)$,
        \item $r_{\Ha\Delta}(e) = r_{\Ha\Gamma}(e)$ for all $e \in E^1(\Ha\Delta)$.
    \end{itemize}
    Then $C^\ast(\Ha\Delta) = C^\ast(\Ha\Gamma)$. In particular, if for all edges $e \in E^1(\Ha\Gamma) \setminus \{f\}$ we have $s_{\Ha\Gamma}(e) \cap s_{\Ha\Gamma}(f) = \emptyset$, then we can delete the edge $f$ from $\Ha\Gamma$ without changing the associated $C^\ast$-algebra.
\end{lemma}

\begin{proof}
    If $\{w\} = s_{\Ha\Gamma}(f)$, then set $\hat{s}_f = 0 \in C^\ast(\Ha\Delta)$.
    Using the universal property of $C^\ast(\Ha\Gamma)$ and $C^\ast(\Ha\Delta)$, respectively, one obtains maps $\varphi: C^\ast(\Ha\Gamma) \to C^\ast(\Ha\Delta)$ and $\psi: C^\ast(\Ha\Delta) \to C^\ast(\Ha\Gamma)$ with 
    \begin{align*}
        \varphi: 
        \left\{
        \begin{aligned}
            p_v &\mapsto \hat{p}_v,              && v \in E^0(\Ha\Gamma), \\
            s_e &\mapsto \hat{s}_e,              && e \in E^1(\Ha\Gamma) \setminus \{f\}, \\
            s_e &\mapsto \hat{s}_f + \hat{p}_w,    && e=f,
        \end{aligned}
        \right.
    \end{align*}
    and
    \begin{align*}
        \psi: 
        \left\{
        \begin{aligned}
            \hat{p}_v &\mapsto p_v,            && v \in E^0(\Ha\Delta), \\
            \hat{s}_e &\mapsto s_e,            && e \in E^1(\Ha\Delta) \setminus \{f\}, \\
            \hat{s}_e &\mapsto s_f - p_w,        && e = f.
        \end{aligned}
        \right.
    \end{align*}
    As usually we write $\hat{p}_v$ and $\hat{s}_e$ for the generators of $C^\ast(\Ha\Delta)$ to avoid confusion with the elements of $C^\ast(\Ha\Gamma)$.
    For both $\varphi$ and $\psi$ the hypergraph relations are checked by routine calculations. Moreover, one easily checks that $\varphi$ and $\psi$ are inverse to each other. Thus, $C^\ast(\Ha\Delta) = C^\ast(\Ha\Gamma)$.

    Finally, if for all edges $e \in E^1(\Ha\Gamma) \setminus \{f\}$ it is $s_{\Ha\Gamma}(e) \cap s_{\Ha\Gamma}(f) = \emptyset$, then we can use the previous result to successively remove every vertex $w \in s_{\Ha\Gamma}(f)$ from the source of $f$ without changing the associated $C^\ast$-algebra. In the end, this deletes the edge $f$.
\end{proof}

\begin{proposition}[edge/vertex deletion] \label{min_deletion}
    Assume that $\Ha\Delta$ is obtained from $\Ha\Gamma$ by 
    \begin{enumerate}
        \item deleting an edge $f$, or
        \item deleting a vertex $w$.
    \end{enumerate}
    Then $C^\ast(\Ha\Delta)$ is the quotient of a subalgebra of $C^\ast(\Ha\Gamma)$.
\end{proposition}

\begin{proof}
    Ad (1): First obtain $\Ha\Gamma^\prime$ from $\Ha\Gamma$ by cutting the edge $f$. Then we have $C^\ast(\Ha\Gamma^\prime) \subset C^\ast(\Ha\Gamma)$ by Proposition \ref{min_edge_cut}. Let $v_1, \dots, v_n \in E^1(\Ha\Gamma^\prime)$ be those vertices in $s_{\Ha\Gamma^\prime}(f)$ which have only $f$ as outgoing edge and obtain $\Ha\Gamma^{\prime\prime}$ from $\Ha\Gamma^\prime$ by removing the $v_i$ from the source of $f$ and leaving everything else invariant. By Lemma \ref{min_edge_deletion_trivial} this does not change the associated $C^\ast$-algebra, i.e. we have $C^\ast(\Ha\Gamma^{\prime\prime}) = C^\ast(\Ha\Gamma^\prime)$.
    If $f$ has been deleted in the process we are done. Otherwise, observe that
    in $\Ha\Gamma^{\prime\prime}$ the set $\{f\} \subset E^0(\Ha\Gamma^{\prime\prime}) \cup E^1(\Ha\Gamma^{\prime\prime})$ is ideally closed. Since $\Ha\Delta$ is obtained from $\Ha\Gamma^{\prime\prime}$ by deleting the edge $f$, it follows from Lemma \ref{min_lemma_delete_ideally_closed_set} that $C^\ast(\Ha\Delta)$ is a quotient of $C^\ast(\Ha\Gamma^{\prime\prime}) \subset C^\ast(\Ha\Gamma)$.

    \noindent Ad (2):
    First, obtain $\Ha\Gamma^\prime$ from $\Ha\Gamma$ by cutting all edges $e$ with $s_{\Ha\Gamma}(e) = \{v\}$ or $r_{\Ha\Gamma}(e) = \{v\}$. Then $C^\ast(\Ha\Gamma^\prime)$ is a subalgebra of $C^\ast(\Ha\Gamma)$ by Lemma \ref{min_edge_cut}. One readily verifies that $\Ha\Delta$ is obtained from $\Ha\Gamma^\prime$ by deleting the vertex $v$ together with all edges $e$ with $s_{\Ha\Gamma^\prime}(e) = \{v\}$. Fortunately, the set $S = \{v\} \cup \{e \in E^1(\Ha\Gamma^\prime): s_{\Ha\Gamma^\prime}(e) = \{v\}\}$ is ideally closed in $\Ha\Gamma^\prime$ and therefore $C^\ast(\Ha\Delta)$ is a quotient of $C^\ast(\Ha\Gamma^\prime) \subset C^\ast(\Ha\Gamma)$ by Lemma \ref{min_lemma_delete_ideally_closed_set}.
\end{proof}

\subsection{Edge Contraction} \label{min_edge_contraction_subsec}

On the $C^\ast$-algebra side, forward and backward edge contraction do not change the $C^\ast$-algebra up to Morita equivalence.

To prepare for the proof, the following two lemmas allow rearranging certain edges of a hypergraph without changing the associated $C^\ast$-algebra. 

\begin{lemma} \label{min_lemma_forward_contraction}
    Assume that $\Ha\Gamma$ is a hypergraph 
    and $f \in E^1(\Ha\Gamma)$ an edge with nonempty range such that 
    \begin{itemize}
        \item $r_{\Ha\Gamma}(f) \cap s_{\Ha\Gamma}(f) = \emptyset$,
        \item $r_{\Ha\Gamma}(e) \cap s_{\Ha\Gamma}(f) \neq \emptyset \implies r_{\Ha\Gamma}(e) = s_{\Ha\Gamma}(f)$ for all $e \in E^1(\Ha\Gamma) \setminus \{f\}$ and
        \item $f$ is the only edge starting from $s_{\Ha\Gamma}(f)$, i.e. $s_{\Ha\Gamma}(e) \cap s_{\Ha\Gamma}(f) = \emptyset$ for all $e \in E^1(\Ha\Gamma) \setminus \{f\}$.
    \end{itemize}
    Further, let $\Ha\Delta$ be given by
    \begin{itemize}
        \item $E^0(\Ha\Delta) = E^0(\Ha\Gamma)$,
        \item $E^1(\Ha\Delta) = E^1(\Ha\Gamma)$,
        \item $r_{\Ha\Delta}(e) = \begin{cases}
            r_{\Ha\Gamma}(e),           &r_{\Ha\Gamma}(e) \neq s_{\Ha\Gamma}(f), \\
            r_{\Ha\Gamma}(f),           &\text{otherwise},
        \end{cases}$
        \item $s_{\Ha\Delta}(e) = s_{\Ha\Gamma}(e)$,
    \end{itemize}
    for all edges $e \in E^1(\Ha\Delta)$.
    Then $C^\ast(\Ha\Delta) = C^\ast(\Ha\Gamma)$.
\end{lemma}

Intuitively, $\Ha\Delta$ is obtained from $\Ha\Gamma$ by changing the range of every edge $e$ with $r_{\Ha\Gamma}(e) = s_{\Ha\Gamma}(f)$ to $r_{\Ha\Gamma}(f)$.

\begin{proof}
    We use induction over the number $n$ of edges $g \in E^1(\Ha\Gamma)$ which satisfy $r_{\Ha\Gamma}(g) = s_{\Ha\Gamma}(f)$. If $n = 0$, then $\Ha\Delta = \Ha\Gamma$ and there is nothing to do. For the induction step let $g$ be an edge with $r_{\Ha\Gamma}(g) = s_{\Ha\Gamma}(f)$ and obtain $\Ha\Gamma^\prime$ from $\Ha\Gamma$ by changing the range of $g$ to $r_{\Ha\Gamma}(f)$ and leaving everything else invariant.
    Then the universal property of $C^\ast(\Ha\Gamma)$ yields a map $\varphi: C^\ast(\Ha\Gamma) \to C^\ast(\Ha\Gamma^\prime)$ with
    \begin{align*}
        \varphi: 
        \left\{
        \begin{aligned}
            p_v &\mapsto \hat{p}_v,                         && v \in E^0(\Ha\Gamma), \\
            s_e &\mapsto \hat{s}_e,                         && e \in E^1(\Ha\Gamma) \setminus \{g\}, \\
            s_e &\mapsto \hat{s}_g \hat{s}_f^\ast,            && e = g,
        \end{aligned}
        \right.
    \end{align*}
    where we write $\hat{p}_v$ and $\hat{s}_e$ for the generators of $C^\ast(\Ha\Gamma^\prime)$ to avoid confusion. Evidently, the $\varphi(p_v)$ are pairwise orthogonal projections and since
    \begin{align*}
        \varphi(s_g) \varphi(s_g)^\ast \varphi(s_g)
            = \hat{s}_g \, (\hat{s}_f^\ast \hat{s}_f) \, \hat{s}_g^\ast \hat{s}_g \hat{s}_f^\ast
            = \hat{s}_g \, (\hat{s}_g^\ast \hat{s}_g) \, \hat{s}_g^\ast \hat{s}_g \hat{s}_f^\ast
            = \hat{s}_g \hat{s}_f^\ast = \varphi(s_g),
    \end{align*}
    the $\varphi(s_e)$ are partial isometries.
    Further, we have
    \begin{align*}
        \varphi(s_g)^\ast \varphi(s_g)
            = \hat{s}_f \hat{s}_g^\ast \hat{s}_g \hat{s}_f^\ast
            = \hat{s}_f \hat{s}_f^\ast \hat{s}_f \hat{s}_f^\ast
            = \hat{s}_f \hat{s}_f^\ast
            = \sum_{v \in s_{\Ha\Gamma^\prime}(f)} \hat{p}_v
            = \varphi\left(\sum_{v \in r_{\Ha\Gamma}(g)} p_v \right)
    \end{align*}
    and
    \begin{align*}
        \varphi(s_g) \varphi(s_g)^\ast
            = \hat{s}_g \hat{s}_f^\ast \hat{s}_f \hat{s}_g^\ast
            = \hat{s}_g \hat{s}_g^\ast \hat{s}_g \hat{s}_g^\ast
            = \hat{s}_g \hat{s}_g^\ast
            \leq \sum_{v \in s_{\Ha\Gamma^\prime}(g)} \hat{p}_v
            = \varphi\left(\sum_{v \in s_{\Ha\Gamma}(g)} p_v \right).
    \end{align*}
    The other hypergraph relations are checked by routine calculations.
    At the same time, one obtains a map $\psi: C^\ast(\Ha\Gamma^\prime) \to C^\ast(\Ha\Gamma)$ with
    \begin{align*}
        \psi: 
        \left\{
        \begin{aligned}
            \hat{p}_v &\mapsto p_v,                         && v \in E^0(\Ha\Gamma^\prime), \\
            \hat{s}_e &\mapsto s_e,                         && e \in E^1(\Ha\Gamma^\prime) \setminus \{g\} \\
            \hat{s}_e &\mapsto s_g s_f,                       && e = g.
        \end{aligned}
        \right.
    \end{align*}
    Indeed, $\psi(\hat{s}_g)$ is a partial isometry since
    \begin{align*}
        \psi(\hat{s}_g) \psi(\hat{s}_g)^\ast \psi(\hat{s}_g)
            = s_g (s_f s_f^\ast) s_g^\ast s_g s_f 
            = s_g (s_g^\ast s_g) s_g^\ast s_g s_f 
            = s_g s_f 
            = \psi(\hat{s}_g).
    \end{align*}
    Moreover,
    \begin{align*}
        \psi(\hat{s}_g)^\ast \psi(\hat{s}_g)
            = s_f^\ast (s_g^\ast s_g) s_f 
            = s_f^\ast (s_f s_f^\ast) s_f 
            = s_f^\ast s_f 
            = \sum_{v \in r_{\Ha\Gamma}(f)} p_v
            = \psi\left(\sum_{v \in r_{\Ha\Gamma^\prime}(g)} \hat{p}_v\right)
    \end{align*}
    and
    \begin{align*}
        \psi(\hat{s}_g) \psi(\hat{s}_g)^\ast 
            = s_g (s_f s_f^\ast) s_g^\ast 
            = s_g (s_g^\ast s_g) s_g^\ast 
            = s_g s_g^\ast 
            \leq \sum_{v \in s_{\Ha\Gamma}(g)} p_v
            = \psi\left(\sum_{v \in s_{\Ha\Gamma^\prime}(g)} \hat{p}_v\right).
    \end{align*}
    Again the other hypergraph relations are checked by routine calculations.
    As
    \begin{align*}
        \varphi(\psi(\hat{s}_g)) = \varphi(s_g s_f) = \hat{s}_g \hat{s}_f^\ast \hat{s}_f = \hat{s}_g \hat{s}_g^\ast \hat{s}_g = \hat{s}_g
    \end{align*}
    and
    \begin{align*}
        \psi(\varphi(s_g)) = \psi(\hat{s}_g \hat{s}_f^\ast) = s_g s_f s_f^\ast = s_g s_g^\ast s_g = s_g
    \end{align*}
    the maps $\varphi$ and $\psi$ are inverse to each other. Thus, $C^\ast(\Ha\Gamma) = C^\ast(\Ha\Gamma^\prime)$ and we may apply the induction hypothesis to obtain $C^\ast(\Ha\Gamma) = C^\ast(\Ha\Delta)$.
\end{proof}

\begin{lemma} \label{min_lemma_backward_contraction}
    Assume that $\Ha\Gamma$ is a hypergraph 
    and $f \in E^1(\Ha\Gamma)$ an edge with nonempty range such that
    \begin{itemize}
        \item $r_{\Ha\Gamma}(f) \cap s_{\Ha\Gamma}(f) = \emptyset$,
        \item $r_{\Ha\Gamma}(f) \cap s_{\Ha\Gamma}(e) \neq \emptyset \implies r_{\Ha\Gamma}(f) \subset s_{\Ha\Gamma}(e)$ for all $e \in E^1(\Ha\Gamma)$,
        \item $f$ is the only edge starting from $s_{\Ha\Gamma}(f)$, i.e. $s_{\Ha\Gamma}(e) \cap s_{\Ha\Gamma}(f) = \emptyset$ for all $e \in E^1(\Ha\Gamma) \setminus \{f\}$.
    \end{itemize}
    Further, let $\Ha\Delta$ be given by
    \begin{itemize}
        \item $E^0(\Ha\Delta) = E^0(\Ha\Gamma)$,
        \item $E^1(\Ha\Delta) = E^1(\Ha\Gamma)$,
        \item $r_{\Ha\Delta}(e) = \begin{cases}
            r_{\Ha\Gamma}(e),          &e \neq f, \\
            s_{\Ha\Gamma}(f),          &e = f,
        \end{cases}$
        \item $s_{\Ha\Delta}(e) = \begin{cases}
            s_{\Ha\Gamma}(e),           &e \neq f \wedge r_{\Ha\Gamma}(f) \not \subset s_{\Ha\Gamma}(e), \\
            r_{\Ha\Gamma}(f),           &e = f, \\
            (s_{\Ha\Gamma}(e) \setminus r_{\Ha\Gamma}(f)) \cup s_{\Ha\Gamma}(f),           &\text{otherwise},
        \end{cases}$
    \end{itemize}
    for all edges $e \in E^1(\Ha\Delta)$.
    Then $C^\ast(\Ha\Delta) = C^\ast(\Ha\Gamma)$.
\end{lemma}

Intuitively, $\Ha\Delta$ is obtained from $\Ha\Gamma$ by inverting the direction of the edge $f$ and by replacing $r_{\Ha\Gamma}(f)$ with $s_{\Ha\Gamma}(f)$ in the source of every edge $e$ different from $f$.

\begin{proof}
    Let $g_1, \dots, g_n \in E^1(\Ha\Gamma) \setminus \{f\}$ be those edges with $r_{\Ha\Gamma}(f) \subset s_{\Ha\Gamma}(g_i)$ and use the universal property of $C^\ast(\Ha\Gamma)$ to obtain a map $\varphi: C^\ast(\Ha\Gamma) \to C^\ast(\Ha\Delta)$ with 
\begin{align*}
    \varphi: 
    \left\{
    \begin{aligned}
        p_v &\mapsto \hat{p}_v,         && v \in E^0(\Ha\Gamma), \\
        s_e &\mapsto \hat{s}_e,         && e \in E^1(\Ha\Gamma) \setminus \{f, g_1, \dots, g_n\}, \\
        s_e &\mapsto \hat{s}_f^\ast,    && e = f, \\
        s_e &\mapsto (\hat{s}_f + \sum_{v \in s_{\Ha\Gamma}(g_i)} \hat{p}_v) \hat{s}_{g_i},  && e = g_i \text{ and } r_{\Ha\Gamma}(g_i) \neq \emptyset, \\
        s_e & \mapsto (\hat{s}_f + \sum_{v \in s_{\Ha\Gamma}(g_i)} \hat{p}_v) \hat{s}_{g_i} (\hat{s}_f^\ast + \sum_{v \in s_{\Ha\Gamma}(g_i)} \hat{p}_v),   && e = g_i \text{ and } r_{\Ha\Gamma}(g_i) = \emptyset,
    \end{aligned}
    \right.
\end{align*}
where we write $\hat{p}_v$ and $\hat{s}_e$ for the generators of $C^\ast(\Ha\Delta)$ to avoid confusion. 
Evidently, the $\varphi(p_v)$ are pairwise orthogonal projections. Using
\begin{align*}
    \hat{s}_f^\ast \left(\sum_{v \in s_{\Ha\Gamma}(g_i)} \hat{p}_v\right) 
    &= \hat{s}_f^\ast \hat{s}_f \hat{s}_f^\ast \left(\sum_{v \in s_{\Ha\Gamma}(g_i)} \hat{p}_v\right)  \\
    &= \hat{s}_f^\ast \left(\sum_{v \in s_{\Ha\Delta}(f)} \hat{p}_v\right) \left(\sum_{v \in s_{\Ha\Gamma}(g_i)} \hat{p}_v\right) \\
    &= \hat{s}_f^\ast \left(\sum_{v \in s_{\Ha\Delta}(f)} \hat{p}_v\right) \\
    &= \hat{s}_f^\ast \hat{s}_f \hat{s}_f^\ast = \hat{s}_f^\ast
\end{align*}
and 
\begin{align*}
        \hat{s}_f^\ast \hat{s}_{g_i}
        = \hat{s}_f^\ast \left(\sum_{v \in s_{\Ha\Delta}(f)} \hat{p}_v\right) \left(\sum_{v \in s_{\Ha\Delta}(g_i)} \hat{p}_v\right) \hat{s}_{g_i}
        = 0
\end{align*}
one obtains 
\begin{align*}
    \varphi(s_{g_i}) &\varphi(s_{g_i})^\ast \varphi(s_{g_i}) \\
        &= (\hat{s}_f + \sum_{v \in s_{\Ha\Gamma}(g_i)} \hat{p}_v) \hat{s}_{g_i} \hat{s}_{g_i}^\ast (\hat{s}_f^\ast + \sum_{v \in s_{\Ha\Gamma}(g_i)} \hat{p}_v) (\hat{s}_f + \sum_{v \in s_{\Ha\Gamma}(g_i)} \hat{p}_v) \hat{s}_{g_i} \\
        &= (\hat{s}_f + \sum_{v \in s_{\Ha\Gamma}(g_i)} \hat{p}_v) \hat{s}_{g_i} \hat{s}_{g_i}^\ast (\hat{s}_f^\ast \hat{s}_f + \sum_{v \in s_{\Ha\Gamma}(g_i)} \hat{p}_v + \hat{s}_f^\ast + \hat{s}_f) \hat{s}_{g_i} \\
        &= (\hat{s}_f + \sum_{v \in s_{\Ha\Gamma}(g_i)} \hat{p}_v) \hat{s}_{g_i} \hat{s}_{g_i}^\ast (\hat{s}_f^\ast \hat{s}_f + \sum_{v \in s_{\Ha\Gamma}(g_i)} \hat{p}_v) \hat{s}_{g_i} \\
        &= (\hat{s}_f + \sum_{v \in s_{\Ha\Gamma}(g_i)} \hat{p}_v) \hat{s}_{g_i} \hat{s}_{g_i}^\ast \left(\sum_{v \in s_{\Ha\Delta}(g_i)} \hat{p}_v\right) \left(\sum_{v \in r_{\Ha\Delta}(f)} \hat{p}_v + \sum_{v \in s_{\Ha\Gamma}(g_i)} \hat{p}_v\right) \hat{s}_{g_i} \\
        &= (\hat{s}_f + \sum_{v \in s_{\Ha\Gamma}(g_i)} \hat{p}_v) \hat{s}_{g_i} \hat{s}_{g_i}^\ast \hat{s}_{g_i} \\
        &= (\hat{s}_f + \sum_{v \in s_{\Ha\Gamma}(g_i)} \hat{p}_v) \hat{s}_{g_i} \\
        &= \varphi(s_{g_i})
\end{align*}
for edges $g_i$ with nonempty range. Otherwise, the same calculation can be used to show that $\varphi(s_{g_i})$ is a projection. Thus, the $\varphi(s_e)$ are partial isometries. We check the hypergraph relations.

        (HR1): For $e, e^\prime \not \in \{f, g_1, \dots, g_n\}$ we have 
        \begin{align*}
            \varphi(s_e^\ast s_e^\prime) = \hat{s}_e^\ast \hat{s}_e^\prime 
            = \begin{cases}
                \delta_{ee^\prime} \sum_{v \in r_{\Ha\Delta}(e)} \hat{p}_v = \delta_{ee^\prime} \varphi\left(\sum_{v \in r_{\Ha\Gamma}(e)} p_v\right), &r_{\Ha\Gamma}(e) \neq \emptyset, \\
                \delta_{ee^\prime} \hat{s}_e = \delta_{ee^\prime} \varphi(s_e),                &\text{otherwise},
            \end{cases}
        \end{align*}
        and for $e = f$ it is 
        \begin{align*}
            \varphi(s_f)^\ast \varphi(s_f) = \hat{s}_f \hat{s}_f^\ast = \sum_{v \in s_{\Ha\Delta}(f)} \hat{p}_v = \varphi\left(\sum_{v \in r_{\Ha\Gamma}(f)} p_v \right)
        \end{align*}
        since $f$ is the only edge in $\Ha\Delta$ that starts from $r_{\Ha\Gamma}(f)( = s_{\Ha\Delta}(f))$. Further, one obtains for $g_i$ with nonempty range the equality
        \begin{align*}
            \varphi(s_{g_i})^\ast \varphi(s_{g_i})
                &= \hat{s}_{g_i}^\ast (\hat{s}_f^\ast + \sum_{v \in s_{\Ha\Gamma}(g_i)} \hat{p}_v) (\hat{s}_f + \sum_{v \in s_{\Ha\Gamma}(g_i)} \hat{p}_v) \hat{s}_{g_i} \\
                &= \hat{s}_{g_i}^\ast \hat{s}_{g_i} \\
                &= \sum_{v \in r_{\Ha\Delta}(g_i)} \hat{p}_v \\
                &= \varphi\left(\sum_{v \in r_{\Ha\Gamma}(g_i)} p_v\right),
        \end{align*}
        and for $g_i$ with empty range the equality
        \begin{align*}
            \varphi(s_{g_i})^\ast \varphi(s_{g_i})
                = \varphi(s_{g_i})
        \end{align*}
        similarly as in the calculation of $\varphi(s_{g_i}) \varphi(s_{g_i})^\ast \varphi(s_{g_i})$ above. For $g_i$ with nonempty range and $e \not \in \{g_1, \dots, g_n,f\}$ observe
        \begin{align*}
            \varphi(s_e)^\ast \varphi(s_{g_i})
                &= \hat{s}_e^\ast (\hat{s}_f + \sum_{v \in s_{\Ha\Gamma}(g_i)} \hat{p}_v) \hat{s}_{g_i} \\
                &= \hat{s}_e^\ast \hat{s}_f + \hat{s}_e^\ast \left(\sum_{v \in s_{\Ha\Gamma}(g_i)} \hat{p}_v\right) \hat{s}_{g_i}  \\
                &= \hat{s}_e^\ast \hat{s}_f + \hat{s}_e^\ast \left(\sum_{v \in s_{\Ha\Delta}(e)} \hat{p}_v\right) \left(\sum_{v \in s_{\Ha\Gamma}(g_i)} \hat{p}_v\right) \left(\sum_{v \in s_{\Ha\Delta}(g_i)} \hat{p}_v\right) \hat{s}_{g_i} \\
                &= \hat{s}_e^\ast \hat{s}_f + \hat{s}_e^\ast \left( \sum_{v \in s_{\Ha\Delta}(e) \cap s_{\Ha\Gamma}(g_i) \cap s_{\Ha\Delta}(g_i)} \hat{p}_v \right) \hat{s}_{g_i} \\
                &= \hat{s}_e^\ast \hat{s}_f + \hat{s}_e^\ast \left( \sum_{v \in s_{\Ha\Delta}(e) \cap s_{\Ha\Delta}(g_i)} \hat{p}_v \right) \hat{s}_{g_i} \\
                &= \hat{s}_e^\ast \hat{s}_f + \hat{s}_e^\ast \hat{s}_{g_i} \\
                &= 0
        \end{align*}
        and for $g_j \neq g_i$ with nonempty range
        \begin{align*}
            \varphi&(s_{g_i})^\ast \varphi(s_{g_j})
                = \hat{s}_{g_i}^\ast \left(\hat{s}_f^\ast +  \sum_{v \in s_{\Ha\Gamma}(g_i)} \hat{p}_v\right) \left(\hat{s}_f + \sum_{v \in s_{\Ha\Gamma}(g_j)} \hat{p}_v\right) \hat{g}_j \\
                &= \hat{s}_{g_i}^\ast \left(\hat{s}_f^\ast \hat{s}_f + \left( \sum_{v \in s_{\Ha\Gamma}(g_i)} \hat{p}_v\right)\left( \sum_{v \in s_{\Ha\Gamma}(g_j)} \hat{p}_v\right) \right. \\
                &\qquad \quad \left. + \hat{s}_f^\ast \left( \sum_{v \in s_{\Ha\Gamma}(g_j)} \hat{p}_v\right) + \left( \sum_{v \in s_{\Ha\Gamma}(g_i)} \hat{p}_v\right) \hat{s}_f\right) \hat{g}_j \\
                &= \hat{s}_{g_i}^\ast \left(\sum_{v \in r_{\Ha\Delta}(f)} \hat{p}_v + \left(\sum_{v \in s_{\Ha\Gamma}(g_i)} \hat{p}_v\right) \left(\sum_{v \in s_{\Ha\Gamma}(g_j)} \hat{p}_v\right) + \hat{f}^\ast + \hat{s}_f\right) \hat{g}_j \\
                &= \hat{s}_{g_i}^\ast \left(\sum_{v \in r_{\Ha\Delta}(f)} \hat{p}_v + \left(\sum_{v \in s_{\Ha\Gamma}(g_i)} \hat{p}_v\right) \left(\sum_{v \in s_{\Ha\Gamma}(g_j)} \hat{p}_v\right)\right) \hat{g}_j \\
                &= \hat{s}_{g_i}^\ast \hat{s}_{g_j} \\
                &= 0.
        \end{align*}
        Similar calculations apply if $g_i$ and/or $g_j$ has empty range.
        Further, for any edge $e \not \in \{g_1, \dots, g_n, f\}$ use $s_{\Ha\Delta}(e) \cap r_{\Ha\Delta}(f) = \emptyset$ to obtain 
        \begin{align*}
            \varphi(s_e)^\ast \varphi(s_f)
                &= \hat{s}_e^\ast \hat{s}_f^\ast
                = \hat{s}_e^\ast \hat{s}_e \hat{s}_e^\ast \hat{s}_f^\ast \hat{s}_f \hat{s}_f^\ast
                = 0.
        \end{align*}
        Finally, observe that for $g_i$ with nonempty range
        \begin{align*}
            \varphi&(s_{g_i})^\ast \varphi(s_f)
                = \hat{s}_{g_i}^\ast \left(\hat{s}_f^\ast + \left(\sum_{v \in s_{\Ha\Gamma}(g_i)} \hat{p}_v\right)\right) \hat{s}_f^\ast \\
                &= \hat{s}_{g_i}^\ast \left(\hat{s}_f^\ast \hat{s}_f^\ast + \left(\sum_{v \in s_{\Ha\Gamma}(g_i)} \hat{p}_v\right) \hat{s}_f^\ast\right)  \\
                &= \hat{s}_{g_i}^\ast \left(\hat{s}_f^\ast \left(\sum_{v \in s_{\Ha\Delta}(f)} \hat{p}_v\right) \left(\sum_{v \in r_{\Ha\Delta}(f)} \hat{p}_v\right) \hat{s}_f^\ast + \left(\sum_{v \in s_{\Ha\Gamma}(g_i)} \hat{p}_v\right) \left(\sum_{v \in r_{\Ha\Delta}(f)} \hat{p}_v\right) \hat{s}_f^\ast\right) \\
        \end{align*}
        is zero,
        using $s_{\Ha\Delta}(f) \cap r_{\Ha\Delta}(f) = \emptyset$ and $s_{\Ha\Gamma}(g_i) \cap r_{\Ha\Delta}(f) = \emptyset$. Similarly, one obtains $\varphi(s_{g_i})^\ast \varphi(s_f) = 0$ for $g_i$ with empty range.

        (HR2): For $e \not \in \{f, g_1, \dots, g_n\}$ the second hypergraph relation is easily checked. For the remaining edges observe 
        \begin{align*}
            &\varphi(s_{g_i}) \varphi(s_{g_i})^\ast \\
                &= \left(\hat{s}_f + \sum_{v \in s_{\Ha\Gamma}(g_i)} \hat{p}_v\right) \hat{g}_i \hat{g}_i^\ast \left(\hat{s}_f^\ast + \sum_{v \in s_{\Ha\Gamma}(g_i)} \hat{p}_v\right) \\
                &\leq \left(\hat{s}_f + \sum_{v \in s_{\Ha\Gamma}(g_i)} \hat{p}_v\right) \left(\sum_{v \in s_{\Ha\Delta}(g_i)} \hat{p}_v\right) \left(\hat{s}_f^\ast + \sum_{v \in s_{\Ha\Gamma}(g_i)} \hat{p}_v\right) \\
                &= \left(\hat{s}_f \hat{s}_f^\ast \hat{s}_f \left(\sum_{v \in s_{\Ha\Delta}(g_i)} \hat{p}_v\right) + \left(\sum_{v \in s_{\Ha\Gamma}(g_i)} \hat{p}_v\right)\left(\sum_{v \in s_{\Ha\Delta}(g_i)} \hat{p}_v\right)\right) \, \\
                &\qquad \cdot \left(\left(\sum_{v \in s_{\Ha\Delta}(g_i)} \hat{p}_v\right)  \hat{s}_f^\ast \hat{s}_f \hat{s}_f^\ast + \left(\sum_{v \in s_{\Ha\Delta}(g_i)} \hat{p}_v\right) \left(\sum_{v \in s_{\Ha\Gamma}(g_i)} \hat{p}_v\right) \right)\\
                &= \left(\hat{s}_f + \left(\sum_{v \in s_{\Ha\Gamma}(g_i)} \hat{p}_v\right)- \hat{s}_f \hat{s}_f^\ast\right) \, \left(\hat{s}_f^\ast + \left(\sum_{v \in s_{\Ha\Gamma}(g_i)} \hat{p}_v\right)- \hat{s}_f \hat{s}_f^\ast\right) \\
                &= \hat{s}_f \hat{s}_f^\ast + \hat{s}_f \left(\sum_{v \in s_{\Ha\Gamma}(g_i)} \hat{p}_v\right)- \hat{s}_f \hat{s}_f \hat{s}_f^\ast \\
                    &\qquad + \left(\sum_{v \in s_{\Ha\Gamma}(g_i)} \hat{p}_v\right)\hat{s}_f^\ast + \left(\sum_{v \in s_{\Ha\Gamma}(g_i)} \hat{p}_v\right)- \left(\sum_{v \in s_{\Ha\Gamma}(g_i)} \hat{p}_v\right)\hat{s}_f \hat{s}_f^\ast \\
                    &\qquad - \hat{s}_f \hat{s}_f^\ast \hat{s}_f^\ast - \hat{s}_f \hat{s}_f^\ast \left(\sum_{v \in s_{\Ha\Gamma}(g_i)} \hat{p}_v\right)+ \hat{s}_f \hat{s}_f^\ast \\
                &= \sum_{v \in s_{\Ha\Gamma}(g_i)} \hat{p}_v\\
                &= \varphi\left(\sum_{v \in s_{\Ha\Gamma}(g_i)} p_v\right)
        \end{align*}
        using 
        \begin{align*}
            \hat{s}_f \hat{s}_f^\ast &= \sum_{v \in s_{\Ha\Delta}(f)} \hat{p}_v = \sum_{v \in r_{\Ha\Gamma}(f)} \hat{p}_v \leq \sum_{v \in s_{\Ha\Gamma}(g_i)} \hat{p}_v \\
            & \perp \sum_{v \in s_{\Ha\Gamma}(f)} \hat{p}_v = \sum_{v \in r_{\Ha\Delta}(f)} \hat{p}_v = \hat{s}_f^\ast \hat{s}_f \leq \sum_{v \in s_{\Ha\Delta}(g_i)} \hat{p}_v \perp \hat{s}_{\Ha\Delta}(f) = \hat{s}_f \hat{s}_f^\ast
        \end{align*}
        Note that these (in)equalities hold for $g_i$ with empty or nonempty range at the same time. Moreover, it is
        \begin{align*}
            \varphi(s_f) \varphi(s_f)^\ast 
                &= \hat{s}_f^\ast \hat{s}_f
                = \sum_{v \in r_{\Ha\Delta}(f)} \hat{p}_v
                = \sum_{v \in s_{\Ha\Gamma}(f)} \hat{p}_v
                = \varphi\left(\sum_{v \in s_{\Ha\Gamma}(f)} p_v\right).
        \end{align*}
        
        (HR3): If $v \not \in r_{\Ha\Gamma}(f) \cup s_{\Ha\Gamma}(f)$ is not a sink in $\Ha\Gamma$, then with 
        \begin{align*} 
            \hat{s}_f \hat{p}_v 
            &= \hat{s}_f \hat{s}_f^\ast \hat{s}_f \hat{p}_v \\
            &= 0 \\
            &= \hat{s}_f \left(\sum_{v \in r_{\Ha\Delta}(f)} \hat{p}_v\right) \left(\sum_{v \in s_{\Ha\Delta}(e)} \hat{p}_v\right) \hat{s}_e \hat{s}_e^\ast \\
            &= \hat{s}_f \hat{s}_f^\ast \hat{s}_f \hat{s}_e \hat{s}_e^\ast \\
            &= \hat{s}_f \hat{s}_e \hat{s}_e^\ast
        \end{align*} 
        for all $e \not \in \{g_1, \dots, g_n, f\}$ one obtains
        \begin{align*}
            \varphi(p_v)
                &= \hat{p}_v \\
                &= (\hat{s}_f + 1 - \hat{s}_f^\ast \hat{s}_f) \hat{p}_v (\hat{s}_f^\ast + 1 - \hat{s}_f^\ast \hat{s}_f) \\
                &\leq (\hat{s}_f + 1 - \hat{s}_f^\ast \hat{s}_f) \left(\sum_{e \in E^1(\Ha\Delta): v \in s_{\Ha\Delta}(e)} \hat{s}_e \hat{s}_e^\ast \right) (\hat{s}_f^\ast + 1 - \hat{s}_f^\ast \hat{s}_f) \\
                &= \sum_{e \in E^1(\Ha\Delta): e \neq g_i \wedge v \in s_{\Ha\Gamma}(e)} \hat{s}_e \hat{s}_e^\ast \\
                    &\qquad + \sum_{g_i: v \in s_{\Ha\Gamma}(g_i)} \left(\hat{s}_f + \sum_{v \in s_{\Ha\Gamma}(g_i)} \hat{p}_v\right) \hat{s}_{g_i} \hat{s}_{g_i}^\ast \left(\hat{s}_f^\ast + \sum_{v \in s_{\Ha\Gamma}(g_i)} \hat{p}_v\right) \\
                &= \sum_{e \in E^1(\Ha\Delta): v \in s_{\Ha\Gamma}(e)} \varphi(s_e) \varphi(s_e)^\ast.
        \end{align*}
        Further, if $v \in s_{\Ha\Gamma}(f)$, then $f$ is the only edge in $\Ha\Gamma$ starting from $v$ and therefore
        \begin{align*}
            \varphi(p_v)
                &= \hat{p}_v
                \leq \sum_{v \in s_{\Ha\Gamma}(f)} \hat{p}_v 
                = \sum_{v \in r_{\Ha\Delta}(f)} \hat{p}_v
                = \hat{s}_f^\ast \hat{s}_f
                = \sum_{e \in E^1(\Ha\Gamma): v \in s_{\Ha\Gamma}(e)} \varphi(s_f) \varphi(s_f)^\ast.
        \end{align*}
        Finally, let $v \in r_{\Ha\Gamma}(f)$ not be a sink in $\Ha\Gamma$. Then in $\Ha\Gamma$ the $g_i$ are all edges that start from $v$ while in $\Ha\Delta$ the $g_i$ are all edges that start from $s_{\Ha\Gamma}(f)( = r_{\Ha\Delta}(f))$. Thus, we get 
        \begin{align*}
            \varphi(p_v)
                = \hat{p}_v
                &\leq \sum_{v \in s_{\Ha\Delta}(f)} \hat{p}_v \\
                &= \hat{s}_f \hat{s}_f^\ast \\
                &= \hat{s}_f \hat{s}_f^\ast \hat{s}_f \hat{s}_f^\ast \\
                &= (\hat{s}_f + 1 - \hat{s}_f^\ast \hat{s}_f) \hat{s}_f^\ast \hat{s}_f (\hat{s}_f^\ast + 1 - \hat{s}_f^\ast \hat{s}_f) \\ 
                &\leq (\hat{s}_f + 1 - \hat{s}_f^\ast \hat{s}_f) \left( \sum_{i=1}^n \hat{s}_{g_i} \hat{s}_{g_i}^\ast \right) (\hat{s}_f^\ast + 1 - \hat{s}_f^\ast \hat{s}_f) \\
                &= \sum_{i=1}^n \left(\hat{s}_f + \sum_{v \in s_{\Ha\Gamma}(g_i)} \hat{p}_v\right) \hat{g}_i \hat{g}_i^\ast \left(\hat{s}_f^\ast + \sum_{v \in s_{\Ha\Gamma}(g_i)} \hat{p}_v\right) \\
                &= \sum_{e \in E^1(\Ha\Gamma): v \in s_{\Ha\Gamma}(e)} \varphi(s_e) \varphi(s_e)^\ast.
        \end{align*}

    It remains to find an inverse map $\psi: C^\ast(\Ha\Delta) \to C^\ast(\Ha\Gamma)$. A close inspection of $\Ha\Delta$ reveals that it satisfies the assumptions required from $\Ha\Gamma$ and indeed we get back $\Ha\Gamma$ from $\Ha\Delta$ by the same procedure that gave us $\Ha\Delta$ from $\Ha\Gamma$. Thus, by the very same arguments as above one obtains a map $\psi$ with 
    \begin{align*}
        \psi: 
        \left\{
        \begin{aligned}
            \hat{p}_v &\mapsto p_v,         && v \in E^0(\Ha\Delta), \\
            \hat{s}_e &\mapsto s_e,         && e \in E^1(\Ha\Delta) \setminus \{f, g_1, \dots, g_n\}, \\
            \hat{s}_e &\mapsto s_f^\ast,    && e = f, \\
            \hat{s}_e &\mapsto \left(s_f + \sum_{v \in s_{\Ha\Delta}(g_i)} p_v\right) s_{g_i},  && e = g_i \text{ and } r_{\Ha\Delta}(g_i) \neq \emptyset, \\
            \hat{s}_e &\mapsto \left(s_f + \sum_{v \in s_{\Ha\Delta}(g_i)} p_v\right) s_{g_i} \left(s_f^\ast + \sum_{v \in s_{\Ha\Delta}(g_i)} p_v\right),  && e = g_i \text{ and } r_{\Ha\Delta}(g_i) = \emptyset. \\
        \end{aligned}
        \right.
    \end{align*}
    Using
    \begin{align*}
        \varphi(\psi(\hat{s}_{g_i})) 
            &= \varphi\left(\left(s_f + \sum_{v \in s_{\Ha\Delta}(g_i)} p_v\right) s_{g_i}\right) \\
            &= \left(\hat{s}_f^\ast + \sum_{v \in s_{\Ha\Delta}(g_i)} \hat{p}_v\right) \left(\hat{s}_f + \sum_{v \in s_{\Ha\Gamma}(g_i)} \hat{p}_v\right) \hat{s}_{g_i} \\
            &= \left(\hat{s}_f^\ast \hat{s}_f + \hat{s}_f^\ast \left(\sum_{v \in s_{\Ha\Gamma}(g_i)} \hat{p}_v\right)\right. \\
            & \qquad \left. + \left(\sum_{v \in s_{\Ha\Delta}(g_i)} \hat{p}_v\right) \hat{s}_f + \left(\sum_{v \in s_{\Ha\Delta}(g_i)} \hat{p}_v\right) \left(\sum_{v \in s_{\Ha\Gamma}(g_i)} \hat{p}_v\right)\right) \hat{s}_{g_i} \\
            &= \left(\hat{s}_f^\ast \hat{s}_f + \hat{s}_f^\ast + 0 + \left(\sum_{v \in s_{\Ha\Delta}(g_i)} \hat{p}_v\right) - \hat{s}_f^\ast \hat{s}_f\right) \hat{s}_{g_i} \\
            &= \hat{s}_{g_i}
    \end{align*}
    for $g_i$ with nonempty range and a similar calculation for $g_i$ with empty range, one readily checks that $\varphi$ and $\psi$ are inverse to each other. This concludes the proof.
\end{proof}

\begin{proposition}[edge contraction] \label{min_edge_contraction}
    Assume that $\Ha\Delta$ is obtained from $\Ha\Gamma$ by 
    \begin{enumerate}
        \item forward contracting an edge $f$ with $s_{\Ha\Gamma}(f) = \{w\}$, or
        \item backward contracting an edge $f$ with $r_{\Ha\Gamma}(f) = \{w\}$.
    \end{enumerate}
    Then $C^\ast(\Ha\Delta) = (1 - p_w) C^\ast(\Ha\Gamma) (1 - p_w)$ and the latter is a full corner in $C^\ast(\Ha\Gamma)$. In particular, $C^\ast(\Ha\Gamma)$ and $C^\ast(\Ha\Delta)$ are Morita equivalent.
\end{proposition}

\begin{proof}
    Ad (1): Because of the range decomposition operation we may assume without loss of generality that
    \begin{align*}
        w \in r_{\Ha\Gamma}(e) \implies \{w\} = r_{\Ha\Gamma}(e) \tag{$\ast$}
    \end{align*}
    holds for all edges $e \in E^1(\Ha\Gamma)$. Otherwise, first decompose the ranges of all edges $e$ with $w \in r_{\Ha\Gamma}(e)$ and then apply forward contraction on $f$. Since $w \in r_{\Ha\Gamma}(e)$ implies $r_{\Ha\Gamma}(e) \cap r_{\Ha\Gamma}(f) = \emptyset$ by assumption, we may apply range decomposition backwards to obtain the desired hypergraph $\Ha\Delta$.

    Now, assume that ($\ast$) is true and obtain $\Ha\Gamma^\prime$ from $\Ha\Gamma$ by changing the range of every edge $e$ with $\{w\} = r_{\Ha\Gamma}(e)$ to $r_{\Ha\Gamma}(f)$. By Lemma \ref{min_lemma_forward_contraction}, $C^\ast(\Ha\Gamma^\prime) = C^\ast(\Ha\Gamma)$. It remains to show that $C^\ast(\Ha\Delta)$ is a full corner of $C^\ast(\Ha\Gamma^\prime)$.

    Observe that $\Ha\Delta$ is obtained from $\Ha\Gamma^\prime$ by first cutting the edge $f$ and then deleting the ideally closed set $\{w, f\}$. Let $\Ha\Gamma^{\prime\prime}$ be the hypergraph obtained after the first operation. Then $C^\ast(\Ha\Gamma^{\prime\prime})$ is a subalgebra of $C^\ast(\Ha\Gamma^\prime)$ by a previous proposition. As $w$ has no incoming or outgoing edges in $\Ha\Gamma^{\prime\prime}$ one has $C^\ast(\Ha\Gamma^{\prime\prime}) = C^\ast(\Ha\Delta) \oplus \C p_w$. Thus, $C^\ast(\Ha\Delta) \subset C^\ast(\Ha\Gamma^{\prime\prime}) \subset C^\ast(\Ha\Gamma^\prime) = C^\ast(\Ha\Gamma)$. A closer look at the embeddings reveals that $C^\ast(\Ha\Delta)$ is the subalgebra generated by the projections $p_v$ for $v \neq w$ and the partial isometries $s_e$ for $e \neq f$. All of these elements are in the corner $(1-p_w) C^\ast(\Ha\Gamma^\prime) (1-p_w)$, and therefore we have $C^\ast(\Ha\Delta) \subset (1-p_w) C^\ast(\Ha\Gamma^\prime) (1-p_w)$.
    
    To show equality, use that a dense subset of $C^\ast(\Ha\Gamma^\prime)$ is spanned by words of the form
    $x = x_1 \dots x_n$ with $x_i \in \{p_v, s_e, s_e^\ast: v \in E^0(\Ha\Gamma^\prime), e \in E^1(\Ha\Gamma^\prime)\}$.
    We claim that if $x$ is in $(1-p_w) C^\ast(\Ha\Gamma^\prime) (1-p_w)$, then it is also in $C^\ast(p_v, s_e: v \neq w, e \neq f)$ which is equal to $C^\ast(\Ha\Delta)$. This implies immediately
    $$ (1-p_w) C^\ast(\Ha\Gamma^\prime) (1-p_w) \subset C^\ast(\Ha\Delta). $$
    
    We prove the claim by induction over the number $N$ of occurrences of $s_f$ or $p_w$ in the word $x$. Without loss of generality, however, $x$ does not contain the letter $p_w$ since we could replace it with $s_f s_f^\ast$.
    If $N = 0$ there is nothing to show. For the induction step distinguish the following cases:

    \noindent \underline{Case 1} $x = s_f x^\prime$ for some $x^\prime$. Then $x \in (1-p_w) C^\ast(\Ha\Gamma^\prime) (1-p_w)$ implies
    \begin{align*}
        x = (1-p_w) x = (1-p_w) s_f x^\prime = (1-p_w) s_f s_f^\ast s_f x^\prime = (1-p_w) p_w s_f x^\prime = 0.
    \end{align*}

    \noindent\underline{Case 2} $x = x^\prime s_e s_f x^{\prime\prime}$ for some $x^\prime, x^{\prime\prime}$ and an edge $e \in E^1(\Ha\Gamma^\prime)$. The properties of $\Ha\Gamma^\prime$ imply that the intersection $r_{\Ha\Gamma^\prime}(e) \cap s_{\Ha\Gamma^\prime}(f)$ is empty. Therefore, $x = 0$.

    \noindent\underline{Case 3} $x = x^\prime s_e^\ast s_f x^{\prime\prime}$ for some $x^\prime, x^{\prime\prime}$ and an edge $e \in E^1(\Ha\Gamma)$. Unless $x$ is zero the intersection $s_{\Ha\Gamma^\prime}(e) \cap s_{\Ha\Gamma^\prime}(f)$ must not be empty. This leaves only the possibility $e = f$. Then
    \begin{align*}
        x = x^\prime s_e^\ast s_f x^{\prime\prime} 
            = x^\prime s_f^\ast s_f x^{\prime\prime}
            = x^\prime \sum_{v \in r_{\Ha\Gamma^\prime}(f)} p_v x^{\prime\prime}
            = \sum_{v \in r_{\Ha\Gamma^\prime}(f)} x^\prime p_v x^{\prime\prime}.
    \end{align*}
    On the last term we may apply the induction hypothesis.

    \noindent\underline{Case 4} $x = x^\prime p_v s_f x^{\prime\prime}$ for some $x^\prime, x^{\prime\prime}$ and a vertex $v \in E^0(\Ha\Gamma)$. By assumption, we have without loss of generality $v \neq w$. Thus,
    $$ x = x^\prime p_v s_f x^{\prime\prime} = 0. $$

    \noindent \underline{Case 5} $x = x^\prime s_f^\ast$, $x = x^\prime s_f^\ast s_e x^{\prime\prime}, x = x^\prime s_f^\ast s_e^\ast x^{\prime\prime}$ or $x = x^\prime s_f^\ast p_v x^{\prime\prime}$. By passing to the adjoint one of the previous cases applies.

    Ad (2): One checks that the edge $f$ satisfies the assumptions from Lemma \ref{min_lemma_backward_contraction}. Thus, we may obtain $\Ha\Gamma^\prime$ by changing the source of every edge $e$ with $w \in s_{\Ha\Gamma}(e)$ from $s_{\Ha\Gamma}(e)$ to $(s_{\Ha\Gamma}(e) \setminus \{w\}) \cup s_{\Ha\Gamma}(f)$ and invert the direction of the edge $f$ without changing the associated $C^\ast$-algebra. Then one readily checks that $\Ha\Delta$ is obtained from $\Ha\Gamma^\prime$ by forward contracting the edge $f$. Now, the claim follows from (1).
\end{proof}

\begin{proof}[Proof of Theorem \ref{min_big_thm}]
For range decomposition we refer to the proof of Theorem 4.1 in [Tri22] which is easily adapted to the
present situation. The other statements in the bullet points are identical to Propositions \ref{min_edge_sep}, \ref{min_edge_cut}, \ref{min_deletion} and \ref{min_edge_contraction}, respectively. 
The last statement follows since the class of exact $C^\ast$-algebras is closed
under taking quotients and subalgebras, and because Morita equivalence preserves exactness.
\end{proof}

\section{Hypergraph Normalization}
\label{norm_sec}


Assume that we have given a hypergraph $\Ha\Gamma$ and its associated $C^\ast$-algebra $C^\ast(\Ha\Gamma)$. In the previous Section \ref{min_sec} we have seen that certain minor operations do not change the $C^\ast$-algebra up to Morita equivalence, see Theorem \ref{min_big_thm}. 
This can be used to put any given hypergraph in a ``normalized'' form.
Let us first define the notion of a normal hypergraph.

\begin{definition}[normal hypergraph\index{normal hypergraph}] \label{mainres_nuc_normal_def}
    A hypergraph $\Ha\Gamma$ is called \emph{normal} if it has the following properties.
    \begin{enumerate}
        \item $|r(e)| \leq 1$ for all edges $e \in E^1(\Ha\Gamma)$.
        \item For every edge $e$ there exists another edge $f$ such that $s(e) \cap s(f) \neq \emptyset$ or $\emptyset \neq r(e) \subset s(e)$.
        \item Whenever $(e,f)$ is a pair of distinct edges with $|s(e) \cap s(f)| = 1$, then one of the following holds:
        \begin{enumerate}[label=\alph*)]
            \item $|s(e)|=|s(f)|=1$.
            \item There is an edge $g \neq e$ with $s(e) \cap s(f) \subsetneq s(e) \cap s(g)$.
        \end{enumerate}
    \end{enumerate}
\end{definition}

The next lemma asserts that, without changing the associated $C^\ast$-algebra up to Morita equivalence, any hypergraph $\Ha\Gamma$ can be normalized by passing to a suitable hypergraph minor.

\begin{lemma} \label{mainres_normal_hyp_thm}
    Let $\Ha\Gamma$ by a hypergraph. Then there is a normal hypergraph $\Ha\Delta \leq \Ha\Gamma$ such that $C^\ast(\Ha\Delta)$ is Morita equivalent to $C^\ast(\Ha\Gamma)$. We call $\Ha\Delta$ a \emph{normalized version}\index{normalized version} of $\Ha\Gamma$.
\end{lemma}




The idea of the proof is to use source separation, backward contraction and range decomposition as often as possible without changing the $C^\ast$-algebra up to Morita equivalence.

\begin{proof}
    Consider the sets $S_i := S_i(\Ha\Gamma)$ given by
    \begin{align*}
        S_1 &:= \{e \in E^1(\Ha\Gamma)| \, e \text{ violates condition (2) from Definition \ref{mainres_nuc_normal_def}} \}, \\
        S_2 &:= \left\{ (e,f) \in E^1(\Ha\Gamma) \times E^1(\Ha\Gamma) \left| \, \begin{aligned}
            &e \neq f \text{ and the pair } (e,f) \text{ violates} \\ 
            &\text{condition (3) from Definition \ref{mainres_nuc_normal_def}}
        \end{aligned} \right. \right\},
    \end{align*}
    and set $n_i := n_i(\Ha\Gamma) := | S_i |$ for $i=1, 2$. 

    \underline{Step 1} First, let us assume $n_2(\Ha\Gamma) = 0$. By applying range decomposition on all edges we may further assume that every edge $e \in E^1(\Ha\Gamma)$ satisfies $|r_{\Ha\Gamma}(e)| \leq 1$. We prove the claim by induction over the number of vertices. 
    
    If $\Ha\Gamma$ has no vertices or $n_1(\Ha\Gamma) = 0$, then $\Ha\Gamma$ is normal and there is nothing to do. Otherwise, choose some edge $e \in S_1(\Ha\Gamma)$. 
    If $r_{\Ha\Gamma}(e) = \emptyset$ then we can delete the edge $e$ without changing the associated $C^\ast$-algebra, see Lemma \ref{min_edge_deletion_trivial}.
    Hence, without loss of generality we may assume $r_{\Ha\Gamma}(e) \neq \emptyset$.
    Then by assumption we have $r_{\Ha\Gamma}(e) \not \subset s_{\Ha\Gamma}(e)$, $|r_{\Ha\Gamma}(e)| = 1$ and $e$ is the only edge starting from $s_{\Ha\Gamma}(e)$, i.e. it is $s_{\Ha\Gamma}(e) \cap s_{\Ha\Gamma}(f) = \emptyset$ for all $f \in E^1(\Ha\Gamma) \setminus \{e\}$.

    Therefore, we may construct a hypergraph $\Ha\Gamma^\prime$ by applying backward contraction on $e$. Further, obtain another hypergraph $\Ha\Gamma^{\prime\prime}$ from $\Ha\Gamma^\prime$  by applying range decomposition on all edges. Evidently, $\Ha\Gamma^{\prime\prime} \leq \Ha\Gamma$ and $\Ha\Gamma^{\prime\prime}$ has fewer vertices than $\Ha\Gamma$. Moreover, we have $C^\ast(\Ha\Gamma^{\prime\prime}) =_M C^\ast(\Ha\Gamma)$ by Theorem \ref{min_big_thm}.

    In order to apply the induction hypothesis we need to check $n_2(\Ha\Gamma^{\prime\prime}) = 0$. Assume that there is a pair $(e^{\prime\prime}, f^{\prime\prime}) \in S_2(\Ha\Gamma^{\prime\prime})$, i.e. $e^{\prime\prime}$ and $f^{\prime\prime}$ are distinct edges in $\Ha\Gamma^{\prime\prime}$ with
    \begin{itemize}
        \item $|s_{\Ha\Gamma^{\prime\prime}}(e^{\prime\prime}) \cap s_{\Ha\Gamma^{\prime\prime}}(f^{\prime\prime})| = 1$,
        \item $|s_{\Ha\Gamma^{\prime\prime}}(e^{\prime\prime})|>1$ or $|s_{\Ha\Gamma^{\prime\prime}}(f^{\prime\prime})|>1$ and
        \item there is no edge $g^{\prime\prime}$ in $\Ha\Gamma^{\prime\prime}$ with $s_{\Ha\Gamma^{\prime\prime}}(e^{\prime\prime}) \cap s_{\Ha\Gamma^{\prime\prime}}(f^{\prime\prime}) \subsetneq s_{\Ha\Gamma^{\prime\prime}}(e^{\prime\prime}) \cap s_{\Ha\Gamma^{\prime\prime}}(g^{\prime\prime})$.
    \end{itemize}
    Clearly, it is not $s_{\Ha\Gamma^{\prime\prime}}(e^{\prime\prime}) = s_{\Ha\Gamma^{\prime\prime}}(f^{\prime\prime})$. A moment's thought shows that there are edges $e^\prime, f^\prime$ in $\Ha\Gamma^\prime$ with $s_{\Ha\Gamma^\prime}(e^\prime) = s_{\Ha\Gamma^{\prime\prime}}(e^{\prime\prime})$ and $s_{\Ha\Gamma^\prime}(f^\prime) = s_{\Ha\Gamma^{\prime\prime}}(f^{\prime\prime})$ such that $(e^\prime, f^\prime) \in S_2(\Ha\Gamma^\prime)$. Recall that $\Ha\Gamma^\prime$ is obtained from $\Ha\Gamma$ by deleting the edge $e$ together with the vertex in $r_{\Ha\Gamma}(e)$ and by replacing $r_{\Ha\Gamma}(e)$ with $s_{\Ha\Gamma}(e)$ in the range or source of every edge different from $e$.
    With this in mind, it is not difficult to see $|s_{\Ha\Gamma}(e^\prime) \cap s_{\Ha\Gamma}(f^\prime)| = 1$ as well as $|s_{\Ha\Gamma}(e^\prime)|>1$ or $|s_{\Ha\Gamma}(f^\prime)| > 1$. Assume that there is an edge $g$ in $\Ha\Gamma$ with
    \begin{align*}
        s_{\Ha\Gamma}(e^\prime) \cap s_{\Ha\Gamma}(f^\prime) \subsetneq s_{\Ha\Gamma}(e^\prime) \cap s_{\Ha\Gamma}(g).
    \end{align*}
    Let $w$ be the unique vertex in $r_{\Ha\Gamma}(e)$. One checks
    \begin{align*}
        s_{\Ha\Gamma^\prime}&(e^\prime) \cap s_{\Ha\Gamma^\prime}(f^\prime) \\
        &= \begin{cases}
            ((s_{\Ha\Gamma}(e^\prime) \cap s_{\Ha\Gamma}(f^\prime)) \setminus \{w\} ) \cup s_{\Ha\Gamma}(e), & w \in s_{\Ha\Gamma}(e^\prime) \cap s_{\Ha\Gamma}(f^\prime), \\
            s_{\Ha\Gamma}(e^\prime) \cap s_{\Ha\Gamma}(f^\prime),
                & \text{otherwise,}
        \end{cases} \\
        & \subsetneq \begin{cases}
            ((s_{\Ha\Gamma}(e^\prime) \cap s_{\Ha\Gamma}(g)) \setminus \{w\} ) \cup s_{\Ha\Gamma}(e),
                & w \in s_{\Ha\Gamma}(e^\prime) \cap s_{\Ha\Gamma}(f^\prime), \\
            s_{\Ha\Gamma}(e^\prime) \cap s_{\Ha\Gamma}(g),
                & \text{otherwise}, \\
        \end{cases} \\
        &= \begin{cases}
            s_{\Ha\Gamma^\prime}(e^\prime) \cap s_{\Ha\Gamma^\prime}(g), & w \in s_{\Ha\Gamma}(e^\prime) \cap s_{\Ha\Gamma}(f^\prime), \\
            ((s_{\Ha\Gamma^\prime}(e^\prime) \cap s_{\Ha\Gamma^\prime}(g)) \setminus s_{\Ha\Gamma}(e)) \cup \{w\},
                & w \in (s_{\Ha\Gamma}(e^\prime) \cap s_{\Ha\Gamma}(g)) \setminus s_{\Ha\Gamma}(f^\prime), \\
            s_{\Ha\Gamma^\prime}(e^\prime) \cap s_{\Ha\Gamma^\prime}(g),
                &\text{otherwise}.
        \end{cases}
    \end{align*}
    Using that the intersection on the left-hand side does not contain $w$, it follows 
    $$ s_{\Ha\Gamma^\prime}(e^\prime) \cap s_{\Ha\Gamma^\prime}(f^\prime) \subsetneq s_{\Ha\Gamma^\prime}(e^\prime) \cap s_{\Ha\Gamma^\prime}(g). $$
    This contradicts the assumption $(e^\prime, f^\prime) \in S_2(\Ha\Gamma^\prime)$. Thus, there is no such edge $g$, and we have $(e^\prime, f^\prime) \in S_2(\Ha\Gamma)$. However, we assumed $S_2(\Ha\Gamma) = \emptyset$. By contradiction, $S_2(\Ha\Gamma^{\prime\prime}) = \emptyset$ and $n_2(\Ha\Gamma^{\prime\prime}) = 0$. 
    
    Since in $\Ha\Gamma^{\prime\prime}$ every edge has at most one vertex in its range we may apply the induction hypothesis on $\Ha\Gamma^{\prime\prime}$ and conclude.

    \underline{Step 2} In the general case, let us use induction over $n_2(\Ha\Gamma)$. If $n_2(\Ha\Gamma) = 0$, then the previous step applies. Otherwise, choose a pair $(e,f) \in S_2(\Ha\Gamma)$ and let $\{w\} = s_{\Ha\Gamma}(e) \cap s_{\Ha\Gamma}(f)$. The negation of condition (3) from Definition \ref{mainres_nuc_normal_def} entails
    \begin{align*}
        w \in s_{\Ha\Gamma}(g) \implies \{w\} = s_{\Ha\Gamma}(e) \cap s_{\Ha\Gamma}(g) \quad \text{for all } g \in E^1(\Ha\Gamma) \setminus \{e\}. \tag{$\ast$}
    \end{align*}
    Construct a hypergraph $\Ha\Gamma^\prime$ by applying source separation on $\{e\}$ at $w$, i.e. $\Ha\Gamma^\prime$ is given by 
    \begin{itemize}
        \item $E^0(\Ha\Gamma^\prime) = E^0(\Ha\Gamma) \cup \{w^\prime\}$
        \item $E^1(\Ha\Gamma^\prime) = E^1(\Ha\Gamma)$,
        \item $r_{\Ha\Gamma^\prime}(g) = \begin{cases}
            r_{\Ha\Gamma}(g),                   &w \not \in r_{\Ha\Gamma}(g), \\
            r_{\Ha\Gamma}(g) \cup \{w^\prime\}, &\text{otherwise},
        \end{cases} \quad$ for all $g \in E^1(\Ha\Gamma^\prime)$,
        \item $s_{\Ha\Gamma^\prime}(g) = \begin{cases}
            s_{\Ha\Gamma}(g),                                         &g \neq e, \\
            (s_{\Ha\Gamma}(g) \setminus \{w\}) \cup \{w^\prime\},     &g = e,
        \end{cases} \quad$ for all $g \in E^1(\Ha\Gamma^\prime)$.
    \end{itemize}
    In view of Proposition \ref{min_edge_sep} and $(\ast)$ we have $C^\ast(\Ha\Gamma^\prime) = C^\ast(\Ha\Gamma)$. In order to apply the induction hypothesis we need to show $n_2(\Ha\Gamma^\prime) < n_2(\Ha\Gamma)$. 
    To do this, let us first show $S_2(\Ha\Gamma^\prime) \subset S_2(\Ha\Gamma)$. For that, assume $(e^\prime, f^\prime) \not \in S_2(\Ha\Gamma)$. We show $(e^\prime, f^\prime) \not \in S_2(\Ha\Gamma^\prime)$.
    As $(e^\prime, f^\prime) \not \in S_2(\Ha\Gamma)$ one of the following cases applies.
    \begin{enumerate}[label=\roman*)]
        \item $|s_{\Ha\Gamma}(e^\prime) \cap s_{\Ha\Gamma}(f^\prime)| \neq 1$. We show that then $|s_{\Ha\Gamma^\prime}(e^\prime) \cap s_{\Ha\Gamma^\prime}(f^\prime)| \neq 1$ holds as well. Indeed, when passing from $\Ha\Gamma$ to $\Ha\Gamma^\prime$ the vertex $w$ in the source of the edge $e$ is replaced with $w^\prime$, but the sources of all other edges remain unchanged. Thus, the only possibility that $|s_{\Ha\Gamma^\prime}(e^\prime) \cap s_{\Ha\Gamma^\prime}(f^\prime)| = 1$, is that either $e^\prime = e$ or $f^\prime = e$, and $w \in s_{\Ha\Gamma}(e^\prime) \cap s_{\Ha\Gamma}(f^\prime)$. Without loss of generality, assume $e^\prime = e$. Then we have
        \begin{align*}
            s_{\Ha\Gamma}(e) \cap s_{\Ha\Gamma}(f)
            = s_{\Ha\Gamma}(e^\prime) \cap s_{\Ha\Gamma}(f) = \{w\}
            \subsetneq s_{\Ha\Gamma}(e^\prime) \cap s_{\Ha\Gamma}(f^\prime)
        \end{align*}
        which contradicts the assumption $(e, f) \in S_2(\Ha\Gamma)$. Thus, 
        $$ |s_{\Ha\Gamma^\prime}(e^\prime) \cap s_{\Ha\Gamma^\prime}(f^\prime)| \neq 1. $$
        \item $|s_{\Ha\Gamma}(e^\prime)| = |s_{\Ha\Gamma}(f^\prime)| = 1$. Since passing from $\Ha\Gamma$ to $\Ha\Gamma^\prime$ does not change the cardinalities of the sources of the edges, it follows directly
        \begin{align*}
            |s_{\Ha\Gamma^\prime}(e^\prime)| = |s_{\Ha\Gamma^\prime}(f^\prime)| = 1.
        \end{align*}
        \item There is some $g^\prime \in E^1(\Ha\Gamma) \setminus \{e^\prime\}$ with $s_{\Ha\Gamma}(e^\prime) \cap s_{\Ha\Gamma}(f^\prime) \subsetneq s_{\Ha\Gamma}(e^\prime) \cap s_{\Ha\Gamma}(g^\prime)$. Without loss of generality, the intersection on the left-hand side has cardinality $1$ since otherwise case (i) applies. Then,
        $
            |s_{\Ha\Gamma}(e^\prime) \cap s_{\Ha\Gamma}(g^\prime)| \geq 2.
        $
        We show
        \begin{align*}
            s_{\Ha\Gamma^\prime}(e^\prime) \cap s_{\Ha\Gamma^\prime}(f^\prime) \subsetneq s_{\Ha\Gamma^\prime}(e^\prime) \cap s_{\Ha\Gamma^\prime}(g^\prime). \tag{$+$}
        \end{align*}
        If $e^\prime, f^\prime, g^\prime \neq e$, we have
        \begin{align*}
            s_{\Ha\Gamma^\prime}(e^\prime) \cap s_{\Ha\Gamma^\prime}(f^\prime)
            &= s_{\Ha\Gamma}(e^\prime) \cap s_{\Ha\Gamma}(f^\prime)
            \subsetneq s_{\Ha\Gamma}(e^\prime) \cap s_{\Ha\Gamma}(g^\prime) \\
            &= s_{\Ha\Gamma^\prime}(e^\prime) \cap s_{\Ha\Gamma^\prime}(g^\prime)
        \end{align*}
        since the involved source sets remain unchanged when passing from $\Ha\Gamma$ to $\Ha\Gamma^\prime$. However, if one of the edges $e^\prime, f^\prime, g^\prime$ equals $e$, then one of the two intersections might lose the vertex $w$ and get smaller. This matters only, if it happens on the right-hand side of ($+$) but not on the left-hand side. In this case, one has $g^\prime = e$ or $e^\prime = e$ and $w \in s_{\Ha\Gamma}(e^\prime) \cap s_{\Ha\Gamma}(g^\prime)$. It follows
        \begin{align*}
            s_{\Ha\Gamma}(e) \cap s_{\Ha\Gamma}(f) 
            = \{w\}
            \subsetneq s_{\Ha\Gamma}(e^\prime) \cap s_{\Ha\Gamma}(g^\prime)
            = \begin{cases}
                s_{\Ha\Gamma}(e) \cap s_{\Ha\Gamma}(g^\prime), &\text{or} \\
                s_{\Ha\Gamma}(e^\prime) \cap s_{\Ha\Gamma}(e),
            \end{cases}
        \end{align*}
        contradicting the assumption $(e,f) \in S_2(\Ha\Gamma)$. Altogether, we get
        \begin{align*}
            s_{\Ha\Gamma^\prime}(e^\prime) \cap s_{\Ha\Gamma^\prime}(f^\prime) 
            \subsetneq s_{\Ha\Gamma^\prime}(e^\prime) \cap s_{\Ha\Gamma^\prime}(g^\prime).
        \end{align*}
    \end{enumerate}
    In any of the above cases it follows $(e^\prime, f^\prime) \not \in S_2(\Ha\Gamma^\prime)$. Thus, we have the implication
    $$(e^\prime, f^\prime) \not \in S_2(\Ha\Gamma) \implies (e^\prime, f^\prime) \not \in S_2(\Ha\Gamma^\prime)$$
    which is equivalent to $S_2(\Ha\Gamma^\prime) \subset S_2(\Ha\Gamma)$. At the same time, one readily checks that the pair $(e,f)$ is in $S_2(\Ha\Gamma)$ but not in $S_2(\Ha\Gamma^\prime)$. Thus, the subset relation is strict, and we have $n_2(\Ha\Gamma^\prime) < n_2(\Ha\Gamma)$. Now, we may apply the induction hypothesis and conclude.
\end{proof}


The proof of Lemma \ref{mainres_normal_hyp_thm} translates directly into an algorithm that produces for any given hypergraph $\Ha\Gamma$ a normalized version $\Ha\Delta \leq \Ha\Gamma$ with $C^\ast(\Ha\Gamma) =_M C^\ast(\Ha\Delta)$. We present this algorithm in pseudocode below.

\begin{algorithm}[H]
    \caption{Hypergraph Normalization} \label{pro_nuc_alg_normal}
    \begin{algorithmic}
        \Procedure{normalize}{hypergraph $\Ha\Gamma$}
        \State $\Ha\Gamma \gets$ take $\Ha\Gamma$ and apply range decomposition on all edges
        \While{True}
            \If {there is a pair $(e,f)$ violating condition (3) from Definition \ref{mainres_nuc_normal_def}}
                \State $\Ha\Gamma \gets$ take $\Ha\Gamma$ and separate the source of $\{e\}$ at $s(e) \cap s(f)$
            \ElsIf {there is an edge $e$ violating condition (2) from Definition \ref{mainres_nuc_normal_def}}
                \If {$r(e) = \emptyset$}
                    \State $\Ha\Gamma \gets$ take $\Ha\Gamma$ and delete the edge $e$
                \Else
                    \State $\Ha\Gamma \gets$ take $\Ha\Gamma$ and apply backward contraction on $e$
                \EndIf
                \State $\Ha\Gamma \gets$ take $\Ha\Gamma$ and apply range decomposition on all edges
            \ElsIf{none of the previous cases applies}
                \State \textbf{break}
            \EndIf
        \EndWhile
        \State \Return $\Ha\Gamma$
        \EndProcedure
    \end{algorithmic}
\end{algorithm}
\section{Hypergraph Reduction}
\label{red_sec}

In this section, we identify situations where cutting some edges of a hypergraph $\Ha\Gamma$ does not change nuclearity of the associated $C^\ast$-algebra. More precisely, given a hypergraph $\Ha\Gamma$ we find three different kinds of sets $S \subset E^1(\Ha\Gamma)$ with the following property: If $\Ha\Delta$ is the hypergraph obtained from $\Ha\Gamma$ by cutting all edges in $S$, then $C^\ast(\Ha\Gamma)$ is nuclear if, and only if, the same holds for $C^\ast(\Ha\Delta)$.

To prepare these results, Section \ref{red_entryexit_closed_edge_sets} introduces the notion of an entry- or exit-closed edge set. In Section \ref{red_easy_edge_sets_subsec} we introduce the notion of an \emph{easy edge set} and show that it has the desired behavior described above. The same result is obtained in Section \ref{red_easy_cycles_subsec} for \emph{easy cycles} and in Section \ref{pro_nuc_simple_quasisinks_subsec} for edges ending in a \emph{simple quasisink}.

\subsection{Entry-/Exit-Closed Edge Sets} \label{red_entryexit_closed_edge_sets}

Let $\Ha\Gamma$ and $\Ha\Delta$ be hypergraphs and let $p \in C^\ast(\Ha\Gamma), q \in C^\ast(\Ha\Delta)$ be projections. Sometimes one observes that the corners $p C^\ast(\Ha\Gamma) p$ and $q C^\ast(\Ha\Delta) q$ are equal although the hypergraphs themselves are different. In the following, we identify two situations where this is true.

First, we need a lemma that describes a dense subset of a hypergraph $C^\ast$-algebra. Its proof is evident from the definitions.

\begin{lemma} \label{pre_hyp_alg_dense_subset_words_lemma}
    Let $\Ha\Gamma$ be a hypergraph. A dense subset of $C^\ast(\Ha\Gamma)$ is spanned by products of the form 
    \begin{align*}
        x = x_1 \dots x_n \qquad \text{with } n \in \N, x_i \in \{p_v, s_e, s_e^\ast: v \in E^0, e \in E^1\},
    \end{align*}
    where for every $i < n$ 
    neither of the following is true:
    \begin{enumerate}[label=\alph*)]
        \item $x_i x_{i+1} = s_e^\ast s_f$ for some edges $e,f \in E^1$.
        \item $x_i x_{i+1} = s_e p_v$ or $x_i x_{i+1} = p_v s_e^\ast$ for some $e \in E^1, v \in E^0$ with $r(e) \neq \emptyset$, and either $v \not \in r(e)$ or $\{v\} = r(e)$.
        \item $x_i x_{i+1} = s_e p_v$ or $x_i x_{i+1} = p_v s_e^\ast$ for some $e \in E^1, v \in E^0$ with $r(e) = \emptyset$, and either $v \not \in s(e)$ or $\{v\} = s(e)$.
        \item $x_i x_{i+1} = p_v s_e$ or $x_i x_{i+1} = s_e^\ast p_v$ for some $e \in E^1, v \in E^0$ with $v \not \in s(e)$ or $\{v\} = s(e)$.
        \item $x_i x_{i+1} = s_e s_f$ or $x_i x_{i+1} = s_f^\ast s_e^\ast$ for some $e, f \in E^1$ with $r(e) \cap s(f) = \emptyset$.
        \item $x_i x_{i+1} = s_e s_f^\ast$ for some $e, f \in E^1$ with $r(e) \neq \emptyset$ and $r(e) \cap r(f) = \emptyset$.
        \item $x_i x_{i+1} = s_e s_f^\ast$ for some $e, f \in E^1$ with $r(e) = \emptyset$ and $s(e) \cap r(f) = \emptyset$.
    \end{enumerate}
\end{lemma}

\begin{definition} \label{pro_nuc_def_entry_exit_closed}
    Let $\Ha\Gamma$ be a hypergraph and let $F \subset E^1(\Ha\Gamma)$ be a set of edges in $\Ha\Gamma$. Then $F$ is \emph{closed under source entries} if
    \begin{align*}
        \forall f \in F, e \in E^1(\Ha\Gamma): s(f) \cap r(e) \neq \emptyset \implies e \in F.
    \end{align*}
    Similarly, $F$ is \emph{closed under range exits} if 
    \begin{align*}
        \forall f \in F, e \in E^1(\Ha\Gamma): r(f) \cap s(e) \neq \emptyset \implies e \in F.
    \end{align*}
\end{definition}

\begin{lemma} \label{pro_nuc_lemma_closed_under_entries_exits}
    Let $\Ha\Gamma$ be a hypergraph, $F \subset E^1(\Ha\Gamma)$ and let $p \in C^\ast(\Ha\Gamma)$ be a projection. Further, obtain $\Ha\Delta$ from $\Ha\Gamma$ by cutting all edges in $F$, and assume that one of the following holds:
    \begin{enumerate}
        \item $F$ is closed under source entries, and $|s_{\Ha\Gamma}(f)| = 1$, $ps_f= 0$ hold for all $f \in F$.
        \item $F$ is closed under range exits, and $|r_{\Ha\Gamma}(f)| = 1$, $s_fp= 0$ hold for all $f \in F$. Further, we have 
        \begin{align*}
            \forall f \in F, e \in E^1(\Ha\Gamma): r_{\Ha\Gamma}(f) \cap r_{\Ha\Gamma}(e) \neq \emptyset \implies e = f. \tag{$\ast$}
        \end{align*}
    \end{enumerate}
    Then $p C^\ast(\Ha\Gamma) p = p C^\ast(\Ha\Delta) p$.
\end{lemma}

\begin{proof}
    Recall from Theorem \ref{min_big_thm} that $C^\ast(\Ha\Delta)$ is a subalgebra of $C^\ast(\Ha\Gamma)$. 
    Further, let $x = x_1 \dots x_n$ be a word with $n \in \N, x_i \in \{p_v, s_e, s_e^\ast: v \in E^0(\Ha\Gamma), e \in E^1(\Ha\Gamma) \}$ as in Lemma \ref{pre_hyp_alg_dense_subset_words_lemma}. Assume $pxp \neq 0$. It suffices to show $x \in C^\ast(\Ha\Delta)$ in the situations (1) and (2).

    Ad (1): 
    We show $x_i \not \in \{s_f, s_f^\ast: f \in F\}$ for all $i \leq n$. Indeed, for any $f \in F$ and all $e \in E^1(\Ha\Gamma), v \in E^0(\Ha\Gamma)$ one observes the following:
    \begin{itemize}
        \item $r_{\Ha\Gamma}(e) \cap s_{\Ha\Gamma}(f) \neq \emptyset$ implies $e \in F$, since $F$ is closed under source entries.
        \item $v \in s_{\Ha\Gamma}(f)$ implies $\{v\} = s_{\Ha\Gamma}(f)$ since $|s_{\Ha\Gamma}(f)| = 1$.
    \end{itemize}
    Combining these observations with the properties of $x$, one checks that as soon as $x_i = s_f$ for some $f \in F$, $i > 1$, then $x_{i-1} = s_g$ for some $g \in F$. Inductively, it follows that either $x_1 \in \{s_f: f \in F\}$ or $x_i \not \in \{s_f: f \in F\}$ holds for all $i \leq n$. However, in the first case we have $pxp = (px_1) x_2 \dots x_n p = 0$ since $ps_f = 0$ holds for all $f \in F$. Hence, it is $x_i \not \in \{s_f: f \in F\}$ for all $i \leq n$. By symmetry, $x_i \not \in \{s_f^\ast: f \in F\}$ for all $i \leq n$ holds as well. As $\{p_v, s_e, s_e^\ast: v \in E^0(\Ha\Gamma), e \in E^1(\Ha\Gamma) \setminus F\}$ is a subset of $C^\ast(\Ha\Delta)$, we obtain $x \in C^\ast(\Ha\Delta)$ as desired.

    Ad (2): Let $x$ be as above and assume $p x p \neq 0$. This time, we show that every occurrence of some $s_f$ with $f \in F$ in the product $x$ is followed by $s_f^\ast$. Indeed, for all $e \in E^1(\Ha\Gamma)$ and $v \in E^0(\Ha\Gamma)$ one observes the following:
    \begin{itemize}
        \item $r_{\Ha\Gamma}(f) \neq \emptyset$ since $|r_{\Ha\Gamma}(f)| = 1$.
        \item $r_{\Ha\Gamma}(f) \cap s_{\Ha\Gamma}(e) \neq \emptyset$ implies $e \in F$ since $F$ is closed under range exits.
        \item $r_{\Ha\Gamma}(f) \cap r_{\Ha\Gamma}(e) \neq \emptyset$ implies $e = f$ due to ($\ast$).
        \item $v \in r_{\Ha\Gamma}(f)$ implies $\{v\} = r_{\Ha\Gamma}(f)$ since $|r_{\Ha\Gamma}(f)|=1$.
    \end{itemize}
    Combining these observations with the properties of $x$, one checks that for all $i < n$, $x_i = s_f$ with $f \in F$ entails $x_{i+1} = s_f^\ast$ or $x_{i+1} \in \{s_g: g \in F\}$.
    So, assume $x_i = s_f$ as above and $x_{i+1} \neq s_f^\ast$ for some $i < n$. Without loss of generality $i < n$ is maximal with this property. Then $x_{i+1} \in \{s_g: g \in F\}$ and by induction one obtains $x_n \in \{s_g: g \in F\}$ as well using maximality of $i$. Then, however, we have $p xp = p x_1 \dots x_{n-1} (x_n p) = 0$ contradicting the assumption that $pxp \neq 0$. Thus, every occurrence of some $s_f$ with $f \in F$ as a factor in the product $x$ is followed by $s_f^\ast$. By symmetry, every occurrence of $s_f^\ast $ with $f \in F$ in $x$ is preceded by $s_f$ as well. Altogether, we get $x \in C^\ast(\Ha\Delta)$ since $C^\ast(\Ha\Delta)$ contains $s_f s_f^\ast$ for all $f \in F$ as well as $\{p_v, s_e, s_e^\ast: v \in E^0(\Ha\Gamma), e \in E^1(\Ha\Gamma) \setminus F\}$.
\end{proof}

\subsection{Elimination of Easy Edge Sets} \label{red_easy_edge_sets_subsec}

In this section, we introduce the notion of an easy edge and the easy edge set it generates. In a hypergraph, all edges in an easy edge set may be cut without changing nuclearity of the associated hypergraph $C^\ast$-algebra.

\begin{definition}[easy edge] \label{pro_nuc_easy_edge_def}
    Let $\Ha\Gamma$ be a hypergraph and let $f_0 \in E^1(\Ha\Gamma)$. Let $F$ be the set
    \begin{align*}
        \{e_1 \in E^1(\Ha\Gamma)| \exists n \in \N, e_2, \dots, e_n \in E^1(\Ha\Gamma): e_n = f_0 \wedge e_1 \dots e_n \text{ is a path in } \Ha\Gamma\}.
    \end{align*}
    The edge $f_0$ is called \emph{easy}\index{easy!edge} if for all $f \in F$ it is
        $|s_{\Ha\Gamma}(f)|=|r_{\Ha\Gamma}(f)|=1$.
    In this case, we call $F$ the \emph{easy edge set}\index{easy!edge set} generated by $f_0$.
\end{definition}

\begin{example}
    In the hypergraph $\Ha\Gamma$ below the edge $f$ is easy, and the edges colored in red form the easy edge set generated by $f$.
    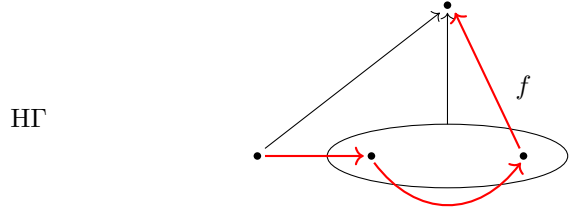
\begin{figure}[H]
        \centering
        \begin{tikzpicture}
            \node[label=0:$\Ha\Gamma$] at (-6, .5) {};
            \node [fill=black, circle, inner sep=1pt] at (1, 0) {};
            \node [fill=black, circle, inner sep=1pt] at (-1, 0) {};
            \node [fill=black, circle, inner sep=1pt] at (0, 2) {};
            \draw[-] (0, 0) ellipse [x radius = 45pt, y radius=12pt];
            \draw[->] (0, .42) -- (0, 1.9);
            \draw[<-, red, thick] (0, 1.75) [partial ellipse=-50:-130:1.5cm and 2.4cm];
            \node [fill=black, circle, inner sep=1pt] at (-2.5, 0) {};
            \draw[->, red, thick] (-2.4, 0) -- (-1.1, 0);
            \draw[->, red, thick] (.95, .1) -- (.1, 1.9);
            \node[label=90:$f$] at (1, .5) {};
            \draw[->] (-2.4, .1) -- (-.1, 1.9);
        \end{tikzpicture}
        \caption{Example of an easy edge}
        \label{pro_nuc_easy_edge_example}
    \end{figure} 
\end{example}

\begin{lemma} \label{pro_nuc_cut_edge_from_source}
    Let $\Ha\Gamma$ be a hypergraph and $f \in E^1(\Ha\Gamma)$ with $|s_{\Ha\Gamma}(f)| = 1$. Further, assume that $s_{\Ha\Gamma}(f)$ is a source, i.e.
    \begin{align*}
        \forall e \in E^1(\Ha\Gamma): \; r_{\Ha\Gamma}(e) \cap s_{\Ha\Gamma}(f) = \emptyset.
    \end{align*}
    Obtain $\Ha\Delta$ from $\Ha\Gamma$ by cutting the edge $f$. Then $C^\ast(\Ha\Gamma)$ is nuclear iff the same holds for $C^\ast(\Ha\Delta)$.
\end{lemma}

\begin{proof}
    Without loss of generality, the edge $f$ has non-empty range. If there is an edge $e \in E^1(\Ha\Gamma) \setminus \{f\}$ with $s_{\Ha\Gamma}(e) \cap s_{\Ha\Gamma}(f) \neq \emptyset$, then obtain $\Ha\Gamma^\prime$ from $\Ha\Gamma$ by separating the source of $f$, i.e.
    \begin{itemize}
        \item $E^0(\Ha\Gamma^\prime) = E^0(\Ha\Gamma) \cup \{w_f\}$,
        \item $E^1(\Ha\Gamma^\prime) = E^1(\Ha\Gamma)$,
        \item $r_{\Ha\Gamma^\prime}(e) = r_{\Ha\Gamma}(e)$ for all $e \in E^1(\Ha\Gamma^\prime)$,
        \item $s_{\Ha\Gamma^\prime}(e) = \begin{cases}
            s_{\Ha\Gamma}(e),       &e \neq f, \\
            \{w_f\},                &e = f,
        \end{cases} \quad$ for all $e \in E^1(\Ha\Gamma^\prime)$.
    \end{itemize}
    Otherwise, set $\Ha\Gamma^\prime := \Ha\Gamma$ and let $w_f \in E^0(\Ha\Gamma^\prime)$ be the vertex with $s_{\Ha\Gamma}(f) = \{w_f\}$.
    By Proposition \ref{min_edge_sep} we have in any case $C^\ast(\Ha\Gamma^\prime) = C^\ast(\Ha\Gamma)$. Further, in $C^\ast(\Ha\Gamma^\prime)$, $f$ is the only edge starting from $w_f$. Hence, applying forward contraction on $f$ does not change the associated $C^\ast$-algebra up to Morita equivalence, see Proposition \ref{min_edge_contraction}. Let $\Ha\Gamma^{\prime\prime}$ be the obtained hypergraph. One readily checks 
    $$ C^\ast(\Ha\Delta) = \C \oplus \C^\ast(\Ha\Gamma^{\prime\prime}) =_M \C \oplus C^\ast(\Ha\Gamma) $$ 
    and this yields the claim.
\end{proof}

\begin{lemma} \label{pro_nuc_lemma_easy_edge_cutting}
    Let $\Ha\Gamma$ be  a hypergraph that contains an easy edge $f_0$. Further, let $F$ be the easy edge set generated by $f_0$ and 
    obtain $\Ha\Delta$ from $\Ha\Gamma$ by cutting all edges in $F$. Then $C^\ast(\Ha\Gamma)$ is nuclear if, and only if, the same holds for $C^\ast(\Ha\Delta)$.
\end{lemma}

\begin{proof}
    Recall from Theorem \ref{min_big_thm} that $C^\ast(\Ha\Delta)$ is a subalgebra of $C^\ast(\Ha\Gamma)$.

    \underline{Step 1} 
    First, assume 
    \begin{align*}
        r_{\Ha\Gamma}(f_0) = s_{\Ha\Gamma}(f) \qquad \qquad \text{ for some } f \in F. \tag{$\ast$}
    \end{align*}
    Let $S \subset E^0(\Ha\Gamma) \cup E^1(\Ha\Gamma)$ be given by
    \begin{align*}
        S := (E^1(\Ha\Gamma) \setminus F) \cup \left(E^0(\Ha\Gamma) \setminus \bigcup_{f \in F} s_{\Ha\Gamma}(f)\right).
    \end{align*}
    We show that $S$ is ideally closed in the sense of Definition \ref{min_ideally_closed_def}. Indeed, we have the following:
    \begin{itemize}
        \item Whenever an edge $e$ is in $S$, then  $r_{\Ha\Gamma}(e)$ is a subset of $S$. Otherwise, there would be an edge $f \in F$ with $r_{\Ha\Gamma}(e) \cap s_{\Ha\Gamma}(f) \neq \emptyset$. By definition of $F$ this implies $e \in F$, i.e. $e \not \in S$.
        \item Whenever an edge $e \in E^1(\Ha\Gamma)$ satisfies $s_{\Ha\Gamma}(e) \subset S$ or $\emptyset \neq r_{\Ha\Gamma}(e) \subset S$, then $e \in S$. Indeed, if $e$ is not in $S$, then we have $e \in F$, and this implies $s_{\Ha\Gamma}(e) \cap S = \emptyset = r_{\Ha\Gamma}(e) \cap S$. For the latter equality, use ($\ast$) to obtain that for every $f \in F$ there is an $f^\prime \in F$ with $r_{\Ha\Gamma}(f) = s_{\Ha\Gamma}(f^\prime)$.
        \item Whenever a vertex $v \in E^0(\Ha\Gamma)$ is not a sink and satisfies 
        $$ v \in s_{\Ha\Gamma}(e) \implies e \in S \text{ for all edges } e \in E^1(\Ha\Gamma), $$ 
        then $v \in S$. Indeed, if $v \not \in S$, then there is an edge $f \in F$ with $v \in s_{\Ha\Gamma}(f)$.
    \end{itemize}

    \underline{Step 2} 
    As $S$ is ideally closed, Lemma \ref{min_lemma_delete_ideally_closed_set} yields the short exact sequence
    \begin{align*}
        0 \to (S) \to C^\ast(\Ha\Gamma) \to C^\ast(\Phi) \to 0,
    \end{align*}
    where $\Phi$ is obtained from $\Ha\Gamma$ by deleting all edges and vertices in $S$. Since all edges that are not in $S$ have exactly one vertex in their range and source, respectively, one verifies that $\Phi$ is an ordinary graph. Thus, $C^\ast(\Phi)$ is nuclear. Since the class of nuclear $C^\ast$-algebras is closed under extensions and ideals, it follows that $C^\ast(\Ha\Gamma)$ is nuclear iff the same holds for $(S)$.

    \underline{Step 3} 
    Set,
    \begin{align*}
        p := 1 - \sum_{f \in F} s_f s_f^\ast 
        = \sum_{e \in S \cap E^1(\Ha\Gamma)} s_e s_e^\ast + \sum_{v \in S \cap E^0(\Ha\Gamma): v \text{ is a sink}} p_v .
    \end{align*}
    Let us show $(S) = (p C^\ast(\Ha\Gamma) p)$. Evidently, $p \in (S)$ and therefore $p C^\ast(\Ha\Gamma) p \subset (S)$. Further, for vertices $v$ in $S$ one has that $v$ is a sink or every edge $e \in E^1(\Ha\Gamma)$ with $v \in s_{\Ha\Gamma}(e)$ is in $S$. In any case, $p_v = p p_v = p_v p = p p_v p$. On the other hand, for every edge $e$ in $S$ with non-empty range, the range $r_{\Ha\Gamma}(e)$ contains only vertices from $S$. Thus, one checks
    \begin{align*}
        p s_e p = s_e s_e^\ast s_e \left( \sum_{w \in r_{\Ha\Gamma}(e)} p_w \right) p = s_e \left( \sum_{w \in r_{\Ha\Gamma}(e)} p_w \right) = s_e \qquad \text{ for all } e \in S.
    \end{align*}
    If $e \in S \cap E^1(\Ha\Gamma)$ has empty range, then
    \begin{align*}
        p s_e p = (s_e s_e^\ast) s_e (s_e s_e^\ast) = s_e.
    \end{align*}
    Altogether, we have $\{p_v, s_e: v \in S \cap E^0(\Ha\Gamma), e \in S \cap E^1(\Ha\Gamma)\} \subset p C^\ast(\Ha\Gamma) p \subset (S)$, and this entails $(S) = (p C^\ast(\Ha\Gamma) p )$.

    \underline{Step 4} 
    Next, use Lemma \ref{pro_nuc_lemma_closed_under_entries_exits} to obtain  $p C^\ast(\Ha\Gamma) p = p C^\ast(\Ha\Delta) p$. Indeed, by definition the set $F$ is closed under source entries and for every edge $f \in F$ it is $|s_{\Ha\Gamma}(f)| = |r_{\Ha\Gamma}(f)| = 1$. Further, for every $f \in F$ one has
    \begin{align*}
        p s_f = \left( 1 - \sum_{f \in F} s_f s_f^\ast \right) s_f s_f^\ast s_f = 0.
    \end{align*}
    Hence, the conditions for Lemma \ref{pro_nuc_lemma_closed_under_entries_exits}(1) are satisfied, and therefore we have the identity $p C^\ast(\Ha\Gamma) p = p C^\ast(\Ha\Delta) p$.
    
    \underline{Step 5}
    Let us show $C^\ast(\Ha\Delta) = \C^{|F|} \oplus p C^\ast(\Ha\Delta)p$. Indeed, from the proof of Proposition \ref{min_edge_cut} it is clear that $C^\ast(\Ha\Delta)$ is the subalgebra of $C^\ast(\Ha\Gamma)$ generated by the elements in 
    \begin{align*}
        \left\{p_v, s_e: v \in E^0(\Ha\Gamma) \setminus \bigcup_{f \in F} s_{\Ha\Gamma}(f), e \in E^1(\Ha\Gamma) \setminus F\right\}
    \end{align*}
    and
    \begin{align*}
        \{s_f s_f^\ast: f \in F\}.
    \end{align*}
    One checks that the elements in the first set are contained in $p C^\ast(\Ha\Delta) p$ while the elements in the latter set are pairwise orthogonal projections. Further, for every $f \in F$ we have $ s_f s_f^\ast \perp p$. From that one gets immediately 
    \begin{align*}
        C^\ast&(\Ha\Delta) \\
        &= C^\ast(s_f s_f^\ast: f \in F) \oplus C^\ast\left(p_v, s_e: v \in E^0(\Ha\Gamma) \setminus \bigcup_{f \in F} s_{\Ha\Gamma}(f), e \in E^1(\Ha\Gamma) \setminus F\right) \\
        &= \C^{|F|} \oplus p C^\ast(\Ha\Delta) p.
    \end{align*}
    In particular, $p C^\ast(\Ha\Delta) p$ is nuclear iff the same holds for $C^\ast(\Ha\Delta)$.
    Putting everything together and using $pC^\ast(\Ha\Delta)p =_M (p C^\ast(\Ha\Delta) p)$, we obtain the following equivalences:
    \begin{align*}
            &C^\ast(\Ha\Gamma) \text{ is nuclear} \\
            \Leftrightarrow \qquad &(S) = (p C^\ast(\Ha\Gamma) p) = (p C^\ast(\Ha\Delta) p) \text{ is nuclear}  \\
            \Leftrightarrow \qquad &p C^\ast(\Ha\Delta) p  \text{ is nuclear} \\
            \Leftrightarrow \qquad &C^\ast(\Ha\Delta)  \text{ is nuclear}.
    \end{align*}

    \underline{Step 6} 
    Finally, let us remove the assumption ($\ast$) that there is an edge $f \in F$ with $r_{\Ha\Gamma}(f_0) = s_{\Ha\Gamma}(f)$. We show the general claim by induction over $|F|$. If $|F|=1$ then the claim follows from Lemma \ref{pro_nuc_cut_edge_from_source}. For $|F|>1$, there are two possibilities: If ($\ast$) is true, then the claim follows from the previous steps. Otherwise, let $f_1, \dots, f_k \in F$ be the edges with $r_{\Ha\Gamma}(f_i) = s_{\Ha\Gamma}(f_0)$. Then the $f_i$ for $i \leq k$ are easy edges in $\Ha\Gamma$ and their generated easy edge sets $F_i$ do not contain $f_0$ since otherwise ($\ast$) would be true. Thus, $|F_i| < |F|$ holds for all $i \leq k$. By induction, we may cut all edges in the sets $F_i$. Afterwards, $s_{\Ha\Gamma}(f_0)$ is a source. By Lemma \ref{pro_nuc_cut_edge_from_source} we may cut the edge $f_0$ without changing nuclearity of the associated $C^\ast$-algebra, and this yields the claim.
\end{proof}

The following corollary is immediate.

\begin{corollary} \label{pro_nuc_easy_paths_corollary}
    Let $\Ha\Gamma$ be a normal hypergraph that contains no easy edge. Then, for all edges $e$ in $\Ha\Gamma$ with non-empty range one of the following holds:
    \begin{enumerate}
        \item $|s_{\Ha\Gamma}(e)| > 1$.
        \item There are edges $e_1, \dots, e_n \in E^1(\Ha\Gamma)$ with $n \geq 2$ and $e_n = e$ such that $e_1 \dots e_n$ is a path in $\Ha\Gamma$ and 
        $$ |s_{\Ha\Gamma}(e_1)| > 1 = |s_{\Ha\Gamma}(e_2)| = |s_{\Ha\Gamma}(e_3)| = \dots = |s_{\Ha\Gamma}(e_n)|. $$
    \end{enumerate}
\end{corollary}


\subsection{Elimination of Easy Cycles} \label{red_easy_cycles_subsec}

In this section, we introduce the notion of an easy cycle. In a hypergraph, all edges from an easy cycle may be cut without changing nuclearity of the associated $C^\ast$-algebra.

\begin{definition}
    Let $\Ha\Gamma$ be a hypergraph. A cycle $\mu = f_1 \dots f_n$ is called \emph{easy}\index{cycle!easy} if for all $i \leq n$ and all $e \in E^1(\Ha\Gamma)$ we have
    \begin{itemize}
        \item $r_{\Ha\Gamma}(f_i) = \{w_i\}$ for suitable vertices $w_i$,
        \item $\{w_i\} \cap r_{\Ha\Gamma}(e) \neq \emptyset \implies e = f_i$,
        \item $\{w_i\} \cap s_{\Ha\Gamma}(e) \neq \emptyset \implies e = f_{i+1}$ if $i < n$,
        \item $\{w_n\} \cap s_{\Ha\Gamma}(e) \neq \emptyset \implies e = f_1$,
    \end{itemize}
    i.e. the edges $f_i$ have exactly one vertex $w_i$ in their range and every vertex $w_i$ has exactly one incoming and exactly one outgoing edge.
\end{definition}

\begin{example}
    Below we present two hypergraphs $\Ha\Delta_1$ and $\Ha\Delta_2$. While in $\Ha\Delta_1$ the edges $f_1$ and $f_2$ form an easy cycle, this is not true in $\Ha\Delta_2$.
    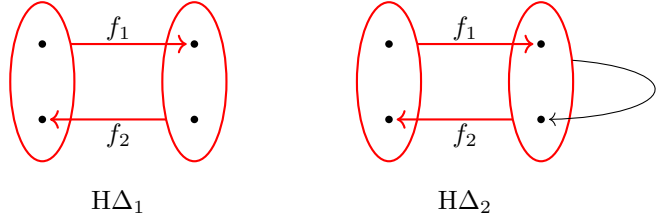
\begin{figure}[H]
        \centering
        \begin{tikzpicture}[baseline=.2]
            \node[label=90:$\Ha\Delta_1$] at (1, -1.5) {};
            \node [fill=black, circle, inner sep=1pt] at (0,1) {};
            \node [fill=black, circle, inner sep=1pt] at (0,0) {};
            \node [fill=black, circle, inner sep=1pt] at (2,1) {};
            \node [fill=black, circle, inner sep=1pt] at (2,0) {};
            \draw[-, red, thick] (0, .5) ellipse [x radius = 12pt, y radius=30pt];
            \draw[->, red, thick] (.37, 1) -- (1.9, 1);
            \node[label=90:$f_1$] at (1, .8) {};
            \draw[-, red, thick] (2, .5) ellipse [x radius = 12pt, y radius=30pt];
            \draw[->, red, thick] (1.63, 0) -- (.1, 0);
            \node[label=-90:$f_2$] at (1, .2) {};
        \end{tikzpicture}
        $\qquad \qquad$
        \begin{tikzpicture}[baseline=.2]
            \node[label=90:$\Ha\Delta_2$] at (1, -1.5) {};
            \node [fill=black, circle, inner sep=1pt] at (0,1) {};
            \node [fill=black, circle, inner sep=1pt] at (0,0) {};
            \node [fill=black, circle, inner sep=1pt] at (2,1) {};
            \node [fill=black, circle, inner sep=1pt] at (2,0) {};
            \draw[-, red, thick] (0, .5) ellipse [x radius = 12pt, y radius=30pt];
            \draw[->, red, thick] (.37, 1) -- (1.9, 1);
            \node[label=90:$f_1$] at (1, .8) {};
            \draw[-, red, thick] (2, .5) ellipse [x radius = 12pt, y radius=30pt];
            \draw[->, red, thick] (1.63, 0) -- (.1, 0);
            \node[label=-90:$f_2$] at (1, .2) {};
            \draw[<-] (2,.4) [partial ellipse=-86:74:1.5cm and .4cm];
        \end{tikzpicture}
        \caption{(Non-)Example of an easy cycle}
        \label{pro_nuc_easy_cycle_example}
    \end{figure} 
\end{example}

\begin{lemma} \label{pro_nuc_easy_cycl}
    Let $\Ha\Gamma$ contain an easy cycle $\mu = f_1 \dots f_n$ and obtain $\Ha\Delta$ from $\Ha\Gamma$ by cutting all edges $f_i$. Then 
    $C^\ast(\Ha\Gamma)$ is nuclear
    if, and only if, the same holds for $C^\ast(\Ha\Delta)$.
\end{lemma}

\begin{proof}
    Let $r_{\Ha\Gamma}(f_i) = \{w_i\}$ for all $i \leq n$ and recall from Theorem \ref{min_big_thm} that $C^\ast(\Ha\Delta)$ is a subalgebra of $C^\ast(\Ha\Gamma)$.

    \underline{Step 1} Define $S := (E^0(\Ha\Gamma) \cup E^1(\Ha\Gamma)) \setminus \bigcup_i \{f_i, w_i\}$. Let us check that $S$ is ideally closed in the sense of Definition \ref{min_ideally_closed_def}.
    \begin{itemize}
        \item If $e \in S$, then $\{w_1, \dots, w_n\} \cap r_{\Ha\Gamma}(e) = \emptyset$ since $\mu$ is an easy cycle. Thus, $r_{\Ha\Gamma}(e) \subset S$.
        \item If $\emptyset \neq r_{\Ha\Gamma}(e) \subset S$ or $s_{\Ha\Gamma}(e) \subset S$ holds for an edge $e$, then $e$ cannot be any of the $f_i$. Thus, $e \in S$.
        \item If $v$ is not a sink and every edge that starts from $v$ is in $S$, then $v$ cannot be any of the $w_i$ since $f_i$ starts from $w_i$ and is not in $S$. Hence, $v \in S$.
    \end{itemize}

    \underline{Step 2} As $S$ is ideally closed, Lemma \ref{min_lemma_delete_ideally_closed_set} yields the short exact sequence 
    \begin{align*}
        0 \to (S) \to C^\ast(\Ha\Gamma) \to C^\ast(\Phi) \to 0,
    \end{align*}
    where $\Phi$ is obtained from $\Ha\Gamma$ by deleting all edges and vertices in $S$.
    It is not hard to verify, that $\Phi$ is an ordinary graph, i.e. $|s_{\Phi}(e)| = |r_{\Phi}(e)| = 1$ for all $e \in E^1(\Phi)$. Indeed, $E^1(\Phi) = \{f_1, \dots, f_n\}$ and $E^0(\Phi) = \{w_1, \dots, w_n\}$.
    Therefore, $C^\ast(\Phi)$ is nuclear as a graph $C^\ast$-algebra. Since the class of nuclear $C^\ast$-algberas is closed under extensions and ideals, it follows that $C^\ast(\Ha\Gamma)$ is nuclear iff the same holds for the ideal $(S)$.

    \underline{Step 3} Set $p := 1 - \sum_{i=1}^n p_{w_i} \in C^\ast(\Ha\Gamma)$. One readily checks $pp_v = p_v$ and $ps_e = s_ep = s_e$ for all vertices $v \in E^0(\Ha\Gamma) \cap S$ and all edges $e \in E^1(\Ha\Gamma) \cap S$. Furthermore, we have
    \begin{align*}
        p = 1 - \sum_{i=1}^n p_{w_i} = \sum_{v \in E^0(\Ha\Gamma)} p_v - \sum_{i=1}^n p_{w_i} = \sum_{v \in E^0(\Ha\Gamma) \cap S} p_v.
    \end{align*}
    Combining both observations, one obtains 
    $$ \{p_v, s_e: v \in S \cap E^0(\Ha\Gamma), e \in S \cap E^1(\Ha\Gamma)\} \subset p C^\ast(\Ha\Gamma)p \subset (S), $$ 
    and therefore $(p C^\ast(\Ha\Gamma)p) = (S)$. In particular, $(S)$ and $p C^\ast(\Ha\Gamma) p$ are Morita-equivalent, so that nuclearity of the former $C^\ast$-algebra is equivalent to nuclearity of the latter $C^\ast$-algebra.

    \underline{Step 4} We show $p C^\ast(\Ha\Gamma) p  = p C^\ast(\Ha\Delta) p$ using Lemma \ref{pro_nuc_lemma_closed_under_entries_exits}(2). Evidently, the set $\{f_i\}$ is closed under range exits. Moreover, $|r_{\Ha\Gamma}(f_i)| = 1$ and $s_{f_i} p = s_{f_i} p_{w_i} p = 0$ hold for all $i \leq n$. Finally, for all $i \leq n$ and $e \in E^1(\Ha\Gamma)$ it is 
    $$ r_{\Ha\Gamma}(f_i) \cap s_{\Ha\Gamma}(e) = \{w_i\} \cap s_{\Ha\Gamma}(e) \neq \emptyset \implies e = f_i.$$ 
    As $\Ha\Delta$ is obtained from $\Ha\Gamma$ by cutting all edges $f_i$, the claim follows from Lemma \ref{pro_nuc_lemma_closed_under_entries_exits}(2). 
    
    Putting the previous steps together, $C^\ast(\Ha\Gamma)$ is nuclear iff the same holds for $p C^\ast(\Ha\Gamma) p = p C^\ast(\Ha\Delta) p$.

    \underline{Step 5} It remains to show that $p C^\ast(\Ha\Delta) p$ is nuclear iff the same holds for $C^\ast(\Ha\Delta)$. However, since in $\Ha\Delta$ the vertices $w_i$ have no incoming edge and the only outgoing edge $f_i$ has empty range, it is not hard to check
    \begin{align*}
        C^\ast(\Ha\Delta) = \C^n \oplus p C^\ast(\Ha\Delta) p.
    \end{align*}
    The claim follows immediately.
\end{proof}

\subsection{Elimination of Simple Quasisinks} \label{pro_nuc_simple_quasisinks_subsec}

In this section, we introduce the notion of a simple quasisink. We will see that in a hypergraph $\Ha\Gamma$ all edges which end in a simple quasisink may be cut without changing nuclearity of the associated $C^\ast$-algebra.

\begin{definition}
    Let $\Ha\Gamma$ be a hypergraph. A vertex $w \in E^0(\Ha\Gamma)$ is called a \emph{simple quasisink}\index{simple quasisink} if
    \begin{itemize}
        \item there is at most one edge $e \in E^1(\Ha\Gamma)$ with $w \in s_{\Ha\Gamma}(e)$ and in this case we have $r_{\Ha\Gamma}(e) = \emptyset$, and
        \item there is at most one edge $e \in E^1(\Ha\Gamma)$ with $w \in r_{\Ha\Gamma}(e)$.
    \end{itemize}
    We say that an edge $f$ ends in a simple quasisink if we have $r_{\Ha\Gamma}(f) = \{w\}$ for a simple quasisink $w$.
\end{definition}

\begin{example}
    The figure below presents three hypergraphs $\Ha\Delta_1, \Ha\Delta_2, \Ha\Delta_3$. While the vertex $w$ is not a simple quasisink in the first two hypergraphs, it is a simple quasisink in $\Ha\Delta_3$.
    \begin{figure}[H]
        \centering
        \begin{tikzpicture}
            \node[label=-90:$\Ha\Delta_1$] at (2, -1) {};
            \node [fill=black, circle, inner sep=1pt] at (.5,1) {};
            \node [fill=black, circle, inner sep=1pt] at (.5,0) {};
            \node [fill=black, circle, inner sep=1pt] at (2,1) {};
            \node [label=90:$w$] at (2, .85) {};
            \node [fill=black, circle, inner sep=1pt] at (2,0) {};
            \draw[->] (.6, 1) -- (1.85, 1);
            \draw[->] (.6, .05) -- (1.85, .85);
            \draw[-] (2, .5) ellipse [x radius = 12pt, y radius=30pt];
            \draw[-] (2.4, .9) -- (3.5, 1.3);
        \end{tikzpicture}
        $\qquad \qquad$
        \begin{tikzpicture}
            \node[label=-90:$\Ha\Delta_2$] at (2, -1) {};
            \node [fill=black, circle, inner sep=1pt] at (.5,1) {};
            \node [fill=black, circle, inner sep=1pt] at (2,1) {};
            \node [label=90:$w$] at (2, .85) {};
            \node [fill=black, circle, inner sep=1pt] at (2,0) {};
            \draw[->] (.6, 1) -- (1.9, 1);
            \draw[-] (2, .5) ellipse [x radius = 12pt, y radius=30pt];
            \draw[-] (2.4, .9) -- (3.5, 1.3);
            \draw[-] (2.4, .1) -- (3.5, -.3);
        \end{tikzpicture}
        $\qquad \qquad$
        \begin{tikzpicture}
            \node[label=-90:$\Ha\Delta_3$] at (2, -1) {};
            \node [fill=black, circle, inner sep=1pt] at (.5,1) {};
            \node [fill=black, circle, inner sep=1pt] at (2,1) {};
            \node [label=90:$w$] at (2, .85) {};
            \node [fill=black, circle, inner sep=1pt] at (2,0) {};
            \draw[->] (.6, 1) -- (1.9, 1);
            \draw[-] (2, .5) ellipse [x radius = 12pt, y radius=30pt];
            \draw[-] (2.4, .9) -- (3.5, 1.3);
        \end{tikzpicture}
        \caption{(Non-)Examples of a simple quasisink}
    \end{figure}
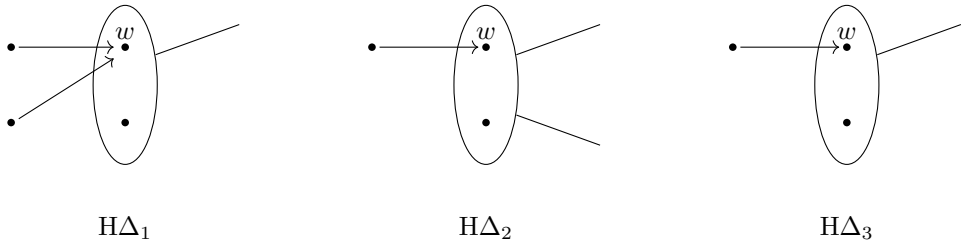 
\end{example}

\begin{lemma} \label{pro_nuc_corner_cut}
    Let $\Ha\Gamma$ be a hypergraph. Assume that $w \in C^\ast(\Ha\Gamma)$ is a simple quasisink in $\Ha\Gamma$ with $\{w\} = r_{\Ha\Gamma}(f)$ for some $f \in E^1(\Ha\Gamma)$ and let $\Ha\Delta$ be obtained from $\Ha\Gamma$ by cutting the edge $f$. Then 
    $C^\ast(\Ha\Gamma)$ is nuclear 
    iff the same holds for $C^\ast(\Ha\Delta)$.
\end{lemma}

\begin{proof}
    First, assume that there is an edge $e \in E^1(\Ha\Gamma)$ with both $r_{\Ha\Gamma}(e) = \emptyset$ and $w \in s_{\Ha\Gamma}(e)$. Then $e$ is the only edge which has $w$ in its source. It follows from Lemma \ref{min_edge_deletion_trivial} that we can remove $w$ from the source of $e$ without changing the associated $C^\ast$-algebra. If $\{w\} = s_{\Ha\Gamma}(e)$, then this means to delete the edge $e$. 

    Thus, without loss of generality $w$ is a sink in $\Ha\Gamma$. Evidently, the set $\{f\}$ is closed under range exits in the sense of Definition \ref{pro_nuc_def_entry_exit_closed}. Now, let
    $$
    p := 1 - p_w \in C^\ast(\Ha\Gamma).
    $$ 
    One readily checks $s_fp = s_f p_w p = 0$. Moreover, it is $|r_{\Ha\Gamma}(f)| = |\{w\}| = 1$ and $f$ is the only edge in $\Ha\Gamma$ which has $w$ in its range. Therefore, the conditions for Lemma \ref{pro_nuc_lemma_closed_under_entries_exits}(2) are satisfied, and we obtain
    \begin{align*}
        p C^\ast(\Ha\Gamma) p = p C^\ast(\Ha\Delta) p.
    \end{align*}
    The corner on the left-hand side is a full corner in $C^\ast(\Ha\Gamma)$ since
    \begin{align*}
        w \not \in s_{\Ha\Gamma}(f)
        &\implies \sum_{v \in s_{\Ha\Gamma}(f)} p_v \in pC^\ast(\Ha\Gamma)p \\
        &\implies s_f = \left(\sum_{v \in s_{\Ha\Gamma}(f)} p_v\right) s_f \in (p C^\ast(\Ha\Gamma) p) \\
        &\implies p_w = s_f^\ast s_f \in (p C^\ast(\Ha\Gamma) p) \\
        &\implies 1 = p + p_w \in (p C^\ast(\Ha\Gamma) p).
    \end{align*}
    Hence, $C^\ast(\Ha\Gamma)$ and $p C^\ast(\Ha\Gamma) p$ are Morita-equivalent. In particular, $C^\ast(\Ha\Gamma)$ is nuclear iff the same holds for $pC^\ast(\Ha\Gamma)p = p C^\ast(\Ha\Delta) p$. 
    
    Finally, in $\Ha\Delta$ the vertex $w$ has neither an incoming nor an outgoing edge. Therefore, it is not hard to check that 
    \begin{align*}
        C^\ast(\Ha\Delta) = \C \oplus pC^\ast(\Ha\Delta)p.
    \end{align*}
    In particular, $pC^\ast(\Ha\Delta)p$ is nuclear iff the same holds for $C^\ast(\Ha\Delta)$. This concludes the proof.
\end{proof}

\subsection{Reduction Algorithm} \label{red_reduction_algo_subsec}

We combine the results of the previous sections with the normalization procedure from Section \ref{norm_sec} to obtain a reduction procedure which transforms any hypergraph $\Ha\Gamma$ into another hypergraph $\Ha\Delta$ such that $C^\ast(\Ha\Gamma)$ is nuclear if, and only if, the same holds for $C^\ast(\Ha\Delta)$.

\begin{algorithm}[H]
    \caption{Hypergraph Reduction} \label{pro_nuc_reduction_algo}
    \begin{algorithmic}[1]
        \Procedure{reduce}{hypergraph $\Ha\Gamma = (E^0, E^1, r, s)$}
        \State $\Ha\Gamma \gets$ normalize $\Ha\Gamma$ 
        \While{True}
            \If{$\exists f \in E^1: f$ is an easy edge}
                \State $\Ha\Gamma \gets$ take $\Ha\Gamma$ and cut all edges in the easy edge set generated by $f$
            \ElsIf{$\exists f_1, \dots, f_n \in E^1: f_1 \dots f_n$ is an easy cycle}
                \State $\Ha\Gamma \gets$ take $\Ha\Gamma$ and cut the edges $f_1, \dots, f_n$
            \ElsIf{$\exists f \in E^1, w \in E^0: r_{\Ha\Gamma}(f) = \{w\} \wedge w$ is a simple quasisink}
                \State $\Ha\Gamma \gets$ take $\Ha\Gamma$ and cut the edge $f$
            \ElsIf{$\exists f \in E^1 \forall e \in E^1 \setminus \{f\}: r_{\Ha\Gamma}(f) = \emptyset = s_{\Ha\Gamma}(f) \cap s_{\Ha\Gamma}(e)$}
                \State $\Ha\Gamma \gets$ take $\Ha\Gamma$ and delete the edge $f$
            \ElsIf{none of the previous cases applies}
                \State \textbf{break}
            \EndIf
        \EndWhile
        \State \Return $\Ha\Gamma$
        \EndProcedure
    \end{algorithmic}
\end{algorithm}

\begin{thm} \label{pro_nuc_reduction_algo_proposition}
    Algorithm \ref{pro_nuc_reduction_algo} terminates for every hypergraph $\Ha\Gamma$. The obtained hypergraph $\Ha\Delta := \mathrm{reduce}(\Ha\Gamma)$ is a normal hypergraph minor of $\Ha\Gamma$ which contains no easy edge, no easy cycle and no edge that ends in a simple quasisink. We call a hypergraph with these properties \emph{reduced}. Further, $C^\ast(\Ha\Gamma)$ is nuclear if, and only if, the same holds for $C^\ast(\Ha\Delta)$.
\end{thm}

\begin{proof}
    Evidently, $\Ha\Delta$ is a hypergraph minor of $\Ha\Gamma$. Moreover, the algorithm terminates since in each application of lines 4 -- 14 either an edge with nonempty range is cut, or an edge is deleted, or the loop breaks. As there are only finitely many edges and vertices in $\Ha\Gamma$, at some point neither of the first two cases applies. Then the algorithm terminates.

    The hypergraph $\Ha\Delta$ contains no easy edge, no easy cycle and no edge that ends in a simple quasisink since the "while" loop only breaks if in all three of these cases the involved edges had been cut. Let us show that $\Ha\Delta$ is normal. Evidently, $\Ha\Gamma$ is normal after line 2. The operations in lines 4 -- 9 only cut some edges and otherwise leave $\Ha\Gamma$ unchanged. Looking at the conditions for normality from Definition \ref{mainres_nuc_normal_def}, there is only one way how this operation can destroy normality of $\Ha\Gamma$: If one cuts an edge $f$ which satisfies $\emptyset \neq r_{\Ha\Gamma}(f) \subset s_{\Ha\Gamma}(f)$ and where there is no edge $e \in E^1(\Ha\Gamma) \setminus \{f\}$ with $s_{\Ha\Gamma}(e) \cap s_{\Ha\Gamma}(f) \neq \emptyset$. However, in this case the operation in line 11 ensures that the edge $f$ is deleted later on which restores normality of $\Ha\Gamma$. By the conditions in line 10, the edge deletion operation in line 11 never destroys normality of $\Ha\Gamma$. Thus, $\Ha\Delta$ is normal.

    Finally, $C^\ast(\Ha\Gamma)$ is nuclear iff the same holds for $C^\ast(\Ha\Delta)$. Indeed, by Lemmas \ref{pro_nuc_lemma_easy_edge_cutting}, \ref{pro_nuc_easy_cycl}, and \ref{pro_nuc_corner_cut} the operations in lines 5, 7, and 9 change the hypergraph $\Ha\Gamma$ so that nuclearity for the original and the modified hypergraph $C^\ast$-algebra are equivalent. By Lemma \ref{min_edge_deletion_trivial} the operation in line 11 does not change the associated hypergraph $C^\ast$-algebra at all. This concludes the proof.
\end{proof}
\section{Reduced Hypergraphs and the Forbidden Minors}
\label{refor_sec}

In this section, we find that a reduced hypergraph $\Ha\Gamma$ has one of the forbidden minors from Table \ref{main_res_for_min_tbl} as soon as $\Ha\Gamma$ contains an edge with nonempty range. Moreover, in Section \ref{refor_hagamma4_leq_hagamma_subsec} we investigate a special situation where the minor $\Ha\Gamma_4$ can be obtained from $\Ha\Gamma$ using only a restricted set of operations which preserve nuclearity of the associated $C^\ast$-algebra. 
Both results together will be the main combinatorial ingredient for the proof of the main result of this article, Theorem \ref{main_res}, which we present in Section \ref{refor_mainres_proof_subsec}.

\subsection{Forbidden Minors}

The following lemma prepares the proof of the next theorem.

\begin{lemma} \label{pro_nuc_path_contraction}
    Let $\Ha\Gamma$ be a normal hypergraph.
    Further, let $f_1 \dots f_n$ be a path in $\Ha\Gamma$ with $|s_{\Ha\Gamma}(f_i)| =  1$ for all $i \geq 2$. Then the hypergraph $\Ha\Delta$ given by
    \begin{itemize}
        \item $E^0(\Ha\Delta) = E^0(\Ha\Gamma)$, 
        \item $E^1(\Ha\Delta) = E^1(\Ha\Gamma) \setminus \{f_2, \dots, f_n\}$,
        \item $s_{\Ha\Delta}(e) = s_{\Ha\Gamma}(e)$ for all $e \in E^1(\Ha\Delta)$,
        \item $r_{\Ha\Delta}(e) = \begin{cases}
            r_{\Ha\Gamma}(e),   & e \neq f_1, \\
            r_{\Ha\Gamma}(f_n), & e = f_1,
        \end{cases} \quad$ for all $e \in E^1(\Ha\Delta)$,
    \end{itemize}
    is a hypergraph minor of $\Ha\Gamma$.
\end{lemma}

\begin{proof}
    Without loss of generality, we have
    \begin{align*}
        f_i f_{i+1} \dots f_{j} \text{ is not a cycle for all } 2 \leq i \leq j \leq n, \tag{$\ast$}
    \end{align*}
    i.e. the path $f_2 \dots f_n$ does not contain a cycle. 
    Indeed, assume that the statement holds under this additional assumption ($\ast$), and let $\mu = f_2 \dots f_n$ contain a cycle. In this case, obtain a shorter path $f_{i_1} \dots f_{i_m}$ by removing all edges from $\mu$ that are part of a cycle, and construct the hypergraph $\Ha\Delta^\prime$ by applying the statement on the path $f_1 f_{i_1} \dots f_{i_m}$. Then, $\Ha\Delta$ is obtained from $\Ha\Delta^\prime$ by deleting all edges in $\{f_2, \dots, f_n\} \setminus \{f_{i_1} \dots f_{i_m}\}$. Hence, we have $\Ha\Delta \leq \Ha\Gamma$ as desired.

    Let us prove the claim under the assumption ($\ast$) by induction over $n$. If $n = 1$, then we have $\Ha\Delta = \Ha\Gamma$ and there is nothing to show. For the induction step, let $f_1, \dots, f_n$ be a path in $\Ha\Gamma$ with $n \geq 2$ and $|s_{\Ha\Gamma}(f_i)| = 1$ for all $i \geq 2$ such that $f_2 \dots f_n$ does not contain a cycle. The induction hypothesis applied on the path $f_2 \dots f_n$ yields a minor $\Ha\Gamma^{(1)} \leq \Ha\Gamma$ with 
    \begin{itemize}
        \item $E^0(\Ha\Gamma^{(1)}) = E^0(\Ha\Gamma)$, 
        \item $E^1(\Ha\Gamma^{(1)}) = E^1(\Ha\Gamma) \setminus \{f_3, \dots, f_n\}$,
        \item $s_{\Ha\Gamma^{(1)}}(e) = s_{\Ha\Gamma}(e)$ for all $e \in E^1(\Ha\Gamma^{(1)})$,
        \item $r_{\Ha\Gamma^{(1)}}(e) = \begin{cases}
            r_{\Ha\Gamma}(e),   & e \neq f_2, \\
            r_{\Ha\Gamma}(f_n), & e = f_2,
        \end{cases} \quad$ for all $e \in E^1(\Ha\Gamma^{(1)})$.
    \end{itemize}
    Let $w \in E^0(\Ha\Gamma^{(1)})$ be the vertex with $\{w\} = s_{\Ha\Gamma^{(1)}}(f_2)$. Now, consider the following constructions:
    \begin{enumerate}
        \item Obtain $\Ha\Gamma^{(2)}$ from $\Ha\Gamma^{(1)}$ by separating the source of $f_2$ in the sense of Remark \ref{min_rem_separate_source_of_e}. Note that $f_2$ must not be a cycle, and therefore due to condition (2) from the definition of normality (Definition \ref{mainres_nuc_normal_def}) there is an edge $e \neq f_2$ in $\Ha\Gamma$ with $w \in s_{\Ha\Gamma}(e)$. As the path $f_2 \dots f_n$ does not contain a cycle, we have $e \not \in \{f_2, \dots, f_n\}$. Therefore, it is $e \in E^1(\Ha\Gamma^{(1)})$. Since $e$ is an edge different from $f_2$ with $w \in s_{\Ha\Gamma^{(1)}}(e)$, the hypergraph $\Ha\Gamma^{(2)}$ contains a new edge $w^\prime \in E^0(\Ha\Gamma^{(2)}) \setminus E^0(\Ha\Gamma^{(1)})$ such that $f_2$ is the only edge with $w^\prime \in s_{\Ha\Gamma^{(2)}}(f_2)$.
        \item Obtain $\Ha\Gamma^{(3)}$ from $\Ha\Gamma^{(2)}$ by applying range decomposition on all edges in the set
        $$ 
        F := \{e \in E^1(\Ha\Gamma^{(2)})| \; w^\prime \in r_{\Ha\Gamma^{(2)}}(e) \}.
        $$
        Since in $\Ha\Gamma^{(1)}$ every edge has exactly one vertex in its range, we can write 
        $$
        E^1(\Ha\Gamma^{(3)}) = E^1(\Ha\Gamma^{(2)}) \cup \{e^\prime: e \in F\}
        $$ 
        where for every $e \in F$ it is $r_{\Ha\Gamma^{(3)}}(e^\prime) = \{w^\prime\}$ and $r_{\Ha\Gamma^{(3)}}(e) = \{w\}$.
        \item Obtain $\Ha\Gamma^{(4)}$ from $\Ha\Gamma^{(3)}$ by deleting all edges in the set
        $$
            \{ e^\prime: e \in F \setminus \{f_1\}\} \cup \{f_1\}.
        $$
        \item Finally, obtain $\Ha\Delta$ from $\Ha\Gamma^{(4)}$ by applying forward contraction on the edge $f_2$.
    \end{enumerate}
    In this way, we get $\Ha\Delta$ from $\Ha\Gamma^{(1)} \leq \Ha\Gamma$ by applying suitable hypergraph minor operations. It follows $\Ha\Delta \leq \Ha\Gamma$.
\end{proof}

Recall that $\Ha\Gamma$ is a reduced hypergraph if $\Ha\Gamma$ is normal and contains no easy edge, no easy cycle and no edge that ends in a simple quasisink.

\begin{thm} \label{pro_nuc_reduction}
    Let $\Ha\Gamma$ be a reduced hypergraph and assume that $\Ha\Gamma$ contains an edge with non-empty range. Then we have $\Ha\Gamma_i \leq \Ha\Gamma$ for some $i=1,2,3,4$, where the $\Ha\Gamma_i$ are the forbidden minors from Table \ref{main_res_for_min_tbl}.
\end{thm}

\begin{proof}
    Let $e \in E^1(\Ha\Gamma)$ and $w \in E^0(\Ha\Gamma)$ satisfy $r_{\Ha\Gamma}(e) = \{w\}$. 
    As $\Ha\Gamma$ does not contain an easy edge, by Corollary \ref{pro_nuc_easy_paths_corollary} we may assume without loss of generality that $|s_{\Ha\Gamma}(e)| \geq 2$. By assumption, the vertex $w$ is not a simple quasisink. Therefore, one of the following cases applies:
    \begin{enumerate}[label=\Alph*)]
        \item It is $w \in s_{\Ha\Gamma}(e)$.
        \item There is an edge $f \neq e$ with nonempty range and $w \in s_{\Ha\Gamma}(f)$.
        \item There are two distinct edges $f, g \in E^1(\Ha\Gamma) \setminus \{e\}$ with empty range and $w \in s_{\Ha\Gamma}(f) \cap s_{\Ha\Gamma}(g)$.
        \item There is an edge $f \neq e$ with $r_{\Ha\Gamma}(f) = \{w\}$.
    \end{enumerate}
    We discuss each of these cases separately. It will be suitable to consider Case (B) last.

    \underline{Case A} 
    Assume that (A) holds and observe that the edge $e$ must not be an easy cycle in $\Ha\Gamma$. Therefore, one of the following three cases (A1) -- (A3) applies.

    Case A1. There is an edge $f \neq e$ with $w \in s_{\Ha\Gamma}(f)$. By condition (3) from the definition of normality (Definition \ref{mainres_nuc_normal_def}), without loss of generality there is at least one vertex $v$ different from $w$ in the intersection $s_{\Ha\Gamma}(e) \cap s_{\Ha\Gamma}(f)$. Now, cut the edge $f$ and delete all edges and vertices except for $e, f, v$ and $w$. This yields the minor $\Ha\Gamma_3$.

    Case A2. There is an edge $f \neq e$ with $\{w\} = r_{\Ha\Gamma}(f)$ and $|s_{\Ha\Gamma}(f)| \geq 2$. In this case, separate the source of the edge $f$ and then delete all edges and vertices except for $e, f$ as well as two vertices in $s(e)$ and $s(f)$, respectively. Afterwards, apply backward contraction on the edge $f$, followed by range decomposition of $e$. Cutting all resulting edges leaves us with the minor $\Ha\Gamma_1$. 
    Below we sketch the involved operations schematically.

    
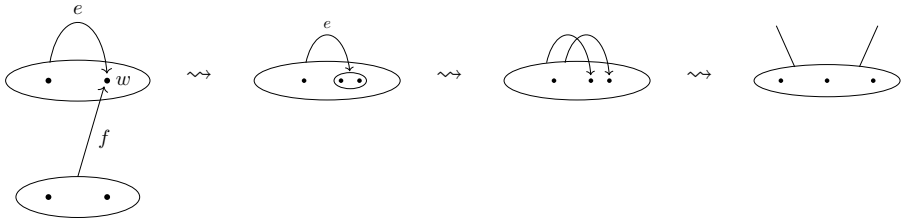
\begin{figure}[H]
        \centering
        \resizebox{2cm}{!}{
        \begin{tikzpicture}[baseline = 2cm]
            \node [fill=black, circle, inner sep=1pt] at (-.5, 2) {};
            \node [fill=black, circle, inner sep=1pt] at (.5, 2) {};
            \node[label=0:$w$] at (.4, 2) {};
            \node [fill=black, circle, inner sep=1pt] at (-.5, 0) {};
            \node [fill=black, circle, inner sep=1pt] at (.5, 0) {};
            \draw[-] (0, 2) ellipse [x radius = 35pt, y radius=10pt];
            \draw[->] (0, 2) [partial ellipse=162:7:.5cm and 1cm];
            \node[label=90:$e$] at (0, 2.9) {};
            \draw[-] (0, 0) ellipse [x radius = 30pt, y radius=10pt];
            \draw[->] (0, .35) -- (.45, 1.9);
            \node[label=0:$f$] at (.1, 1) {};
        \end{tikzpicture}}
        $ \quad \rightsquigarrow \quad $
        \resizebox{2cm}{!}{
        \begin{tikzpicture}[baseline = 2cm]
            \node [fill=black, circle, inner sep=1pt] at (-.5, 2) {};
            \node [fill=black, circle, inner sep=1pt] at (.3, 2) {};
            \node [fill=black, circle, inner sep=1pt] at (.7, 2) {};
            \draw[-] (0, 2) ellipse [x radius = 45pt, y radius=12pt];
            \draw[->] ( 0, 2) [partial ellipse=157:12:.5cm and 1cm];
            \draw[-] (.5, 2) ellipse [x radius = 10pt, y radius=5pt];
            \node[label=90:$e$] at (0, 2.9) {};
        \end{tikzpicture}}
        $ \quad \rightsquigarrow \quad $
        \resizebox{2cm}{!}{
        \begin{tikzpicture}[baseline = 2cm]
            \node [fill=black, circle, inner sep=1pt] at (-.5, 2) {};
            \node [fill=black, circle, inner sep=1pt] at (.3, 2) {};
            \node [fill=black, circle, inner sep=1pt] at (.7, 2) {};
            \draw[-] (0, 2) ellipse [x radius = 45pt, y radius=12pt];
            \draw[->] (-.2, 2) [partial ellipse=157:7:.5cm and 1cm];
            \draw[->] ( .2, 2) [partial ellipse=157:7:.5cm and 1cm];
        \end{tikzpicture}}
        $ \quad \rightsquigarrow \quad $
        \resizebox{2cm}{!}{
        \begin{tikzpicture}[baseline=2cm]
            \node [fill=black, circle, inner sep=1pt] at (-.5, 2) {};
            \node [fill=black, circle, inner sep=1pt] at (.5, 2) {};
            \node [fill=black, circle, inner sep=1pt] at (1.5, 2) {};
            \draw[-] (.5, 2) ellipse [x radius = 45pt, y radius=10pt];
            \draw[-] (-.2, 2.3) -- (-.6, 3.2);
            \draw[-] (1.2, 2.3) -- (1.6, 3.2);
        \end{tikzpicture}}
    \end{figure}

    Case A3. There is an edge $f \neq e$ with $\{w\} = r_{\Ha\Gamma}(f)$ and $|s_{\Ha\Gamma}(f)| = 1$. The edge $f$ must not be easy and therefore by Corollary \ref{pro_nuc_easy_paths_corollary} there is a path $f_1 \dots f_n$ with $f_n = f$ and $|s_{\Ha\Gamma}(f_1)| > 1 = |s_{\Ha\Gamma}(f_i)|$ for all $i \geq 2$. Using that $\Ha\Gamma$ is normal one checks $f_i \neq e$ for all $i \leq n$. Use Lemma \ref{pro_nuc_path_contraction} with the path $f_1 \dots f_n$ to obtain a minor $\Ha\Gamma^\prime$ where $|s(f_1)| \geq 2$ and $r(f_1) = \{w\}$ hold. Then the construction from Case (A2) applied on the hypergraph $\Ha\Gamma^\prime$ yields the minor $\Ha\Gamma_1 \leq \Ha\Gamma^\prime \leq \Ha\Gamma$.

    \underline{Case C} 
    Assume that (C) holds. Without loss of generality case (A) does not apply. Using condition (2) from the definition of normality (Definition \ref{mainres_nuc_normal_def}) one finds an edge $e^\prime \neq e$ with $s_{\Ha\Gamma}(e) \cap s_{\Ha\Gamma}(e^\prime) \neq \emptyset$. Combining with condition (3) from the same definition, without loss of generality there are at least two vertices $v_1$ and $v_2$ in the intersection $s_{\Ha\Gamma}(e) \cap s_{\Ha\Gamma}(e^\prime)$. Now, separate the source of $\{f\}$ at $w$ and afterwards apply range decomposition on the edge $e$. This operation replaces the edge $e$ with two edges $e_1, e_2$ that have the same source as $e$. Finally, cut all edges, and then delete all edges and vertices except for $v_1, v_2, e_1, e_2$ and $e^\prime$. This yields the minor $\Ha\Gamma_2$. Note that $e^\prime = f$ or $e^\prime = g$ is allowed.
    Below we sketch the involved operations schematically.

    
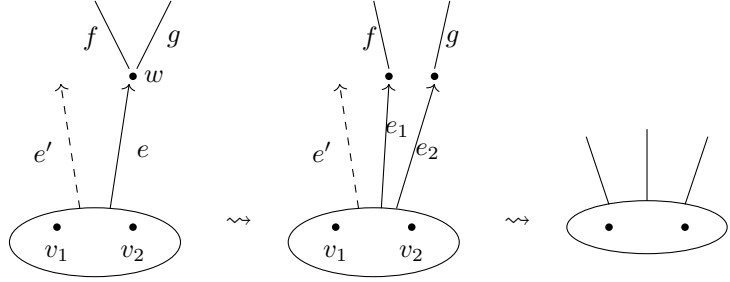
\begin{figure}[H]
        \centering
        \begin{tikzpicture}[baseline = 0cm]
            \node [fill=black, circle, inner sep=1pt] at (.5, 2) {};
            \node[label=0:$w$] at (.4, 2) {};
            \node [fill=black, circle, inner sep=1pt] at (-.5, 0) {};
            \node [fill=black, circle, inner sep=1pt] at (.5, 0) {};
            \node[label=-90:$v_1$] at (-.5, 0) {};
            \node[label=-90:$v_2$] at (.5, 0) {};
            \draw[-] (.45, 2.1) -- (0, 3);
            \draw[-] (.55, 2.1) -- (1, 3);
            \node[label=180:$f$] at (.3, 2.5) {};
            \node[label=0:$g$] at (.7, 2.5) {};
            \draw[-] (0, -.2) ellipse [x radius = 32pt, y radius=13pt];
            \draw[->] (.2, .25) -- (.45, 1.9);
            \draw[->, dashed] (-.2, .25) -- (-.45, 1.9);
            \node[label=0:$e$] at (.3, 1) {};
            \node[label=180:$e^\prime$] at (-.3, 1) {};
        \end{tikzpicture}
        $ \quad \rightsquigarrow \quad $
        \begin{tikzpicture}[baseline = 0cm]
            \node [fill=black, circle, inner sep=1pt] at (.2, 2) {};
            \node [fill=black, circle, inner sep=1pt] at (.8, 2) {};
            \node [fill=black, circle, inner sep=1pt] at (-.5, 0) {};
            \node [fill=black, circle, inner sep=1pt] at (.5, 0) {};
            \node[label=-90:$v_1$] at (-.5, 0) {};
            \node[label=-90:$v_2$] at (.5, 0) {};
            \draw[-] (.2, 2.1) -- (0, 3);
            \draw[-] (.8, 2.1) -- (1, 3);
            \node[label=180:$f$] at (.3, 2.5) {};
            \node[label=0:$g$] at (.7, 2.5) {};
            \draw[-] (0, -.2) ellipse [x radius = 32pt, y radius=13pt];
            \draw[->] (.1, .25) -- (.2, 1.9);
            \draw[->] (.3, .25) -- (.8, 1.9);
            \draw[->, dashed] (-.2, .25) -- (-.45, 1.9);
            \node[label=0:$e_1$] at (-.1, 1.3) {};
            \node[label=0:$e_2$] at (.3, 1) {};
            \node[label=180:$e^\prime$] at (-.3, 1) {};
        \end{tikzpicture}
        $ \quad \rightsquigarrow \quad $
        \begin{tikzpicture}[baseline = .3cm]
            \node [fill=black, circle, inner sep=1pt] at (-.5,.3) {};
            \node [fill=black, circle, inner sep=1pt] at (.5, .3) {};
            \draw[-] (0, .3) ellipse [x radius = 30pt, y radius=10pt];
            \draw[-] (-.5, .6) -- (-.8, 1.5);
            \draw[-] (0, .65) -- (0, 1.6);        
            \draw[-] (.5, .6) -- (.8, 1.5);
        \end{tikzpicture}
    \end{figure}

    \underline{Case D} 
    Assume that (D) holds and observe that $f$ must not be an easy edge. Further, we may assume $w \not \in s_{\Ha\Gamma}(e)$ since otherwise Case (A) applies. This leaves the following three possibilities (D1) -- (D3).

    Case D1. It is $|s_{\Ha\Gamma}(f)| \geq 2$ and $s_{\Ha\Gamma}(e) \cap s_{\Ha\Gamma}(f) = \emptyset$. We may assume $w \not \in s_{\Ha\Gamma}(f)$ since otherwise Case (A) applies for $f$ in the place of $e$. By conditions (2) and (3) from the definition of normality (Definition \ref{mainres_nuc_normal_def}) there is another edge $e^\prime \neq e$ with $|s_{\Ha\Gamma}(e) \cap s_{\Ha\Gamma}(e^\prime)| \geq 2$. Now, transform the hypergraph $\Ha\Gamma$ as follows: First, delete all edges and vertices except for $e, e^\prime, f$ and two vertices in $s_{\Ha\Gamma}(e) \cap s_{\Ha\Gamma}(e^\prime)$ and $s_{\Ha\Gamma}(f)$, respectively. Afterwards, apply backward contraction on the edge $f$, and then decompose the range of $e$. This replaces the edge $e$ with two new edges $e_1, e_2$ that have the same source as $e$. Cutting all edges and deleting all vertices except for those in $s(e) \cap s(e^\prime)$ gives the minor $\Ha\Gamma_2$.
    Below we sketch the involved operations schematically.

    
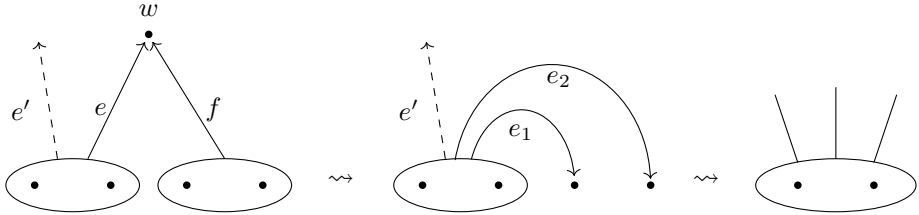
\begin{figure}[H]
        \centering
        \begin{tikzpicture}[baseline = 0cm]
            \node [fill=black, circle, inner sep=1pt] at (1, 2) {};
            \node[label=90:$w$] at (1, 2) {};
            \node [fill=black, circle, inner sep=1pt] at (-.5, 0) {};
            \node [fill=black, circle, inner sep=1pt] at (.5, 0) {};
            \node [fill=black, circle, inner sep=1pt] at (1.5, 0) {};
            \node [fill=black, circle, inner sep=1pt] at (2.5, 0) {};
            \draw[-] (0, 0) ellipse [x radius = 25pt, y radius=10pt];
            \draw[->] (.2, .35) -- (.95, 1.9);
            \draw[->, dashed] (-.2, .35) -- (-.45, 1.9);
            \node[label=180:$e$] at (.7, 1) {};
            \node[label=180:$e^\prime$] at (-.3, 1) {};
            \draw[-] (2, 0) ellipse [x radius = 25pt, y radius=10pt];
            \draw[->] (2, .35) -- (1.05, 1.9);
            \node[label=0:$f$] at (1.5, 1) {};
        \end{tikzpicture}
        $ \quad \rightsquigarrow \quad $
        \begin{tikzpicture}[baseline = 0cm]
            \node [fill=black, circle, inner sep=1pt] at (-.5, 0) {};
            \node [fill=black, circle, inner sep=1pt] at (.5, 0) {};
            \node [fill=black, circle, inner sep=1pt] at (1.5, 0) {};
            \node [fill=black, circle, inner sep=1pt] at (2.5, 0) {};
            \draw[-] (0, 0) ellipse [x radius = 25pt, y radius=10pt];
            \draw[->, dashed] (-.2, .35) -- (-.45, 1.9);
            \node[label=180:$e^\prime$] at (-.3, 1) {};
            \draw[->] (.8, 0) [partial ellipse=160:7:.7cm and 1cm];
            \node[label=180:$e_1$] at (1.2, .7) {};
            \draw[->] (1.2, 0) [partial ellipse=168:3:1.3cm and 1.6cm];
            \node[label=180:$e_2$] at (1.7, 1.4) {};
        \end{tikzpicture}
        $ \quad \rightsquigarrow \quad $
        \begin{tikzpicture}[baseline = .3cm]
            \node [fill=black, circle, inner sep=1pt] at (-.5,.3) {};
            \node [fill=black, circle, inner sep=1pt] at (.5, .3) {};
            \draw[-] (0, .3) ellipse [x radius = 30pt, y radius=10pt];
            \draw[-] (-.5, .6) -- (-.8, 1.5);
            \draw[-] (0, .65) -- (0, 1.6);        
            \draw[-] (.5, .6) -- (.8, 1.5);
        \end{tikzpicture}
    \end{figure}

    Case D2. It is $|s_{\Ha\Gamma}(f)| \geq 2$ and $s_{\Ha\Gamma}(e) \cap s_{\Ha\Gamma}(f) \neq \emptyset$. Again we may assume $w \not \in s_{\Ha\Gamma}(f)$. By condition (3) from the definition of normality (Definition \ref{mainres_nuc_normal_def}) there are two possibilities:
    \begin{itemize}
        \item It is $|s_{\Ha\Gamma}(f) \cap s_{\Ha\Gamma}(e)| \geq 2$.
        \item There is an edge $g \neq e, f$ such that $s_{\Ha\Gamma}(f) \cap s_{\Ha\Gamma}(e) \subsetneq s_{\Ha\Gamma}(g) \cap s_{\Ha\Gamma}(e)$.
    \end{itemize}
    In the latter case, separate the source of $f$ and then use the same construction as in Case (D1) to obtain the minor $\Ha\Gamma_2$. In the first case, there are at least two vertices $v_1, v_2$ in the intersection $s_{\Ha\Gamma}(e) \cap s_{\Ha\Gamma}(f)$. Delete all edges and vertices except for $e, f, v_1, v_2$ and $w$. This yields the minor $\Ha\Gamma_4$.

    Case D3. It is $|s_{\Ha\Gamma}(f)| = 1$. Since $f$ must not be an easy edge, there is a path $f_1 \dots f_n$ in $\Ha\Gamma$ with $n \geq 2, f_n = f$ and $|s_{\Ha\Gamma}(f_1)| > 1 = |s_{\Ha\Gamma}(f_i)|$ for all $i \geq 2$ (see Corollary \ref{pro_nuc_easy_paths_corollary}). Let us distinguish two cases:

    Case D3.1. We have $s_{\Ha\Gamma}(f_1) \cap s_{\Ha\Gamma}(e) = \emptyset$. One easily checks, that then $f_i \neq e$ holds for all $i \leq n$. Use Lemma \ref{pro_nuc_path_contraction} to obtain a minor $\Ha\Gamma^\prime$ where $r(f_1) = \{w\}$. Then Case (D1) applies and yields the minor $\Ha\Gamma_2$.

    Case D3.2. We have $f_1 = e$.
    Then Lemma \ref{pro_nuc_path_contraction} applied on the path $f_2 \dots f_n$ yields a hypergraph minor $\Ha\Gamma^\prime \leq \Ha\Gamma$ where $\{w\} = s(f_2) = r(f_2)$. Delete all edges and vertices except for $e, f_2, w$ and two vertices in $s(e)$. 
    Finally, apply backward contraction on the edge $e$, range decomposition on the edge $f_2$ and afterwards cut one of the obtained edges. This yields the minor $\Ha\Gamma_3$.
    Below we sketch the involved operations schematically.

    
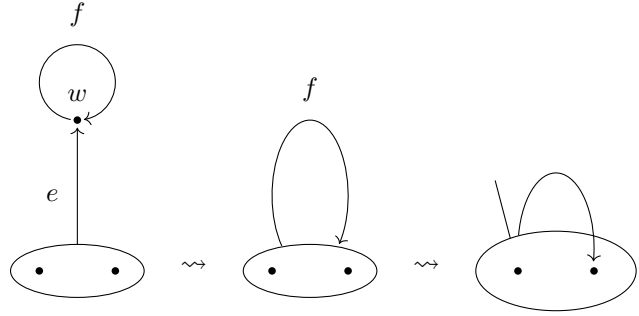
\begin{figure}[H]
        \centering
        \begin{tikzpicture}[baseline = 0cm]
            \node [fill=black, circle, inner sep=1pt] at (0, 2) {};
            \node[label=90:$w$] at (0, 2) {};
            \node [fill=black, circle, inner sep=1pt] at (-.5, 0) {};
            \node [fill=black, circle, inner sep=1pt] at (.5, 0) {};
            \draw[-] (0, 0) ellipse [x radius = 25pt, y radius=10pt];
            \draw[->] (0, .35) -- (0, 1.9);
            \node[label=180:$e$] at (0, 1) {};
            \draw[<-] (0, 2.5) [partial ellipse=-80:260:.5cm and .5cm];
            \node[label=90:$f$] at (0, 3) {};
        \end{tikzpicture}
        $ \quad \rightsquigarrow \quad $
        \begin{tikzpicture}[baseline = 0cm]
            \node [fill=black, circle, inner sep=1pt] at (-.5, 0) {};
            \node [fill=black, circle, inner sep=1pt] at (.5, 0) {};
            \draw[-] (0, 0) ellipse [x radius = 25pt, y radius=10pt];
            \draw[<-] (0, 1) [partial ellipse=-40:222:.5cm and 1cm];
            \node[label=90:$f$] at (0, 2) {};
        \end{tikzpicture}
        $ \quad \rightsquigarrow \quad $
        \begin{tikzpicture}[baseline = .3cm]
            \node [fill=black, circle, inner sep=1pt] at (-.5,.3) {};
            \node [fill=black, circle, inner sep=1pt] at (.5, .3) {};
            \draw[-] (0, .3) ellipse [x radius = 30pt, y radius=15pt];
            \draw[-] (-.6, .73) -- (-.8, 1.5);
            \draw[<-] (0,.6) [partial ellipse=-10:170:.5cm and 1cm];
        \end{tikzpicture}
    \end{figure}

    Case D3.3. It is $f_1 \neq e$ and $s_{\Ha\Gamma}(f_1) \cap s_{\Ha\Gamma}(e) \neq \emptyset$. 
    Due to condition (3)  from the definition of normality (Definition \ref{mainres_nuc_normal_def}) there are two possibilities:
    \begin{itemize}
        \item It is $|s_{\Ha\Gamma}(f_1) \cap s_{\Ha\Gamma}(e)| \geq 2$.
        \item There is an edge $g \neq e, f_1$ such that $s_{\Ha\Gamma}(f_1) \cap s_{\Ha\Gamma}(e) \subsetneq s_{\Ha\Gamma}(g) \cap s_{\Ha\Gamma}(e)$.
    \end{itemize}
    In the latter case, separate the source of $f_1$ and then use the same construction as in Case (D3.1) to obtain the minor $\Ha\Gamma_2$. 
    Otherwise, there are at least two vertices in the intersection $s_{\Ha\Gamma}(f_1) \cap s_{\Ha\Gamma}(e)$. Moreover, without loss of generality $r_{\Ha\Gamma}(f_1) \cap s_{\Ha\Gamma}(f_1) = \emptyset$ since otherwise Case (A) applies for the edge $f_1$. Similarly, we may assume without loss of generality that $r_{\Ha\Gamma}(f_2) \cap s_{\Ha\Gamma}(f_2) = \emptyset$ since otherwise Case (D3.2) applies for the edge $f_1$ in the place of $e$. Now, by conditions (2) and (3) from the definition of normality, there is an edge $f_2^\prime \neq f_2$ in $\Ha\Gamma$ with $s_{\Ha\Gamma}(f_2) = s_{\Ha\Gamma}(f_2^\prime)$. Separate the source of $f_2^\prime$ and afterwards apply range decomposition on $f_1$. This operation replaces $f_1$ with two new edges $f_1^{(1)}$ and $f_1^{(2)}$. Finally, delete all edges and vertices except for $e, f_1^{(1)}, f_1^{(2)}$ and two vertices in $s_{\Ha\Gamma}(e) \cap s_{\Ha\Gamma}(f_1)$. This yields the minor $\Ha\Gamma_2$.
    Below we sketch the involved operations schematically.

    
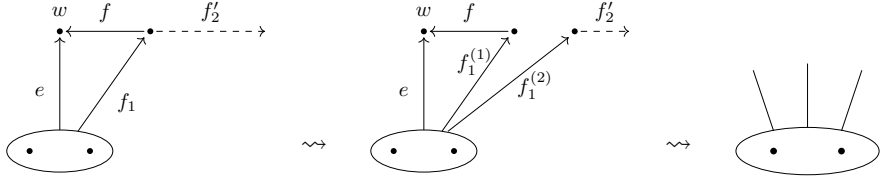
\begin{figure}[H]
        \centering
        \resizebox{3.5cm}{!}{
        \begin{tikzpicture}[baseline = 0cm]
            \node [fill=black, circle, inner sep=1pt] at (0, 2) {};
            \node[label=90:$w$] at (0, 2) {};
            \node [fill=black, circle, inner sep=1pt] at (-.5, 0) {};
            \node [fill=black, circle, inner sep=1pt] at (.5, 0) {};
            \node [fill=black, circle, inner sep=1pt] at (1.5, 2) {};
            \draw[-] (0, 0) ellipse [x radius = 25pt, y radius=10pt];
            \draw[->] (0, .35) -- (0, 1.9);
            \node[label=180:$e$] at (0, 1) {};
            \draw[->] (.3, .33) -- (1.4, 1.9);
            \node[label=0:$f_1$] at (.7, .8) {};
            \draw[->] (1.4, 2) -- (.1, 2);
            \node[label=90:$f$] at (.75, 1.9) {};
            \draw[->, dashed] (1.6, 2) -- (3.4, 2);
            \node[label=90:$f_2^\prime$] at (2.5, 1.9) {};
        \end{tikzpicture}}
        $ \quad \rightsquigarrow \quad $
        \resizebox{3.5cm}{!}{
        \begin{tikzpicture}[baseline = 0cm]
            \node [fill=black, circle, inner sep=1pt] at (0, 2) {};
            \node[label=90:$w$] at (0, 2) {};
            \node [fill=black, circle, inner sep=1pt] at (-.5, 0) {};
            \node [fill=black, circle, inner sep=1pt] at (.5, 0) {};
            \node [fill=black, circle, inner sep=1pt] at (1.5, 2) {};
            \node [fill=black, circle, inner sep=1pt] at (2.5, 2) {};
            \draw[-] (0, 0) ellipse [x radius = 25pt, y radius=10pt];
            \draw[->] (0, .35) -- (0, 1.9);
            \node[label=180:$e$] at (0, 1) {};
            \draw[->] (.3, .33) -- (1.4, 1.9);
            \node[label=0:$f_1^{(1)}$] at (.3, 1.5) {};
            \draw[->] (.4, .33) -- (2.4, 1.9);
            \node[label=0:$f_1^{(2)}$] at (1.3, 1.1) {};
            \draw[->] (1.4, 2) -- (.1, 2);
            \node[label=90:$f$] at (.75, 1.9) {};
            \draw[->, dashed] (2.6, 2) -- (3.4, 2);
            \node[label=90:$f_2^\prime$] at (3, 1.9) {};
        \end{tikzpicture}}
        $ \quad \rightsquigarrow \quad $
        \resizebox{2cm}{!}{
        \begin{tikzpicture}[baseline = .3cm]
            \node [fill=black, circle, inner sep=1pt] at (-.5,.3) {};
            \node [fill=black, circle, inner sep=1pt] at (.5, .3) {};
            \draw[-] (0, .3) ellipse [x radius = 30pt, y radius=10pt];
            \draw[-] (-.5, .6) -- (-.8, 1.5);
            \draw[-] (0, .65) -- (0, 1.6);        
            \draw[-] (.5, .6) -- (.8, 1.5);
        \end{tikzpicture}}
    \end{figure}

    \underline{Case B} 
    Finally, assume that (B) holds and distinguish the following two cases (B1) -- (B2).

    Case B1. It is $|s_{\Ha\Gamma}(f)| = 1$. Then there are two possibilities. If $s_{\Ha\Gamma}(f) = r_{\Ha\Gamma}(f)$, then Case (D) applies. Otherwise, by conditions (2) and (3) from the definition of normality (Definition \ref{mainres_nuc_normal_def}), there is another edge $f^\prime \neq f$ with $s_{\Ha\Gamma}(f^\prime) = \{w\}$. After cutting the edges $f$ and $f^\prime$ one is in the same situation as in Case (C). Similarly as above, one obtains the minor $\Ha\Gamma_2$.

    Case B2. None of the previous cases (A), (C), (D), (B1) applies for any edge with nonempty range. Then there is an edge $e_2$ with $\{w\} \subsetneq s_{\Ha\Gamma}(e_2)$ and $r_{\Ha\Gamma}(e_2) \neq \emptyset$. Let $\{w_2\} := r_{\Ha\Gamma}(e_2)$. Due to the fact that none of the cases (A), (C), (D), (B1) applies for $e_2$, there is an edge $e_3$ and a vertex $w_3$ with $\{w_2\} \subsetneq s_{\Ha\Gamma}(e_3)$ and $\{w_3\} = r_{\Ha\Gamma}(e_2) \neq \emptyset$. Inductively repeating this argument and using that $\Ha\Gamma$ has only finitely many edges, one finds a cycle $f_1 \dots f_n \in \Ha\Gamma$ and vertices $v_1, \dots, v_n$ such that
    \begin{align*}
        r_{\Ha\Gamma}(f_n) = \{v_n\} &\subsetneq s_{\Ha\Gamma}(f_1), \\
        r_{\Ha\Gamma}(f_1) = \{v_1\} &\subsetneq s_{\Ha\Gamma}(f_2), \\
        &\dots, \\
        r_{\Ha\Gamma}(f_{n-1}) = \{v_{n-1}\} &\subsetneq s_{\Ha\Gamma}(f_n).
    \end{align*}
    As $f_1, \dots, f_n$ must not be an easy cycle, there is an $i \leq n$ such that the vertex $v_i$ has two different incoming edges or $v_i$ has two different outgoing edges. Without loss of generality, $v_1$ has this property. However, $v_1$ must not have an incoming edge different from $f_1$ since then Case (D) would apply for the edge $f_1$. Hence, $v_1$ has an outgoing edge different from $f_2$ which we call $f_2^\prime$. After cutting the edges $f_2$ and $f_2^\prime$ the edge $f_1$ has the same property as the edge $e$ in Case (C). Therefore, the argument from the discussion of Case (C) yields the minor $\Ha\Gamma_2$.
\end{proof}

\subsection{Special Situation for \texorpdfstring{$\Ha\Gamma_4 \leq \Ha\Gamma$}{HGamma\_4 leq HGamma}} \label{refor_hagamma4_leq_hagamma_subsec}

From Theorem \ref{pro_nuc_reduction} we know that any reduced hypergraph $\Ha\Gamma$ which contains an edge with non-empty range has one of the forbidden minors $\Ha\Gamma_1, \dots, \Ha\Gamma_4$. In this section, we investigate the case where we have $\Ha\Gamma_4 \leq \Ha\Gamma$ and $\Ha\Gamma_i \not \leq \Ha\Gamma$ for all $i \leq 3$ at the same time.

\begin{proposition} \label{pro_nuc_HGamma4_minor}
    Let $\Ha\Gamma$ be a reduced hypergraph. Assume that $\Ha\Gamma_4$ is a minor of $\Ha\Gamma$ and that $\Ha\Gamma_i \not \leq \Ha\Gamma$ holds for $i \leq 3$. Then $\Ha\Gamma_4$ can be obtained from $\Ha\Gamma$ using only the following operations:
    \begin{itemize}
        \item deletion of an ideally closed set in the sense of Definition \ref{min_ideally_closed_def}
        \item removing a vertex from the source of an edge as in Lemma \ref{min_edge_deletion_trivial}
    \end{itemize}
    Both operations preserve nuclearity of the associated $C^\ast$-algebra.
\end{proposition}

\begin{proof}
    \underline{Step 1} If every edge $e \in E^1(\Ha\Gamma)$ has empty range, then $\Ha\Gamma$ cannot have the minor $\Ha\Gamma_4$. Using Corollary \ref{pro_nuc_easy_paths_corollary}, there is an edge $e$ with $r_{\Ha\Gamma}(e) \neq \emptyset$ and $|s_{\Ha\Gamma}(e)| \geq 2$. A close investigation of the case distinction from the proof of Theorem \ref{pro_nuc_reduction} reveals that Case (D2) must apply since in all other cases $\Ha\Gamma$ has one of the hypergraphs $\Ha\Gamma_1, \Ha\Gamma_2, \Ha\Gamma_3$ as a minor. Hence, for every edge $e$ with $r_{\Ha\Gamma}(e) \neq \emptyset$ and $|s_{\Ha\Gamma}(e)| \geq 2$ there is another edge $e^\prime \neq e$ with $r_{\Ha\Gamma}(e) = r_{\Ha\Gamma}(e^\prime)$ and $|s_{\Ha\Gamma}(e) \cap s_{\Ha\Gamma}(e^\prime)| \geq 2$.

    \underline{Step 2} 
    Let $F := \{f \in E^1(\Ha\Gamma): |s_{\Ha\Gamma}(f)| \geq 2 \text{ and } r_{\Ha\Gamma}(f) \neq \emptyset\}$. By the previous step, the set $F$ is nonempty. We show that there is an edge $f \in F$ such that $r_{\Ha\Gamma}(f) \cap s_{\Ha\Gamma}(e) = \emptyset$ holds for all edges $e \in E^1(\Ha\Gamma)$ with $r_{\Ha\Gamma}(e) \neq \emptyset$. 
    Indeed, assume that this is not true, and let $f_1 \in F$. By assumption there is another edge $f_2 \in E^1(\Ha\Gamma) \setminus \{f_1\}$ such that $r_{\Ha\Gamma}(f_1) \subset s_{\Ha\Gamma}(f_2)$ and $r_{\Ha\Gamma}(f_2) \neq \emptyset$. We prove $f_2 \in F$. First, assume $|s_{\Ha\Gamma}(f_2)| = 1$. There are two possibilities:
    \begin{itemize}
        \item It is $r_{\Ha\Gamma}(f_2) = s_{\Ha\Gamma}(f_2)$. Then it is not difficult to obtain the minor $\Ha\Gamma_3$ similarly as in Case (D3.2) from the proof of the previous theorem.
        \item It is not $r_{\Ha\Gamma}(f_2) = s_{\Ha\Gamma}(f_2)$. Then conditions (2) and (3) from the definition of normality (Definition \ref{mainres_nuc_normal_def}) yield another edge $f_2^\prime$ with $s_{\Ha\Gamma}(f_2) = s_{\Ha\Gamma}(f_2^\prime)$. Using the construction from Case (C) in the proof of Theorem \ref{pro_nuc_reduction} we get the minor $\Ha\Gamma_2$.
    \end{itemize}
    In any event, this contradicts the assumption $\Ha\Gamma_i \not \leq \Ha\Gamma$ for $i \leq 3$.
    Hence, $|s_{\Ha\Gamma}(f_2)| \geq 2$ and $f_2$ is in the set $F$. It follows that there is a path $f_1 \dots f_{|E^1(\Ha\Gamma)|+1}$ in $\Ha\Gamma$ which contains only edges from $F$. Clearly, this path has a closed subpath. By removing superfluous edges one obtains a cycle $g_1 \dots g_n$ with $g_i \in F$ for all $i \leq n$.
    Now, it is not difficult to obtain the hypergraph minor $\Ha\Gamma_3$ from $\Ha\Gamma$. By contradiction this proves that there is an edge $f \in F$ such that $r_{\Ha\Gamma}(f) \cap s_{\Ha\Gamma}(e) = \emptyset$ holds for all edges $e \in E^1(\Ha\Gamma)$ with nonempty range.

    \underline{Step 3} 
    By the previous steps there are $v_1, v_2, w \in E^0(\Ha\Gamma)$ and $f, f^\prime \in E^1(\Ha\Gamma)$ such that 
    \begin{align*}
        \{v_1, v_2\} \subset s_{\Ha\Gamma}(f) \cap s_{\Ha\Gamma}(f^\prime), \quad \{w\} = r_{\Ha\Gamma}(f) = r_{\Ha\Gamma}(f^\prime), \quad r_{\Ha\Gamma}(f) \cap s_{\Ha\Gamma}(e) = \emptyset
    \end{align*}
    hold for all $e \in E^1(\Ha\Gamma)$ with $r_{\Ha\Gamma}(e) \neq \emptyset$. We show that there is no edge $e \in E^1(\Ha\Gamma) \setminus \{f, f^\prime\}$ such that $r_{\Ha\Gamma}(e) = \{w\}$. Assume that this is not true and let $e \in E^1(\Ha\Gamma)$ have range $\{w\}$. There are two possibilities: If $|s_{\Ha\Gamma}(e)| \geq 2$, then the construction from Case (D1) in the proof of Theorem \ref{pro_nuc_reduction} yields the minor $\Ha\Gamma_2$. Otherwise, the argument from Case (D3.1) yields the same minor. However, by assumption $\Ha\Gamma_2 \not \leq \Ha\Gamma$, and therefore we obtain the claim by contradiction.

    \underline{Step 4} 
    Let us show that there is no edge $e \in E^1(\Ha\Gamma)$ with $r_{\Ha\Gamma}(e) \subset \{v_1, v_2\}$. Assume the opposite and, without loss of generality, let $e \neq f, f^\prime$ be an edge with $r_{\Ha\Gamma^\prime}(e) = \{v_1\}$. Distinguish the following cases:

    Case 1. It is $|s_{\Ha\Gamma^\prime}(e)| \geq 2$. In this case, separate the source of $e$, and delete all edges and vertices except for $e, f, f^\prime, v_1, v_2, w$ as well as two vertices in $s(e)$. Afterwards, apply backward contraction on the edge $e$. This yields the minor $\Ha\Gamma_1$. 

    Case 2. It is $|s_{\Ha\Gamma^\prime}(e)| = 1$. Since the edge $e$ must not be easy, by Corollary \ref{pro_nuc_easy_paths_corollary} there is a path $e_1 \dots e_n$ in $\Ha\Gamma^\prime$ with $e_n = e$ and $|s_{\Ha\Gamma^\prime}(e_1)| > 1 = |s_{\Ha\Gamma^\prime}(e_i)|$ for all $i \geq 2$. Apply Lemma \ref{pro_nuc_path_contraction} to obtain a minor where $|s_{\Ha\Gamma^\prime}(e_1)| > 1$ and $s_{\Ha\Gamma^\prime}(e_1) = r_{\Ha\Gamma^\prime}(e_1) = \{v_1\}$. Now, the construction from Case (1) yields the minor $\Ha\Gamma_1$.

    Summarizing, as soon as there is an edge $e \neq f, f^\prime$ with $r_{\Ha\Gamma^\prime}(e) \subset \{v_1, v_2\}$, then $\Ha\Gamma_1$ is a minor of $\Ha\Gamma$. By contradiction, 
    it follows that there are no edges $e$ with $r_{\Ha\Gamma^\prime}(e) \subset \{v_1, v_2\}$.

    \underline{Step 5} Next, let us show that there is at most one edge $e \in E^1(\Ha\Gamma)$ with $w \in s_{\Ha\Gamma}(e)$. Assume the opposite, and let $e, e^\prime \in E^1(\Ha\Gamma)$ be edges with $w \in s_{\Ha\Gamma}(e) \cap s_{\Ha\Gamma}(e^\prime)$. Then a similar construction as in Case (C) of the proof of Theorem \ref{pro_nuc_reduction} yields the minor $\Ha\Gamma_2$. This proves the claim by contradiction.

    Assume that $e \in E^1(\Ha\Gamma)$ is an edge with $w \in s_{\Ha\Gamma}(e)$. By Step (2) we know that $e$ has empty range. Construct a hypergraph $\Ha\Gamma^\prime$ by removing the vertex $w$ from the source of $e$ as in Lemma \ref{min_edge_deletion_trivial}. One easily checks that the assumptions for this lemma are satisfied. Hence, we have $C^\ast(\Ha\Gamma^\prime) = C^\ast(\Ha\Gamma)$. In $\Ha\Gamma^\prime$ the vertex $w$ is a sink.

    \underline{Step 6} 
    Set
    \begin{align*}
        S := (E^0(\Ha\Gamma^\prime) \cup E^1(\Ha\Gamma^\prime)) \setminus \{v_1, v_2, w, f, f^\prime\}.
    \end{align*}
    We show that $S$ is ideally closed.
    We check the three conditions from Definition \ref{min_ideally_closed_def}.
    \begin{itemize}
        \item Assume that $e$ is an edge in $S$. Then it is $e \not \in \{f, f^\prime\}$. Combining Steps (3) and (4) one observes $r_{\Ha\Gamma^\prime}(e) \subset E^0(\Ha\Gamma^\prime) \setminus \{v_1, v_2, w \} \subset S$.
        \item Assume that $e \in E^1(\Ha\Gamma^\prime)$ satisfies $s_{\Ha\Gamma^\prime}(e) \subset S$ or $\emptyset \neq r_{\Ha\Gamma^\prime}(e) \subset S$. Both claims are not true for $f, f^\prime$. Therefore, $e \in E^1(\Ha\Gamma) \setminus \{f, f^\prime\} \subset S$.
        \item Finally, assume that $v \in E^0(\Ha\Gamma^\prime)$ is not a sink and satisfies $v \in s_{\Ha\Gamma^\prime}(e) \implies e \in S$ for all edges $e \in E^1(\Ha\Gamma^\prime)$. Clearly, this is not true for neither $v_1, v_2$ nor $w$ and therefore $v \in E^0(\Ha\Gamma^\prime) \setminus \{v_1, v_2, w\} \subset S$.
    \end{itemize}
    Evidently, $\Ha\Gamma_4$ is obtained from $\Ha\Gamma^\prime$ by deleting the set $S$. This concludes the proof.
\end{proof}

\subsection{Proof of Theorem \ref{main_res}} \label{refor_mainres_proof_subsec}

Finally, let us prove Theorem \ref{main_res}, the main result of this article.

\begin{proof}[Proof of Theorem \ref{main_res}]
    Let $\Ha\Delta := \mathrm{reduce}(\Ha\Gamma)$ be the reduced version of $\Ha\Gamma$ obtained by Algorithm \ref{pro_nuc_reduction_algo}. By Theorem \ref{pro_nuc_reduction_algo_proposition}, $\Ha\Delta$ is a normal hypergraph minor of $\Ha\Gamma$, and $C^\ast(\Ha\Gamma)$ is nuclear if, and only if, the same holds for $C^\ast(\Ha\Delta)$.

    Ad (1): From Theorem \ref{min_big_thm} we know that the minor operations preserve exactness of the associated hypergraph $C^\ast$-algebra. By Proposition \ref{main_res_forbidden_min_algs} the $C^\ast$-algebras $C^\ast(\Ha\Gamma_1), C^\ast(\Ha\Gamma_3)$ and $C^\ast(\Ha\Gamma_4)$ are not exact. Thus, $\Ha\Gamma_i \leq \Ha\Delta \leq \Ha\Gamma$ for some $i \leq 3$ implies that $C^\ast(\Ha\Gamma)$ is not exact.

    Ad (2): If $\Ha\Gamma_i \leq \Ha\Gamma$ holds for some $i \leq 3$, then the claim follows from (1) using that nuclearity implies exactness. Otherwise, by Proposition \ref{pro_nuc_HGamma4_minor} we may obtain $\Ha\Gamma_4$ from $\Ha\Delta$ using only two operations which preserve nuclearity of the associated hypergraph $C^\ast$-algebra. As $C^\ast(\Ha\Gamma_4)$ is not nuclear (see Proposition \ref{main_res_forbidden_min_algs}) it follows that $C^\ast(\Ha\Delta)$ is not nuclear. Then $C^\ast(\Ha\Gamma)$ is not nuclear as well.

    Ad (3): The hypergraph $\Ha\Delta$ satisfies the conditions for Theorem \ref{pro_nuc_reduction}. Thus, if $\Ha\Delta$ has none of the forbidden minors, then it must be an undirected hypergraph.
\end{proof}
\section{Examples}
\label{exa_sec}

In this section, we use the previous results to determine whether the $C^\ast$-algebra associated to a particular hypergraph is nuclear. In the first example, we retain the fact that every graph $C^\ast$-algebra associated to a finite graph is nuclear.

\begin{example}
    Assume that $\Gamma$ is a finite directed graph, i.e. $|r_{\Gamma}(e)| = |s_{\Gamma}(e)| = 1$. One verifies that the normalization procedure will only apply backward contraction on every edge $e$ with $s_{\Gamma}(e) \cap s_{\Gamma}(f) = \emptyset$ for all $f \in E^1(\Gamma) \setminus \{e\}$ and $s_{\Gamma}(e) \not \subset r_{\Gamma}(e)$. Clearly, the obtained normalized version of $\Gamma$ is again an ordinary graph which we may call $\Gamma^\prime$. Now, in $\Gamma^\prime$ every edge is easy in the sense of Definition \ref{pro_nuc_easy_edge_def}. Therefore, the reduction procedure will cut every edge of $\Gamma^\prime$. This leaves us with an undirected (hyper)graph $\Delta$ where every edge has exactly one vertex in its source. It is not hard to check that $C^\ast(\Delta) = \C^{n+m}$, where $n$ is the number of edges and $m$ the number of sinks in $\Delta$. In particular, $C^\ast(\Delta)$ is nuclear. By Theorem \ref{main_res} it follows that $C^\ast(\Gamma)$ is nuclear as well. 
\end{example}

\begin{example}
    Consider the hypergraph $\Ha\Gamma$ sketched below.
    
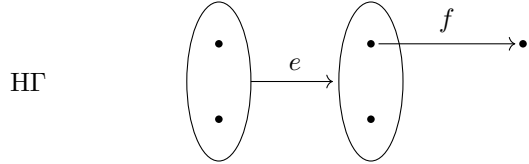
\begin{figure}[H]
        \centering
        \begin{tikzpicture}
            \node[label=0:$\Ha\Gamma$] at (-3, 0) {};
            \node [fill=black, circle, inner sep=1pt] at (0,.5) {};
            \node [fill=black, circle, inner sep=1pt] at (0,-.5) {};
            \draw[-] (0, 0) ellipse [x radius = 12pt, y radius=30pt];
            \node [fill=black, circle, inner sep=1pt] at (2,.5) {};
            \node [fill=black, circle, inner sep=1pt] at (2,-.5) {};
            \draw[-] (2, 0) ellipse [x radius = 12pt, y radius=30pt];
            \node [fill=black, circle, inner sep=1pt] at (4,.5) {};
            \draw[->] (.42, 0) -- (1.5, 0);
            \node[label=90:$e$] at (1, -.1) {};
            \draw[->] (2.1, .5) -- (3.9, .5);
            \node[label=90:$f$] at (3, .4) {};
        \end{tikzpicture}
    \end{figure} 
    \noindent The normalization procedure applies first range decomposition on the edge $e$ and then backward contraction on the edge $f$. We sketch these two steps below.
    
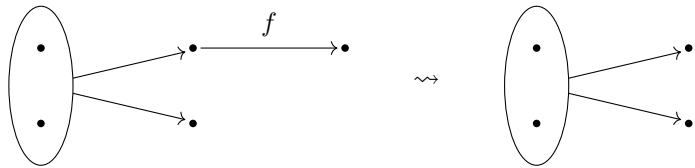
\begin{figure}[H]
        \centering
        \begin{tikzpicture}[baseline=0cm]
            \node [fill=black, circle, inner sep=1pt] at (0,.5) {};
            \node [fill=black, circle, inner sep=1pt] at (0,-.5) {};
            \draw[-] (0, 0) ellipse [x radius = 12pt, y radius=30pt];
            \node [fill=black, circle, inner sep=1pt] at (2,.5) {};
            \node [fill=black, circle, inner sep=1pt] at (2,-.5) {};
            \node [fill=black, circle, inner sep=1pt] at (4,.5) {};
            \draw[->] (.42, .1) -- (1.9, .45);
            \draw[->] (.42, -.1) -- (1.9, -.45);
            \draw[->] (2.1, .5) -- (3.9, .5);
            \node[label=90:$f$] at (3, .4) {};
        \end{tikzpicture}
        $ \qquad \leadsto \qquad$
        \begin{tikzpicture}[baseline=0cm]
            \node [fill=black, circle, inner sep=1pt] at (0,.5) {};
            \node [fill=black, circle, inner sep=1pt] at (0,-.5) {};
            \draw[-] (0, 0) ellipse [x radius = 12pt, y radius=30pt];
            \node [fill=black, circle, inner sep=1pt] at (2,.5) {};
            \node [fill=black, circle, inner sep=1pt] at (2,-.5) {};
            \draw[->] (.42, .1) -- (1.9, .45);
            \draw[->] (.42, -.1) -- (1.9, -.45);
        \end{tikzpicture}
    \end{figure} 
    \noindent Let us call the hypergraph on the right-hand side above $\Ha\Gamma^\prime$. One checks that $\Ha\Gamma^\prime$ is normal. The reduction procedure will cut both edges in $\Ha\Gamma^\prime$ since they end in a simple quasisink. Thus, one obtains the reduced hypergraph $\Ha\Delta$ sketched below.
    
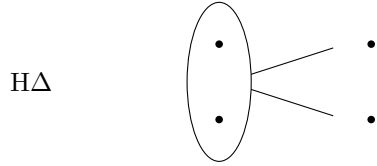
\begin{figure}[H]
        \centering
        \begin{tikzpicture}
            \node[label=0:$\Ha\Delta$] at (-3, 0) {};
            \node [fill=black, circle, inner sep=1pt] at (0,.5) {};
            \node [fill=black, circle, inner sep=1pt] at (0,-.5) {};
            \draw[-] (0, 0) ellipse [x radius = 12pt, y radius=30pt];
            \node [fill=black, circle, inner sep=1pt] at (2,.5) {};
            \node [fill=black, circle, inner sep=1pt] at (2,-.5) {};
            \draw[-] (.42, .1) -- (1.5, .45);
            \draw[-] (.42, -.1) -- (1.5, -.45);
        \end{tikzpicture}
    \end{figure} 
    \noindent One verifies $C^\ast(\Ha\Delta) = (\C^2 *_\C \C^2) \oplus \C^2$. Since, both factors in this direct sum are nuclear, $C^\ast(\Ha\Delta)$ is nuclear as well. By Theorem \ref{main_res} it follows that $C^\ast(\Ha\Gamma)$ is nuclear.
\end{example}

\begin{example}
    Let $\Ha\Gamma$ be the hypergraph sketched below.
    
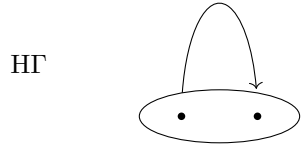
\begin{figure}[H]
        \centering
        \begin{tikzpicture}
            \node[label=0:$\Ha\Gamma$] at (-3, 1) {};
            \node [fill=black, circle, inner sep=1pt] at (-.5,.3) {};
            \node [fill=black, circle, inner sep=1pt] at (.5, .3) {};
            \draw[-] (0, .3) ellipse [x radius = 30pt, y radius=10pt];
            \draw[->] (0,.3) [partial ellipse=168:14:.5cm and 1.5cm];
        \end{tikzpicture}
    \end{figure}
    The normalization procedure will apply range decomposition on the only edge of this hypergraph, which leads to the hypergraph $\Ha\Delta$ below.
    
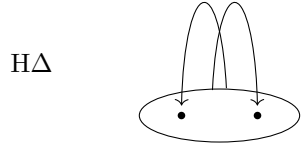
\begin{figure}[H]
        \centering
        \begin{tikzpicture}
            \node[label=0:$\Ha\Delta$] at (-3, 1) {};
            \node [fill=black, circle, inner sep=1pt] at (-.5,.3) {};
            \node [fill=black, circle, inner sep=1pt] at (.5, .3) {};
            \draw[-] (0, .3) ellipse [x radius = 30pt, y radius=10pt];
            \draw[->] (.2,.3) [partial ellipse=167:5:.3cm and 1.5cm];
            \draw[<-] (-.2,.3) [partial ellipse=175:14:.3cm and 1.5cm];
        \end{tikzpicture}
    \end{figure}
    Clearly, $\Ha\Delta$ is normal. Moreover, one checks that $\Ha\Delta$ contains no easy edge, no easy cycle and no edge ending in a simple quasisink. Thus, $\Ha\Delta$ is a reduced hypergraph. As it contains an edge with nonempty range we know from Theorem \ref{pro_nuc_reduction} that $\Ha\Delta$ has one of the forbidden minors $\Ha\Gamma_1, \dots, \Ha\Gamma_4$. What's more, we may investigate the proof of Theorem \ref{pro_nuc_reduction} and find that for $\Ha\Delta$ Case (A1) applies. In particular, $\Ha\Delta$ has the minor $\Ha\Gamma_3$. Indeed, one obtains $\Ha\Gamma_3$ from $\Ha\Delta$ by simply cutting one of the two edges. With Theorem \ref{main_res} it follows that $C^\ast(\Ha\Gamma)$ is not exact. 
    
    Indeed, $C^\ast(\Ha\Gamma) = \C^2 *_\C C(S^1) \supset C^\ast(\mathbb{F}_2)$ has been observed in \cite[Proposition 4.2]{triebHypergraphAlgebras2024} as first example of a non-nuclear hypergraph $C^\ast$-algebra.
\end{example}
\section{Problems}
\label{pro_sec}

In this section, we present three problems which came up during our research.

\begin{problem} \label{pro_bipartite_graph_algebra}
    The main theorem of this article, Theorem \ref{main_res}, gives a criterion that guarantees a hypergraph $C^\ast$-algebra to be non-nuclear. If this criterion is not met, then the problem of nuclearity of the given hypergraph $C^\ast$-algebra is reduced to the problem of nuclearity of a undirected hypergraph $C^\ast$-algebra. Thus, to completely characterize nuclearity of hypergraph $C^\ast$-algebras one needs to answer the following question: For which \emph{undirected} hypergraphs $\Ha\Delta$ is the $C^\ast$-algebra $C^\ast(\Ha\Delta)$ nuclear?
    

    The problem might be better phrased using the following definition of a $C^\ast$-algebra associated to a bipartite graph.
    Let $G = (U \sqcup V, E)$ be a bipartite graph, where $U$ and $V$ are the two sets of vertices and $E \subset \mathcal{P}(U \sqcup V)$ is the set of edges. For vertices $u \in U$ and $v \in V$ write $u \sim v$ if the edge $\{u, v\}$ exists in $E$. Associated to $G$ consider the $C^\ast$-algebra
    \begin{align*}
        C^\ast(G) = \C^U *_\C \C^V / (p_u p_v: u \not \sim v),
    \end{align*}
    where we write $p_u$ ($p_v$) for the projections generating $\C^U$ ($C^V$). If $\Ha\Delta$ is an undirected hypergraph, then $C^\ast(\Ha\Delta) = C^\ast(G)$ where $G = (E^0(\Ha\Delta) \sqcup E^1(\Ha\Delta), \sim)$ is a bipartite graph given by $v \sim e :\Leftrightarrow v \in s_{\Ha\Delta}(e)$ for all $v \in E^0(\Ha\Delta), e \in E^1(\Ha\Delta)$. As discussed in Remark \ref{mainres_remark_undirected_hyp_algs}, we believe that $C^\ast(G)$ is nuclear if, and only if, the following holds: For any vertex $u \in U \sqcup V$ there are at most two distinct vertices $v_1, v_2 \in U \sqcup V$ with $u \sim v_1$ and $u \sim v_2$.
\end{problem}

\begin{problem}
    A special case of a $C^\ast$-algebra associated to a bipartite graph as in Problem \ref{pro_bipartite_graph_algebra} is the unital free product $\C^2 *_\C \C^2$. For this $C^\ast$-algebra \cite{pedersenMeasureTheoryAlgebras1968} gives an explicit description as quotient of $C([0,1], M_2)$ (see also \cite{raeburnAlgebraGeneratedTwo1989} or the survey \cite{bottcherGentleGuideBasics2010}). Is there a similar description of the $C^\ast$-algebras associated to other bipartite graphs?
\end{problem}

\begin{problem}
    Compute the K-theory of a hypergraph $C^\ast$-algebra. For separated graph $C^\ast$-algebras \cite{araAlgebrasSeparatedGraphs2011} determines their K-theory using a free product description (similar to \cite{duncanCertainFreeProducts2010}) of the $C^\ast$-algebra in combination with a powerful six term exact sequence from \cite{thomsenKKtheoryEtheoryAmalgamated2003}. Although hypergraph $C^\ast$-algebras achieve certain free product constructions, they are generally not as easily described as free product $C^\ast$-algebra.
\end{problem}

\appendix

\printbibliography

\end{document}